\documentclass[10pt]{article}
\usepackage{hyperref}
\usepackage{color}
\usepackage[usenames,dvipsnames]{xcolor}

\hypersetup{
    colorlinks,
    linkcolor={RoyalBlue!55!black},
    citecolor={RoyalBlue!55!black},
    urlcolor={RoyalBlue!55!black}
}

\usepackage{bm}
\usepackage{enumitem}
\usepackage{graphicx}
\usepackage{tikz}
\usetikzlibrary{arrows.meta,calc}
\usepackage{varwidth}%
\usepackage{mathrsfs}
\usepackage{mathtools}
\usepackage[font=small,labelfont=bf]{caption}
\usepackage{subcaption}
\usepackage{wrapfig}
\usepackage{overpic}
\usepackage{multirow}

\usepackage[numbers,square,sort&compress]{natbib}

\usepackage{tcolorbox}
\usepackage[normalem]{ulem}
\usepackage{xfrac}
\usepackage[margin=1in]{geometry}
\usepackage{authblk}
\usepackage{amsmath,amsthm,amssymb,colonequals,etoolbox}
\usepackage{thmtools}
\usepackage{thm-restate}
\usepackage{url}
\usepackage{cleveref}

\usepackage{algorithm}
\usepackage{algpseudocode}
\usepackage{MnSymbol}
\usepackage{mdframed}
\usepackage{cancel}
\usepackage{booktabs}
\usepackage{tabularx}

\allowdisplaybreaks

\algrenewcommand\algorithmicrequire{\textbf{Require:}}
\algnewcommand\Input{\item[\textbf{Input:}]}
\algrenewcommand\algorithmicensure{\textbf{Output:}}

\renewcommand{\ge}{\geqslant}
\renewcommand{\geq}{\geqslant}
\renewcommand{\le}{\leqslant}
\renewcommand{\leq}{\leqslant}

\newtheorem{thm}{Theorem}[section]
\newtheorem{mydef}[thm]{Definition}
\newtheorem{myprop}[thm]{Proposition}
\newtheorem{mylemma}[thm]{Lemma}
\newtheorem{mycorollary}[thm]{Corollary}
\newtheorem{myfact}[thm]{Fact}

\newtheorem{mythm}[thm]{Theorem}
\newtheorem{assmp}[thm]{Assumption}

\newtheorem{rmk}[thm]{Remark}

\theoremstyle{definition}
\newtheorem{myex}[thm]{Exercise}

\newtheorem{problem}{Problem}
\crefname{problem}{problem}{problems}
\Crefname{problem}{Problem}{Problems}

\crefname{assmp}{assumption}{assumptions}
\Crefname{assmp}{Assumption}{Assumptions}

\mdfdefinestyle{theoremstyle}{%
linecolor=black,linewidth=1pt,%
frametitlerule=true,%
frametitlebackgroundcolor=gray!20,%
innertopmargin=0em,%
innerbottommargin=0.4em,%
nobreak=true,%
}

\DeclareMathOperator*{\argmin}{arg\,min}

\DeclareMathOperator*{\rank}{rank}

\DeclareMathOperator*{\diag}{diag}
\DeclareMathOperator*{\Tr}{\mathrm{tr}}

\newcommand{\calL}{\mathcal{L}}
\newcommand{\calD}{\mathcal{D}}

\newcommand{\R}{\ensuremath{\mathbb{R}}}
\newcommand{\C}{\ensuremath{\mathbb{C}}}

\newcommand{\N}{\ensuremath{\mathbb{N}}}

\newcommand{\norm}[1]{\lVert #1 \rVert}
\newcommand{\bignorm}[1]{\left\lVert #1 \right\rVert}

\newcommand{\ip}[2]{\ensuremath{\langle #1, #2 \rangle}}

\newcommand{\Var}{\mathrm{Var}}
\newcommand{\E}{\mathbb{E}}
\newcommand{\abs}[1]{\ensuremath{| #1 |}}
\newcommand{\bigabs}[1]{\ensuremath{\left| #1 \right|}}

\newcommand{\ceil}[1]{\lceil #1 \rceil}

\newcommand{\T}{\mathsf{T}}

\newcommand{\calG}{\mathcal{G}}

\newcommand{\calE}{\mathcal{E}}

\newcommand{\bbS}{\mathbb{S}}

\newcommand{\bmattwo}[4]{\begin{bmatrix} #1 & #2 \\ #3 & #4 \end{bmatrix}}

\newcommand{\cvectwo}[2]{\begin{bmatrix} #1 \\ #2 \end{bmatrix}}
\newcommand{\rvectwo}[2]{\begin{bmatrix} #1 & #2 \end{bmatrix}}

\numberwithin{equation}{section}

\newcommand{\calK}{\mathcal{K}}

\DeclarePairedDelimiterX{\infdivx}[2]{(}{)}{%
  #1\;\delimsize\|\;#2%
}

\newcommand{\opnorm}[1]{\norm{#1}_{\mathrm{op}}}

\newcommand{\e}{\varepsilon}

\newcommand{\rmd}{\mathrm{d}}

\newcommand{\sfN}{\mathsf{N}}

\newcommand{\wh}{\widehat}
\newcommand{\calS}{\mathcal{S}}

\DeclareMathOperator*{\eig}{eig}

\allowdisplaybreaks

\title{Learning Stabilizing Latent Feedback Controllers via \\ Multi-Step Backward Losses}
\author{
  Paul Lutkus\textsuperscript{1}\textsuperscript{\thanks{\hangindent=1.8em Corresponding Author. Email: \href{mailto:lutkus@usc.edu}{lutkus@usc.edu}\ \newline
  Project code available at: \href{https://github.com/paullutkus/Learning-Stabilizing-Latent-Feedback-Controllers}{github.com/paullutkus/Learning-Stabilizing-Latent-Feedback-Controllers}}}, 
  Lars Lindemann\textsuperscript{1,2}, 
  Stephen Tu\textsuperscript{3}\\
  \small\textsuperscript{1} Thomas Lord Department of Computer Science, University of Southern California \\
  \small\textsuperscript{2} Automatic Control Laboratory, ETH Z\"urich \\
  \small\textsuperscript{3} Ming-Hsieh Department of Electrical and Computer Engineering, University of Southern California
}
\date{\today}
\begin{document}

\maketitle

\begin{abstract}

Can learned latent dynamics models support the design of feedback controllers that provably stabilize the original system?
For linear time-invariant systems, we show that learning linear autoencoders and latent dynamics with multi-step \emph{prediction-reconstruction losses}---which we term \emph{backward losses} as they encourage latent trajectories to reconstruct full state trajectories---provably recovers the structure needed for stabilization.
Specifically, a sufficiently long reconstruction horizon ensures that a minimum ${H}_\infty$-norm controller designed using the learned latent model stabilizes the original system through output-feedback via the corresponding encoder, 
provided the latent dimension is large enough to capture all non-stable modes of the underlying system. 
We also show why both reconstruction and a multi-step horizon are essential: forward prediction losses that require only low prediction error along latent trajectories admit global loss minimizers that discard unstable modes regardless of horizon length, while single-step reconstruction losses can yield latent models whose stabilizing controllers fail to stabilize the original system.
Finally, experiments on a non-linear flexible-arm system illustrate the practical implications of our results for learning non-linear
autoencoders and latent dynamics.

\end{abstract}

\section{Introduction}
\label{sec:introduction}

Many control applications require stabilizing systems with high-dimensional, non-linear dynamics, including 
synchronizing large power grids after faults, 
suppressing aeroelastic flutter in flexible aircraft,
stabilizing a pneumatically actuated silicone arm under external loads, and 
maintaining stable humanoid locomotion on uneven or slippery terrain.
For such complex systems, where the complete state representation may be too high-dimensional to directly calculate an optimal state-feedback controller,
a natural approach is to utilize a low-dimensional surrogate model of the system, \emph{e.g.,} a \textit{latent dynamics model}, \textit{reduced-order model}, or \textit{latent world model}, which tracks the evolution of the original state $x$
under a (generally lossy) coordinate transformation $z=E(x)$,
where the map $E$ is commonly referred to as an \textit{encoder}.
In an ideal setting, such a latent surrogate model simultaneously satisfies
two key properties:
(i) it is rich enough to capture the key dynamical behavior of the system, and 
(ii) it is sufficiently low-dimensional to enable the computation of an optimal---or at least stabilizing---feedback controller $K(z)$.
Given properties (i) and (ii), an intuitive and computationally efficient 
feedback control strategy
for the original high-dimensional system is to 
compose the latent state feedback controller $K(z)$ with encoder $E(x)$, yielding
a form of output-feedback control.

The standard approach to designing such transformations from the full state to latent coordinates is to apply classical model reduction techniques, such as balanced truncation methods \cite{lall2002subspace, flinois2015projection}, proper orthogonal decompositions \cite{ravindran2000reduced,willcox2002balanced,flinois2015projection}, spectral submanifold theory \cite{haller2016nonlinear, alora2022data}, or partial feedback linearization \cite{spong1994partial,de2018symmetry}. 
In many cases, implementations of these techniques either strictly require, or are augmented by, knowledge of the true dynamics, an assumption that may not be practical in modern applications. 
More importantly, it is not always clear whether the application of these methods will yield \emph{stable zero dynamics}; without such a guarantee, the behavior of the system on the \emph{zero manifold} (i.e., kernel of the encoder $E$) could be 
unstable, and not controllable from only observing the latent coordinates.
This phenomenon is commonly encountered in partial feedback linearization, where the stability of zero dynamics, particularly for \textit{non-minimum phase systems}, is not immediately apparent when choosing an output map. In proper orthogonal decomposition, the reduced coordinates are chosen by applying singular value decomposition to trajectory data, which can potentially lead to unstable zero dynamics via the truncation of weakly excited unstable modes that happen to contribute low energy over a finite observation window. 
Finally, balanced truncation requires more sophisticated extensions to handle unstable systems~\cite{flinois2015projection}, as the infinite horizon Gramians become unbounded for unstable systems. 
While there exists techniques in the literature that address some of the
aforementioned limitations~\citep[see e.g.,][]{flinois2015projection, flinois2016feedback, cohen2024constructive, fatmawati2011balanced, buurmeijer2025taming},
we take a broader perspective in this work.
Namely, we ask whether it is possible
design a reduced-order modeling framework that is simultaneously 
(i) \emph{model-free}, requiring only trajectory snippets from the underlying
dynamics instead of a full description of the system dynamics,
(ii) employs a \emph{unified algorithm} that applies ``out of the box'' to
both linear and non-linear systems, and
(iii) {provably enables} \emph{stable feedback control design} 
using the output-feedback design strategy described previously.

In this paper, we propose such a modeling framework by building on 
algorithms from the
machine learning (ML) and reinforcement learning (RL) communities
for learning latent dynamics models from data,
where such methods fall under the general umbrella of
\emph{world modeling}~\cite{watter2015embed,hafner2019learning,hafner2019dream,zhou2025dinowm,balestriero2025lejepa,hansen2024tdmpc2,maes2026leworldmodel}. 
The advantages of utilizing an ML/RL world modeling approach 
are evident, as they immediately satisfy the first two of our 
desiderata: (i) they operate directly on trajectory data (model-free),
and (ii) they provide general algorithms that make little assumptions 
on the underlying dynamics (unified algorithms).
However, the last desiderata---enabling stable feedback control design---remains an open question.
Indeed, while the literature demonstrates impressive empirical
results with planning algorithms that utilize learned world models in 
the loop~\cite[see e.g.,][]{hafner2019learning,hansen2024tdmpc2,assran2025vjepa2,maes2026leworldmodel}, the theoretical foundations of such design
procedures, including basic questions regarding 
closed-loop stability, safety, and performance, remain largely unaddressed
in the literature. 
Our work initiates a study of these foundations, 
focusing on the following fundamental question:
\begin{center}
\emph{Can latent dynamics models and encoders 
learned with ML/RL approaches
be used to design feedback controllers that 
are guaranteed to be closed-loop stable on the true underlying system?} 
\end{center}

To focus on the central challenges of model reduction, we assume access to the full system state. Although partial observability is important in practice, this assumption allows us to isolate the effects of model reduction from the separate complications introduced by filtering.

\subsection{A Motivating Example}
\label{sec:introduction:example}

\begin{figure}[t]
    \centering
    \includegraphics[width=\linewidth]{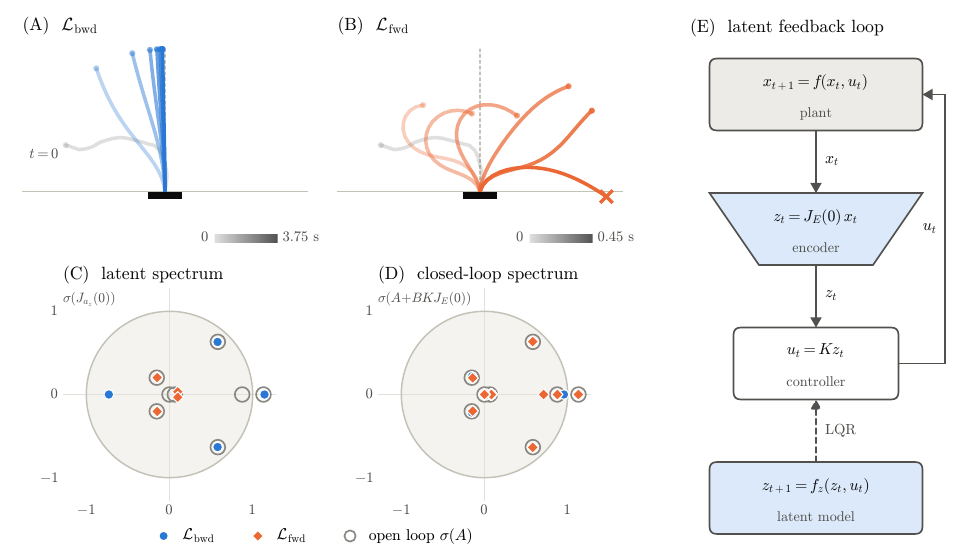}
    \caption{Stabilization of a flexible arm using linear latent output feedback. Only $\calL_\mathrm{bwd}$, 
    the multi-step prediction-reconstruction loss, leads to a stabilizing latent output-feedback controller.
    In (A), (B), the light gray arm denotes the shared initial state, and opacity denotes the position 
    of the arm in time (darker denotes later). In (B), the $\times$ denotes the arm impacting the ground.
    Panels (C) and (D) are linearizations about the upright equilibrium.
    }
    \label{fig:arm_model}
\end{figure}

Before we describe our contributions in detail, we first motivate them with a numerical example inspired by soft robotics applications.
We consider a multi-link flexible arm model, which we instantiate with $20$ sections to obtain an $n=40$ dimensional dynamical system; after
time-discretization we obtain a non-linear discrete-time system
$x_{t+1} = f(x_t, u_t)$.
The control task is to stabilize the arm around its unstable upright 
equilibrium; see \Cref{fig:arm_model}, panels (A) and (B).
For illustrative purposes, we set system parameters such that the number of unstable modes
of the 
Jacobian $\partial_x f(x, u)$ 
evaluated at the upright equilibrium 
is equal to one; 
physically, this corresponds to an arm that can lose balance through one independent bending motion, while perturbations in its remaining deformation modes decay.
This suggests that, at least near the equilibrium, there exists an
encoder $z=E(x)$ and latent dynamics $f_z(z, u)$ operating in an $m$
dimensional space, with $m \ll n$, which still captures the relevant dynamics of $f(x, u)$ needed
for stabilization.

With this setup in place, we explore two predominant learning strategies in the ML/RL literature on latent dynamics modeling for learning
the encoder $E$, latent dynamics $f_z$, and optional decoder $D$:
\begin{align}
    \calL_{\mathrm{bwd}}(E, D, f_z) &:= \text{(Latent Prediction-Reconstruction) + (Latent Prediction) + (SIGReg)}, \label{eq:loss_combo_bwd}\\
    \calL_{\mathrm{fwd}}(E, f_z) &:= \text{(Latent Prediction) + (SIGReg)}, \label{eq:loss_combo_sigreg}
\end{align}
where the different loss functions above are commonly used in the literature:
\begin{enumerate}[label=(\alph*)]
    \item \emph{Latent Prediction-Reconstruction (Backward Conjugacy):} $\sum_{t=0}^T \|D(f_z^t(E(x_0), \bm{u}_{0:t-1}))- x_t \|^2$,
    \item \emph{Latent Prediction (Forward Conjugacy):} $\sum_{t=1}^T\|f_z^t(E(x_0), \bm{u}_{0:t-1})-E(x_t)\|^2$,
    \item \emph{SIGReg:} $\E_{a}\left[\mathcal{T}\left(\{ \ip{a}{E(x_i)} \}_{i=1}^{N}\right)\right]$, where $a$ is drawn uniformly from the unit sphere $\bbS^{m-1}$ in $\R^m$, $\mathcal{T}$ is the Epps-Pulley statistic for testing agreement with $\sfN(0,1)$~\cite{balestriero2025lejepa},
    and $\{x_i\}_{i=1}^{N}$ denotes
    $N$ mini-batch samples of trajectory states.
\end{enumerate}
For loss functions (a) and (b) above,
the notation $f^t_z(z_0, \bm{u}_{0:t-1})$ denotes the $t$-step flow map
of the latent dynamics starting from $z_0 = E(x_0)$ and applying control inputs $\bm{u}_{0:t-1} = (u_0, \dots, u_{t-1})$; for both losses we operate them in the \emph{multi-step} regime, where  $T > 1$.
The learning strategies in \eqref{eq:loss_combo_bwd}
and \eqref{eq:loss_combo_sigreg} represent
two broad approaches.
The backward reconstruction loss \eqref{eq:loss_combo_bwd} jointly learns a decoder $D$
alongside $(E, f_z)$, which is used to reconstruct latent predictions back
into the original space, mirroring approaches such 
as E2C~\cite{watter2015embed}, 
deep Koopman embeddings~\cite{lusch2018deep},
PlaNet~\cite{hafner2019learning}, and Dreamer~\cite{hafner2019dream,hafner2025mastering}.
The forward prediction loss \eqref{eq:loss_combo_sigreg} on the other hand follows recent \emph{decoder-free} approaches
such as DINO-WM~\cite{zhou2025dinowm}, LeJEPA~\cite{balestriero2025lejepa},
and TD-MPC2~\cite{hansen2024tdmpc2},
which measures the discrepancy between the latent
model trajectories and the encoded
system trajectories.
In order to prevent \emph{representation collapse}, where
trivial models and encoders are learned, forward losses
must be paired with a regularizer; following LeJEPA~\cite{balestriero2025lejepa} and LeWorldModel~\cite{maes2026leworldmodel},
we utilize SIGReg, a commonly used approach
which encourages the encoded latents $E(x_i)$ to match
an isotropic Gaussian (and is therefore non-degenerate).

Given a learned $(E, f_z)$,
we then employ a simple control design strategy visualized in \Cref{fig:arm_model}(E).
We choose the latent dynamics to be control-affine, i.e., $z_{t+1}= f_z(z_t, u_t) = a_z(z_t)+g_z(z_t)u_t$,
where $a_z$, $g_z$, the encoder $E$, and the decoder $D$ are all MLPs.
Let $J_q(\bar{x}) := \frac{\rmd}{\rmd x} q(x)|_{x=\bar{x}} $ denote the Jacobian of a function $q(x)$ at $\bar{x}$.
We compute an LQR gain $K$
from the linearization $(J_{a_z}(0), g_z(0))$ about the equilibrium point,
and construct the corresponding output-feedback controller
$\pi(x)=KJ_E(0)x$ using the linearization of the MLP encoder as the output function.
We train two models, using \eqref{eq:loss_combo_bwd} and \eqref{eq:loss_combo_sigreg} respectively, 
from $2000$ trajectories of length $T=40$ steps (corresponding to $2$\,s rollouts) from the non-linear arm.
Our evaluation metric is whether or not the controller $\pi(x)$
stabilizes the original non-linear system $f(x, u)$. 

The results of this experiment are displayed in \Cref{fig:arm_model},
which illustrates drastically different closed-loop behavior between
the models learned from \Cref{eq:loss_combo_bwd} and \Cref{eq:loss_combo_sigreg}:
released from rest with sections bent up to $122^\circ$ from upright (tip $70$\,cm to the side), 
the arm returns upright when controlled by the latent model learned via 
$\calL_\mathrm{bwd}$, whereas under $\calL_\mathrm{fwd}$'s model it reaches the ground within $0.46$\,s.
Spectral analysis of the eigenvalues of the Jacobian $J_{a_z}(0)$
evaluated at the equilibrium provides some intuition for the observed
behavior. \Cref{fig:arm_model} (panels (C) and (D)) shows that the Jacobian $J_{a_z}(0)$ for the latent model learned via loss \eqref{eq:loss_combo_bwd} 
(shown in blue) includes an unstable eigenvalue nearly identical to the unstable eigenvalue of the true system, 
whereas the spectrum of $J_{a_z}(0)$ learned by the decoder-free forward loss \eqref{eq:loss_combo_sigreg} includes only stable eigenvalues.
Informally, this implies that the zero dynamics learned from \eqref{eq:loss_combo_sigreg} must contain an unstable eigenvalue, which is
inaccessible to the latent controller.

This simple experiment motivates the following hypothesis: when a multi-step prediction-reconstruction loss is utilized, and the latent model has sufficient capacity, the latent dynamics must converge to a representation that captures the unstable subsystem, enabling the design of a stabilizing output-feedback controller from the latent dynamics model. 
The remainder of this paper is dedicated towards analyzing this
phenomenon, using linear dynamical systems as a tractable setting for studying its underlying mechanisms.

\subsection{Contributions} 
\label{sec:introduction:contributions}

With the motivating example from \Cref{sec:introduction:example} in place,
we now state our main contributions. 
While the losses considered in \Cref{sec:introduction:example}
are written for non-linear systems,
our analysis focuses on the specialized setting where 
the true system is a stabilizable
LTI system $f(x, u) = A x + Bu$.
We utilize loss $\calL_{\mathrm{bwd}}$ in \eqref{eq:loss_combo_bwd}
to learn a linear autoencoder pair $(E, D)$ alongside a 
linear latent dynamics model $f_z(z, u) = Lz + Hu$; thus our learnable
parameters are the tuple $\Theta := (E, D, L, H)$.
We show the following results (\Cref{sec:results}),
whenever the latent dimension $m$ has sufficient capacity to capture the non-stable modes of $A$:
\begin{enumerate}[label=(C\arabic*)]
    \item \textbf{Latent stabilization problem formulation, and exact solution with \emph{known} dynamics $(A, B)$:}
    We provide a concrete formulation of
    the latent stabilization problem 
    (\Cref{prob:latent_stabilization}),
    and provide an exact solution in \Cref{sec:problem_formulation:invariance_subspaces}
    when the LTI dynamics
    $(A, B)$ are \emph{known}.
    Two solution concepts emerge that 
    form the cornerstone
    of our subsequent analysis
    for model-free algorithms:
    (i) \emph{detectability of $(A, E)$}, 
    requiring that every mode of $A$ that is unobservable through $E$ should be asymptotically stable;
    and (ii) \emph{encoder invariance}, requiring that the kernel of $E$ should be $A$-invariant.
    Indeed, we show that both conditions (i) and (ii) ensure that the zero dynamics defined on the kernel of $E$ are stable.

    \item \textbf{Optimizers of long horizon multi-step backward losses satisfy \emph{approximate} detectability and invariance conditions:}
    We next show (\Cref{thm:informal_recovery}) that by taking the multi-step horizon $T$
    for the latent prediction-reconstruction loss in $\calL_{\mathrm{bwd}}$ (cf.~\Cref{eq:loss_combo_bwd}) \emph{sufficiently long} (depending on system-theoretic quantities involving only $(A, B)$ which we make explicit in our analysis), then
    any global minimizer $\Theta$ of $\calL_{\mathrm{bwd}}$ 
    satisfies an \emph{approximate} 
    version of the detectability and encoder invariance and conditions described in (C1). 
    We then use these approximate conditions as the structural requirements for a robust analysis of closed-loop latent control design 
    which we describe next.

    \item \textbf{A robust analysis of closed-loop latent control design, and an end-to-end latent stabilization guarantee:}
    We utilize the approximate detectability and invariance conditions described in (C2) and
    show (\Cref{thm:resolvent_gain_bound}) that under these conditions, there exists a control design procedure
    $K = K(L, E)$ which yields a stable closed-loop
    system $A + BKE$ through the output-feedback
    controller $u(x) = KEx$.
        Combining \Cref{thm:informal_recovery} with
    \Cref{thm:resolvent_gain_bound}
    then yields our end-to-end result
    (\Cref{cor:latent_stability_transfer}) that minimizing
    the multi-step $\calL_{\mathrm{bwd}}$ loss solves the latent stabilization problem (\Cref{prob:latent_stabilization}).
    We prove our results for $K(L, E)$ that minimizes the $H_\infty$-norm gain
    of the latent closed-loop, although empirically (\Cref{sec:experiments})
    we also find similar results when utilizing a latent
    LQR control design.

    \item \textbf{Insufficiency of multi-step forward losses and single-step backward losses for the latent stabilization problem:}
    Finally, we show that the multi-step backward loss $\calL_{\mathrm{bwd}}$ is actually
    \emph{necessary} in the following sense:
    (i) for the $\calL_{\mathrm{fwd}}$ loss in \eqref{eq:loss_combo_sigreg}, regardless of how long the multi-step horizon $T$ is set to, there will always exist global minimizers 
    where the encoder $E$ results in that $A+BKE$ is \emph{always unstable}
    for any $K$ (\Cref{thm:fwd_cycle_insufficient});
    (ii) for the $\calL_{\mathrm{bwd}}$ loss, 
    we construct specific problem instances where when $T=1$,
    but the set of gains $K$ that stabilize $L+HK$ is \emph{completely
    disjoint} from the set of gains $K$ that stabilize $A+BKE$ (\Cref{thm:bwd_single_step_nonnormal}).

\end{enumerate}

Furthermore, while our theory applies to LTI systems, we provide
a more detailed numerical study in \Cref{sec:experiments}
for the high-dimensional non-linear system from \Cref{sec:introduction:example}, where we further illustrate the benefits of multi-step backward reconstruction losses
when learning latent dynamics for non-linear control design.
We conclude in \Cref{sec:conclusion} with some directions for future work inspired from our results.

\section{Related Work}
\label{sec:related_work}

\paragraph{Classical model reduction.}
Reducing an representation of a dynamical system while maintaining input/output behavior is a classical problem in control theory.
Balanced truncation (BT) starts from a state-space model $(A,B,C)$,
and uses the system's Hankel operator to select coordinates that preserve the system's input--output behavior up to a tolerance that depends on the sum of the unmodelled Hankel singular values (HSVs)~\cite{moore1981principal}.
While the most classical version of balanced truncation assumes knowledge of the 
state-space representation, 
this can be relaxed to assuming access to state trajectories from
forward and adjoint impulse responses~\cite{rowley2005model},
or 
from output impulse-response data using the 
eigensystem realization algorithm (ERA)~\cite{juang1985eigensystem,ma2011reduced}.
For unstable systems, \citet{flinois2015projection} analyze the direct
application of unmodified snapshot-based BT using finite-horizon
Gramians, and \citet{flinois2016feedback} study the corresponding
ERA-based approach.
The key intuition is that controllable and observable unstable dynamics
eventually dominate the finite-horizon Hankel singular spectrum, and hence BT in principle will 
capture the unstable modes in its reduced order model
(ROM), assuming both the specified order is large enough and the trajectory data is long enough. This recovery mechanism, at a high level, is similar to our analysis strategy. However there are 
several key differences.
First, the convergence guarantees of \citet{flinois2015projection}
are asymptotic as the snapshot horizon tends to infinity, whereas
our work provides explicit upper bounds on the horizon length
required for multi-step backward reconstruction.
Second, while non-linear extensions for BT  exist~\cite{lall2002subspace,scherpen1993balancing,kramer2023nonlinear},
the reconstruction losses studied here provide a different formulation
that can be \emph{directly applied} to non-linear dynamics and non-linear
coordinate reductions without formulating
new algorithms.
Finally, classical BT error guarantees~\cite{enns1984model,glover1984all} bound the input--output
discrepancy between the original system and the ROM.
A controller that stabilizes the ROM need not stabilize the original
system solely by virtue of this approximation guarantee;
we provide a specific example in \Cref{sec:appendix:balanced_truncation}.

Closely related to BT is proper orthogonal decomposition (POD), which selects a low-dimensional linear subspace that minimizes the $\ell_2$-projection error of a collection of state snapshots in a prescribed inner product~\cite{sirovich1987turbulence,berkooz1993proper,rowley2005model};
for the Euclidean inner product, this is equivalent to PCA on the data matrix of state snapshots.
POD therefore provides a natural connection between classical model reduction and reconstruction-based representation learning: its coordinates optimally reconstruct the observed snapshots among orthogonal projections of a fixed rank.
However, the optimality of the state reconstruction by itself does not necessarily guarantee preservation of
dynamical behavior.
Balanced proper orthogonal decomposition (balanced POD) addresses this limitation by using the observability Gramian in the inner product, 
which relates balanced POD to
BT~\cite{rowley2005model}.
Our analysis, on the other hand, 
studies an alternate route to encoding dynamical
information into reconstruction, namely through
multi-step reconstruction, rather than
a specific inner-product choice.

\paragraph{Geometric representations and non-linear stabilization.}
Spectral sub-manifold (SSM) theory characterizes invariant manifolds that provide non-linear continuations of spectral subspaces near an equilibrium~\cite{haller2016nonlinear}. 
For model reduction, slow SSMs are particularly useful when the discarded modes decay faster than the retained modes~\cite{alora2022data,buurmeijer2025taming}, but establishing these conditions can be challenging when the dynamics are known only through data. Most importantly, extending SSM reductions to controlled system is non-trivial and requires additional structure on the control signal~\cite{mahlknecht2022using,kaundinya2025data,haller2024nonlinear}, since general control inputs need not preserve the invariance of a SSM of the unforced dynamics.

Koopman operator methods represent non-linear dynamics through the linear evolution of observable functions of the state~\cite{brunton2022modern}. 
Practical implementations learn finite-dimensional approximations, which need not represent the dynamics exactly~\cite{williams2015data,brunton2016koopman}; consequently, stabilizing the learned Koopman approximation need not stabilize the original non-linear system. 
Existing results address this gap under additional conditions on the observables, approximation error, and controller synthesis. 
\citet{mamakoukas2023learning} construct Lyapunov certificates from stable Koopman approximations when the residual error satisfies a proportional bound sufficiently small relative to the Lyapunov stability margin; however, their learning procedure does not enforce this bound uniformly over a region. 
SafEDMD~\cite{straesser2026safedmd} incorporates proportional error bounds into robust controller synthesis, providing probabilistic guarantees of exponential stability within a certified region of attraction. 
However, feasibility of 
relevant controller-synthesis LMIs
is assumed rather than established from stabilizability of the original non-linear system.
Thus, while these approaches provide sufficient conditions for non-linear stabilization, neither small empirical prediction errors nor stability of a learned Koopman approximation alone ensures that the requisite conditions hold.

The templates-and-anchors framework relates a low-dimensional model of target behavior to a more detailed model of the physical system~\cite{full1999templates,holmes2006dynamics}.
In its dynamical systems formulation, anchoring requires both reproducing the template dynamics on an invariant manifold, and ensuring that nearby motions approach that manifold; the framework however does not prescribe a general method for constructing such a realization.

Partial feedback linearization (PFL) provides a complementary approach for underactuated mechanical systems~\cite{spong1994partial,khalil2002nonlinear}.
While the resulting decomposition permits local asymptotic stability 
to be established from stable linear output dynamics and locally asymptotically stable zero dynamics,
PFL requires both exact knowledge of the model and a careful selection of the 
output map. In general the zero dynamics
are not guaranteed a-priori to be stable, and can be sensitive to model-mismatch error.

\paragraph{Bisimulation.}
Bisimulation formalizes relationships between a full system and an abstraction.
For control systems, exact bisimulation characterizes state-space reduction~\cite{vanderschaft2004equivalence},
and approximate bisimulation bounds discrepancies between observed trajectories~\cite{girard2007approximate}.
In reinforcement learning, bisimulation connects state aggregation to value approximation~\cite{ferns2004metrics}, reward and transition prediction~\cite{gelada2019deepmdp}, and invariance to task-irrelevant visual information~\cite{zhang2021invariant}.
These approaches preserve behavior relative to specified outputs or rewards, but choosing these quantities
to capture the information needed for e.g., stabilization is itself non-trivial; in particular, exact output equivalence can conceal unstable internal dynamics, while a non-zero approximation permits errors that persist over time.
Furthermore, value approximation guarantees alone do not necessarily establish Lyapunov stability or asymptotic convergence of the closed-loop state to a target equilibrium.

\paragraph{Latent world models.}
World models~\cite{ha2018recurrent} support planning and policy learning by predicting dynamics in a latent space.
As discussed in \Cref{sec:introduction}, 
there are two broad types of latent world models.
The first class of models decode latent predictions to match trajectories in the original space, including
E2C~\cite{watter2015embed},
SOLAR~\cite{zhang2019solar}, 
PlaNet~\cite{hafner2019learning}, and
Dreamer~\cite{hafner2019dream,hafner2025mastering}.
On the other hand, \emph{decoder-free}
approaches build latent representations and models without explicitly reconstructing observations,
including DINO-WM~\cite{zhou2025dinowm},
LeJEPA~\cite{balestriero2025lejepa},
LeWorldModel~\cite{maes2026leworldmodel}, and 
TD-MPC2~\cite{hansen2024tdmpc2}.
Such approaches must rely on careful regularization
procedures to avoid representation collapse,
including SIGReg~\cite{balestriero2025lejepa}
and earlier covariance regularization approaches in self-supervised learning~\cite{bardes2022vicreg,zbontar2021barlow}.
The goal of our work is to provide an analysis 
of a representative approach from both classes of world models 
applied to LTI systems,
illustrating that multi-step
reconstruction losses actually play a critical part
for enabling stabilization from latent control design.

\paragraph{Latent analysis and control design.}
Latent analysis and control design with ML approaches using autoencoders and world models
have become a increasingly
popular research direction.
For analysis of an existing controller, MORALS~\cite{vieira2024morals} learns an autoencoder and latent closed-loop dynamics from trajectories, then constructs Morse graphs in the latent space to estimate attractors and their regions of attraction. 
\citet{nakamura2025generalizing} apply Hamilton--Jacobi reachability analysis in the latent space of a world model to construct safety filters, demonstrating their effectiveness in practice;
however, provably transferring such latent certificates
to the real system is left as an open question.
Our prior work \cite{lutkus2025latent} gives approximate conjugacy and geometric conditions under which latent Lyapunov and barrier certificates imply stability and invariance guarantees for the original non-linear dynamics.
In the case of stabilization, these results only show
convergence to the \emph{kernel} of the encoder, and leave open the
question of whether or not the zero dynamics are actually stable.
The non-linear simulations 
however suggest that stable closed-loop behavior is learned, which provided the initial motivations of our study.
\citet{mootz2026latent} show that an exact latent model can conceal safety-relevant directions, and derive barrier-transfer bounds with separate margins for model error and variation of safety across states sharing the same encoding.
Their analysis relates safety transfer to detectability of the safety-relevant dynamics;
as we will shortly see, detectability through the encoder
also plays a key role in our analysis.
\citet{toso2025learning} 
learn to stabilize LTI systems with policy gradient methods by first learning a basis for the
left unstable subspace of $A$, followed
by a sequence of discounted LQR policy gradient optimizations that operate on the dynamics
restricted to the unstable subspace.
Our analysis also crucially lies on showing
the encoder recovers the left unstable subspace of $A$, but unlike \cite{toso2025learning} we do not require rollouts from the adjoint system.
Most related to this paper is the concurrent work of \citet{toso2026preserving}, 
who show the insufficiency
of JEPA world models for recovering unstable dynamics,
and how 
action-reconstruction (i.e., inverse dynamics prediction) 
recovers a detectable encoder;
we defer a more detailed
comparison to \Cref{sec:appendix:comparison_to_toso}.

\section{Problem Setup and Exact Solution with Known Dynamics}
\label{sec:problem_formulation}

We consider discrete-time, linear time-invariant (LTI) systems of the form:
\begin{align}
    x_{t+1}=Ax_t+Bu_t, \qquad x_t \in \R^n, \,\, u_t \in \R^p. \label{eq:LDS}
\end{align}
In this work, we restrict our attention to linear dynamics $(A, B)$ that are is stabilizable.
\begin{assmp}[Stabilizability]
\label{assmp:stabilizability}
We assume the pair $(A, B)$ in \eqref{eq:LDS} is stabilizable.
\end{assmp}
Next, we assume that we do not have knowledge of the dynamics matrices $(A, B)$, but instead have access to a dataset of $N$ length-$T$ trajectories of state-control pairs generated under the following data model:
\begin{assmp}[Data model]
    \label{assmp:data_model}
    We assume access to a dataset $\calD_N := \{ (x_0^{(i)}, u_0^{(i)}, \dots, x_{T-1}^{(i)}, u_{T-1}^{(i)}, x_T^{(i)}) \}_{i=1}^{N}$ generated from the LTI system $(A, B)$ in \eqref{eq:LDS} via independent initial conditions and control inputs: %
    \begin{align}
        x_{t}^{(i)} = A^tx_0^{(i)} + \sum_{k=0}^{t-1}A^{t-1-k}Bu_k^{(i)}, \qquad x_0^{(i)}\sim\sfN(0, \sigma_x^2 \cdot I_n),\ u_t^{(i)}\sim\sfN(0, \sigma_u^2 \alpha_T(t) \cdot I_p), \label{eq:data_model}
    \end{align}
    where
    both $x_0^{(i)}$ and $\bm{u}_{0:T-1}^{(i)} := ( u_0^{(i)}, \dots, u_{T-1}^{(i)} )$ are sampled mutually independent
    across time, 
    $\sigma_x^2,\sigma_u^2 > 0$ are positive
    variance parameters,
    and
    $\{ \alpha_T(k) \}_{k=0}^{T-1}$ is a sequence of positive scaling factors that modulate the energy input into the system.
\end{assmp}

The scaling factors $\alpha_T(k)$ are introduced in the control input distribution \eqref{eq:data_model} in order to fix the total amount of energy that is propagated through the system 
\eqref{eq:LDS} over the horizon $T$.
Indeed, we have
that 
$\calE(\{ \alpha_T(k) \}_{k=0}^{T-1}) := \sum_{k=0}^{T-1} \E\norm{u_k^{(i)}}^2 = \sigma_u^2 \sum_{k=0}^{T-1} \alpha_T(k)$.
A common choice of $\alpha_T(k)$ we will consider is to set
$\alpha_T(k) \equiv 1/T$ so that the total input energy equals $\sigma_u^2$,
although other choices such as decaying rates 
$\alpha_T(k) \asymp 1/(k+1)^{\eta}$ are certainly possible as well (and will also be covered by our subsequent analysis).

With the LTI system \eqref{eq:LDS}
and data model \eqref{eq:data_model} in place, we now state the main problem we study.
\begin{tcolorbox}[colback=gray!10, colframe=gray!50, boxrule=0.5pt,
    sharp corners, colbacktitle=gray!20,
    coltitle=black, fonttitle=\bfseries]
\begin{problem}[Latent stabilization]\label{prob:latent_stabilization}
Suppose that the LTI system $(A, B)$
in \eqref{eq:LDS} 
is stabilizable (\Cref{assmp:stabilizability}),
and we have access to a dataset $\calD_N$ of trajectories
(\Cref{assmp:data_model}).
Choose a latent dimension $m \in [n]$.
\emph{When possible}, learn a 
\begin{enumerate}[label=(\roman*), nosep, topsep=4pt]
    \item \emph{latent linear encoder} 
$z = E x$ with
$E \in \R^{m \times n}$ and
    \item \emph{latent linear dynamics} $z_{t+1} = L z_t + H u_t$ with $L \in \R^{m \times m}$
and $H \in \R^{m \times p}$,
\end{enumerate}
such that there exists a \emph{latent control synthesis} procedure $K = K(L, H) \in \R^{p \times m}$ that ensures
that both $L + HK$ and $A + BKE$ are Schur stable.
\end{problem}
\end{tcolorbox}

Here, we utilize the standard terminology 
\emph{Schur stable} to refer to a square matrix $M$ with spectral radius $\rho(M)$ strictly less than one, i.e., $\rho(M) < 1$.
\Cref{prob:latent_stabilization} makes no mention of any corresponding decoder model; this is intentional, as the decoder is an algorithmic detail of our proposed solution,
but not required to state the problem formulation.
Furthermore, \Cref{prob:latent_stabilization} involves the
qualifier phrase \emph{when possible}, implying that our desired goal may not always be achievable.
Indeed, the choice of the latent dimension $m$ plays an important role in determining when such a problem
is feasible, as we will see shortly (cf.~\Cref{prop:latent_stability_iff_output_feedback_stability}).
Before we start discussing our solution, we first discuss how we can solve the latent stabilization problem given \emph{exact knowledge} of the dynamics $(A, B)$, as this motivates the key concepts for what follows.

\paragraph{Notation.}
We collect the notation used in this paper, which is standard.
Let $\mathbb{D}:=\{z\in\C\mid\abs{z}<1\}$ denote the open unit disk, and let $\mathbb{D}^c:=\C\setminus\mathbb{D}=\{z\in\C\mid\abs{z}\geq1\}$ denote its complement, including the unit circle.
For a vector $x \in \C^n$, we let $\norm{x}_p$ denote its $\ell_p$-norm, with $\norm{x} = \norm{x}_2$.
For a matrix $M \in \C^{n \times m}$, the notation
$\opnorm{M}$ denotes its operator norm,
$\norm{M}_F$ its Frobenius norm, and
$\sigma_{\min}(M)$ denotes its $\min\{n, m\}$-th largest
singular value.
Furthermore, $\mathrm{col}(M)$,
$\mathrm{row}(M)$, and $\mathrm{kern}(M)$
denote the column, row, and nullspace of $M$, respectively.
When $M$ is square, $\rho(M)$ denotes its spectral radius,
$\sigma(M)$ denotes its eigenvalues
and $\Tr(M)$ its trace.
For $n, m \in \N_+$, denote 
$\bbS^{n-1} := \{ x \in \R^n \mid \norm{x} = 1 \}$,
$\mathrm{GL}(n) := \{ M \in \R^{n \times n} \mid \det(M) \neq 0 \}$,
$B_{n \times m}^{(\mathrm{op})}(\gamma) := \{ M \in \R^{n \times m} \mid \opnorm{M} \leq \gamma \}$, and 
$O(n, m) := \{ M \in \R^{n \times m} \mid M^\T M = I_m \}$ as the set of $n \times m$ real-valued matrices with
orthonormal columns; the latter requires that $n \geq m$.
For a subspace $M$, the matrix $P_M$ (resp.\ $P^\perp_M$) denotes the orthogonal projector onto $M$ (resp.\ onto the orthogonal complement of $M$); when $M$ is a matrix, we let
$P_M = P_{\mathrm{col}(M)}$.

\subsection{Latent Stabilization with Known Dynamics}
\label{sec:problem_formulation:invariance_subspaces}

Fix a latent encoder $E \in \R^{m \times n}$
with $m \in [n]$.
Given a state-input trajectory $(\bm{x}_{0:T}, \bm{u}_{0:T-1})$
from \eqref{eq:LDS}, 
the encoder induces a corresponding latent trajectory
$(\bm{z}_{0:T}, \bm{u}_{0:T-1})$,
with $z_t := E x_t$.
A natural question which arises is
the following:
what are conditions on $E$ which ensure that
the latent trajectory is \emph{memoryless}, i.e.,
there exists a latent dynamics model $f_z$ such that
$z_{t+1} = f_z(z_t, u_t)$ for all $t \in \{0, \dots, T-1\}$?
For $M \in \R^{n \times n}$ and $W \subseteq \R^n$ a subspace, recall that $W$ is $M$-invariant if $MW \subseteq W$.
Invariant subspaces provide an exact characterization of when a memoryless representation of the latent trajectory is possible, as described in the following proposition.

\begin{restatable}[Memoryless latent trajectories $\iff$ $\mathrm{kern}(E)$ is $A$-invariant]{myprop}{invarianceequivlatentdynamics}
\label{prop:invariance_equiv_latent_dynamics}
Fix an encoder $E \in \R^{m \times n}$.
There exists a function $f_z : \mathrm{col}(E) \times \R^p \mapsto \R^m$ such that for any 
$\tau \in \N_+$ and state-input trajectory
$(\bm{x}_{0:\tau}, \bm{u}_{0:\tau-1})$ generated by \eqref{eq:LDS},
the equality $z_{t+1} = f_z(z_t, u_t)$ holds
for $z_t = E x_t$, $t \in \{0, \dots, \tau-1\}$
if and only if $\mathrm{kern}(E)$ is $A$-invariant.
Furthermore, such a function $f_z$ must be a linear function.
\end{restatable}

The proof of \Cref{prop:invariance_equiv_latent_dynamics},
along with the proofs of other results in this section, follow from standard linear algebra arguments, and are deferred to 
\Cref{sec:appendix:proofs_for_exact_invariance}.
While \Cref{prop:invariance_equiv_latent_dynamics} provides
 a necessary and sufficient characterization of memoryless latent trajectories, for addressing 
 \Cref{prob:latent_stabilization} the backwards implication will be more
 useful. 
 In particular, suppose that $\mathrm{kern}(E)$ is $A$-invariant, i.e., there exists a matrix $L \in \R^{m \times m}$ such that $EA = LE$.\footnote{This holds since $\mathrm{kern}(E)$ is $A$-invariant if and only if
 $\mathrm{row}(E)$ is $A^\T$-invariant.} Then, it is straightforward to check that for $H := EB \in \R^{m \times p}$, the linear function $f_z(z, u) = L z + H u$ is exactly the latent
 dynamics function prescribed by \Cref{prop:invariance_equiv_latent_dynamics}.
 Hence when $m \ll n$, the latent dynamics $(L, H)$ provides
  a compact representation of the original dynamics $(A, B)$.
  
  The next question we turn to is when we can use
  the latent dynamics for control design,
  and in particular when a stabilizing controller $K$ for $(L, H)$ yields a stabilizing output-feedback controller
  $KE$ for $(A, B)$.
  Here, detectability plays an important role,
  as described in the following proposition.

\begin{restatable}[Stabilization equivalence]{myprop}{latentstabilityiffoutputfeedbackstability}
\label{prop:latent_stability_iff_output_feedback_stability}
Let $E \in \R^{m \times n}$ satisfy
(i) the pair $(A, E)$ is detectable,
(ii) $\mathrm{kern}(E)$ is $A$-invariant, and
(iii) $\rank(E) = m$.
Let $L \in \R^{m \times m}$ be the matrix such that $EA=LE$ from condition (ii) and define $H := EB \in \R^{m \times p}$. 
For any feedback gain $K \in \R^{p \times m}$, we have:
\begin{align}
    K \,\,\emph{stabilizes}\,\, (L, H) \iff KE \,\,\emph{stabilizes}\,\, (A, B). \label{eq:stabilization_equality}
\end{align}
Furthermore,
$(L, H)$ is stabilizable
if and only if $(A, B)$ is stabilizable.
\end{restatable}

Conceptually, under the output-feedback controller $u_t=KEx_t$,
\Cref{prop:latent_stability_iff_output_feedback_stability}
shows that the closed-loop dynamics admit a block-triangular
representation with diagonal blocks corresponding to the
latent dynamics $z_{t+1}=(L+HK)z_t$ and the zero dynamics
$x_{t+1}=Ax_t$ on $\ker(E)$.
Invariance of $\ker(E)$ and detectability of $(A,E)$ ensure
that the zero dynamics are Schur stable.
Moreover, full row rank of $E$ ensures that every latent
trajectory can be realized as the encoding of a full-state
trajectory under the same inputs,
which together with the previous conditions yields the stability equivalence in
\Cref{eq:stabilization_equality}.
Hence, \Cref{prop:latent_stability_iff_output_feedback_stability}
forms the basis of our solution to 
\Cref{prob:latent_stabilization}, including providing a concrete definition of what is meant
by \emph{when possible} in \Cref{prob:latent_stabilization}:
the latent dimension $m \in [n]$ should allow the 
existence of an encoder $E \in \R^{m \times n}$
satisfying the three conditions
(i) $(A, E)$ is detectable,
(ii) $\mathrm{kern}(E)$ is $A$-invariant, and (iii)
$\rank(E) = m$.

We next provide the following characterization of when these three conditions
can be satisfied.
Let 
$\sigma(A)$ denote the spectrum of $A$, which we partition into the
non-stable (i.e., unstable plus marginally stable eigenvalues) and stable eigenvalues:
\begin{align}
    \sigma_+(A) := \{ \lambda \in \sigma(A) \mid \abs{\lambda} \geq 1 \}, \qquad \sigma_-(A) := \sigma(A) \setminus \sigma_+(A). \label{eq:eigenvalue_set_definitions}
\end{align}
Next, for an eigenvalue $\lambda \in \C$ of $A$, we let
$\calG_\lambda(A) \subseteq \C^n$ denote the \emph{generalized eigenspace} of $A$ corresponding to the eigenvalue $\lambda$.
We then define the generalized eigenspaces $\calG_+(A)$ corresponding to $\sigma_+(A)$ as 
the direct sum over $\calG_\lambda(A)$ for $\lambda \in \sigma_+(A)$:
\begin{align}
    \calG_+(A) := \bigoplus_{\lambda \in \sigma_+(A)} \calG_\lambda(A). \label{eq:generalized_unstable_eigenspace}
\end{align}
Furthermore, we define $m_+ \in [n]$ and $0 \leq \gamma_- < 1$ to be:
\begin{align}
    m_+ := \sum_{\lambda \in \sigma_+(A)} \mathrm{alg}_A(\lambda) = \dim \calG_+(A), \qquad \gamma_- := \max\{ \abs{\lambda} \mid \lambda \in \sigma_-(A) \},
\end{align}
where $\mathrm{alg}_A(\lambda)$ denotes the \emph{algebraic multiplicity} of the eigenvalue $\lambda$ for matrix $A$.
\begin{restatable}[Equivalent characterization of \Cref{prop:latent_stability_iff_output_feedback_stability} hypothesis]{myprop}{compatibilityassumptionsequiv}
\label{prop:compatibility_assumptions_equiv}
Fix $m \in [n]$. There exists an $E \in \R^{m \times n}$ satisfying the hypothesis of \Cref{prop:latent_stability_iff_output_feedback_stability} if and only if $m \geq m_+$ and
either condition below holds:
\begin{enumerate}[label=(\roman*), nosep, topsep=4pt]
    \item There exists a $\lambda \in \sigma_-(A)$ such that $\lambda \in \R$.
    \item The difference $m - m_+$ is even.
\end{enumerate}
\end{restatable}

With \Cref{prop:latent_stability_iff_output_feedback_stability} and \Cref{prop:compatibility_assumptions_equiv} in place, 
we can now describe an algorithmic solution to
\Cref{prob:latent_stabilization} \emph{given dynamics} $(A, B)$.
The key tool will be the following real Schur decomposition method.
For a \emph{multi-set} $\Lambda$, we let $\mathrm{mult}_{\Lambda}(\lambda)$ denote the multiplicity of $\lambda \in \Lambda$.
Now given $M \in \R^{n \times n}$ 
and a {multi-set} $\Lambda$ of eigenvalues of $M$ where each $\lambda \in \Lambda$ satisfies
(i) $\mathrm{mult}_\Lambda(\lambda) \leq \mathrm{alg}_M(\lambda)$
and (ii) conjugate-invariance 
(i.e., $\mathrm{mult}_\Lambda(\lambda) = \mathrm{mult}_\Lambda(\bar{\lambda})$),
let $(Q, R) = \mathrm{RealSchur}(M, \Lambda)$ denote a pair $(Q, R)$ that satisfies:
\begin{subequations}
\label{eq:real_schur_decomposition}
\begin{align}
    M &= Q \bmattwo{R_{11}}{R_{12}}{0}{R_{22}} Q^\T, \qquad Q = \rvectwo{Q_1}{Q_2} \in O(n), \,\, Q_1 \in \R^{n \times \abs{\Lambda}}, \,\, Q_2 \in \R^{n \times (n-\abs{\Lambda})}, \\
    R_{11} &\in \R^{\abs{\Lambda} \times \abs{\Lambda}}, \,\, \mathrm{eig}(R_{11}) = \Lambda, \qquad R_{22} \in \R^{(n-\abs{\Lambda}) \times (n-\abs{\Lambda})}, \,\, \mathrm{eig}(R_{22}) = \mathrm{eig}(M) \setminus \Lambda,
\end{align}
\end{subequations}
where $R_{11}, R_{22}$ are quasi-upper triangular, 
$\mathrm{eig}(\cdot)$ denotes the multi-set of eigenvalues, and the set equality and difference operators are interpreted
with multiplicity;
given the assumptions on the multi-set $\Lambda$, such a decomposition must exist.\footnote{Note that actually implementing
such a construction in practice, such as with \texttt{scipy.linalg.schur}, may require special care;
for instance, \texttt{scipy}'s API interface only exposes a binary predicate
which cannot distinguish between repeated eigenvalues.
}

Assuming that $m$ satisfies the conditions in
\Cref{prop:compatibility_assumptions_equiv},
\Cref{alg:exact_solution_given_A_B} constructs the encoder and latent dynamics.
Analogous to the set-valued eigenvalue definitions in
\Cref{eq:eigenvalue_set_definitions}, we define the corresponding eigenvalue {multi-sets}:
\begin{align}
    \mathrm{eig}_+(A) := \{ \lambda \in \eig(A) \mid \abs{\lambda} \geq 1 \}, \qquad \mathrm{eig}_-(A) := \mathrm{eig}(A) \setminus \mathrm{eig}_+(A).
\end{align}
All eigenvalue collections and operations in \Cref{alg:exact_solution_given_A_B} are interpreted
with multiplicity.
\begin{algorithm}[H]
\caption{Exact latent dynamics and encoder from known system matrices $(A, B)$}
\label{alg:exact_solution_given_A_B}
\begin{algorithmic}[1]
    \Input $A\in\R^{n\times n}$, $B\in\R^{n\times p}$, and latent dimension $m\in[n]$
    \Require (i) $m \geq m_+$ and (ii) either $\exists \lambda \in \sigma_-(A)$ with $\lambda \in \R$ or $m-m_+$ is even (cf.~\Cref{prop:compatibility_assumptions_equiv})
    \Ensure $E\in\R^{m\times n}$, $L\in\R^{m\times m}$, and $H\in\R^{m\times p}$
    \State $\Lambda_r \gets \mathrm{SelectStable}(\eig(A), n-m)$ \verb| #| Ensures $\Lambda_r \subseteq \eig_-(A)$, $\abs{\Lambda_r} = n-m$, and $\Lambda_r$ is conjugate-invariant\footnotemark
    \State $(Q,R)\gets\mathrm{RealSchur}(A^\T,\eig(A) \setminus \Lambda_r)$, with $Q_1\in\R^{n\times m}, R_{11} \in \R^{m \times m}$ (cf.~\Cref{eq:real_schur_decomposition}) \label{algline:exact_real_schur}
    \State $E\gets Q_1^\T$, $L\gets R_{11}^\T$, and $H\gets EB$
    \State \Return $(E,L,H)$
\end{algorithmic}
\end{algorithm}
\footnotetext{The proof of \Cref{prop:compatibility_assumptions_equiv} provides a specific procedure to selecting $\Lambda_r$.}

One can readily check that conditions in
\Cref{prop:latent_stability_iff_output_feedback_stability} hold
for the output $(E, L, H)$ from \Cref{alg:exact_solution_given_A_B}.
Indeed, orthogonality of $Q$ gives $\ker(E)=\mathrm{col}(Q_2)$,
and transposing the Schur decomposition from Line~\ref{algline:exact_real_schur} gives
$AQ_2=Q_2R_{22}^\T$. Thus $\ker(E)$ is $A$-invariant (condition (ii)),
and also the unobservable subspace of $(A,E)$ with
dynamics represented by $R_{22}^\T$.
Since $\eig(R_{22}) = \Lambda_r$, all eigenvalues 
of $R_{22}$ have modulus strictly less than one,
and therefore $(A,E)$ is detectable (condition (i)).
Finally, since the columns of $Q_1$ are orthonormal, $E = Q_1^\T$ has full row rank (condition (iii)).
Thus \Cref{prop:latent_stability_iff_output_feedback_stability} allows us to conclude that \emph{any} stabilizing control design $K(L, H)$
using $(L, H)$ from \Cref{alg:exact_solution_given_A_B} 
yields a stable controller for $(A, B)$ via 
the feedback gain $K(L, H) E$, as desired
in \Cref{prob:latent_stabilization}.

With \Cref{alg:exact_solution_given_A_B} in place, we now return to the
latent stabilization problem posed in
\Cref{prob:latent_stabilization}, where we do not assume knowledge of $(A, B)$, but instead access to a dataset
$\calD_N$ of trajectories (cf.~\Cref{assmp:data_model}).
We note that one obvious approach is the \emph{certainty-equivalence} or \emph{model-based} approach, where we first
perform system identification, using
$\calD_N$ to estimate $(\widehat{A}_N, \widehat{B}_N)$
via e.g., least-squares regression~\cite{ljung1999system,simchowitz2018learning}, followed
by using the estimates $(\widehat{A}_N, \widehat{B}_N)$
in place of the true $(A, B)$ in
\Cref{alg:exact_solution_given_A_B}.
While such a solution is reasonable, 
it lacks a clear generalization to non-linear systems, as it heavily exploits 
linear structure \emph{in the actual algorithm};
generalizing \Cref{alg:exact_solution_given_A_B} for learning non-linear latent dynamics
and encoders is not obvious.
What we will develop in the remainder of the paper, on the other hand,
is an algorithm procedure that is general purpose (e.g., can be directly applied to learning non-linear models from data collected from non-linear systems), 
and furthermore when the procedure is applied
to a linear system, provably recovers the
stabilization properties 
stated in \Cref{prob:latent_stabilization}.

\section{Non-Linear Forward and Backward Conjugacy Losses}
\label{sec:problem_formulation:losses}

With the results of \Cref{sec:problem_formulation:invariance_subspaces} in place, we now describe in detail
the latent prediction and reconstruction losses introduced in
\Cref{sec:introduction:example}, 
which form the cornerstone of our approach.
To both emphasize the generality of the following losses,
and for setting the
stage for the non-linear experiments (cf.~\Cref{sec:experiments:nonlinear}),
we will assume (unless otherwise stated) in this section that the underlying dynamics are a non-linear system
$x_{t+1} = f(x_t, u_t)$ for $f : \R^n \times \R^p \mapsto \R^n$,
and that given a latent dimension $m \in [n]$, 
we are interested in learning both a non-linear 
encoder $E : \R^n \mapsto \R^m$ and non-linear latent dynamics
$z_{t+1} = f_z(z_t, u_t)$ for $f_z : \R^m \times \R^p \mapsto \R^m$.
For $t\geq 1$ and $\bm{u}_{0:t-1}:=(u_0,\ldots,u_{t-1})$,
we define the following discrete-time flow maps:
\begin{align}
    \Phi^t_f(x_0, \bm{u}_{0:t-1})
    &:= f\left(\cdots f(f(x_0,u_0),u_1)\cdots,u_{t-1}\right), \\
    \Phi^t_{f_z}(z_0, \bm{u}_{0:t-1})
    &:= f_z\left(\cdots f_z(f_z(z_0,u_0),u_1)\cdots,u_{t-1}\right).
\end{align}
Here, each expression contains $t$ applications of its dynamics map;
for $t=0$, we use the convention
that $\Phi^0_f(x_0,\varnothing)=x_0$
and $\Phi^0_{f_z}(z_0,\varnothing)=z_0$.
With this notation in place, we first start with 
the following multi-step forward conjugacy loss.
\begin{mydef}[Multi-step forward conjugacy]
\label{def:mstep_forward_conjugacy_loss_general}
For $T \in \N_+$, the (population-level) 
multi-step forward conjugacy loss 
is defined as:
\begin{align}
    J_{\mathrm{fwd}}^{(T)}(E, f_z) := \frac{1}{T}\sum_{t=1}^{T}\mathbb{E}_{x_0,\bm{u}_{0:t-1}}\left[ \norm{E \circ \Phi^t_f(x_0,\bm{u}_{0:t-1})-\Phi^t_{f_z}(E(x_0),\bm{u}_{0:t-1})}^2 \right], \label{eq:mstep_forward_conjugacy_loss_general}
\end{align}
where the expectation is taken jointly over draws of initial conditions and input sequences $(x_0, \bm{u}_{0:T-1})$.
\end{mydef}
\Cref{def:mstep_forward_conjugacy_loss_general} measures how accurately
the discrete-time flows of the latent model $f_z$ are able to predict 
the latent encodings of the actual trajectories of $f$, over a horizon length of size $T$.
As stated, such a loss is clearly insufficient on its own, as trivially setting
$E \equiv 0$ and $f_z \equiv 0$ achieves zero loss, a form of \emph{representation collapse}; we will return to the
question of preventing this collapse in \Cref{sec:encoder_only_regularization}.

The next loss we consider is the following multi-step backward conjugacy loss,
which introduces a non-linear decoder $D : \R^m \mapsto \R^n$ and measures the error of reconstruction along the latent flows.
\begin{mydef}[Multi-step backward conjugacy]
\label{def:mstep_backwards_conjugacy_loss_general}
For $T \in \N_+$, the (population-level) 
multi-step backward conjugacy loss 
is defined as:
\begin{align}
    J_{\mathrm{bwd}}^{(T)}(E, D, f_z) := \frac{1}{T+1}\sum_{t=0}^{T}\mathbb{E}_{x_0,\bm{u}_{0:t-1}}\left[ \norm{\Phi^t_f(x_0,\bm{u}_{0:t-1})-D \circ \Phi^t_{f_z}(E(x_0),\bm{u}_{0:t-1})}^2 \right], \label{eq:mstep_backwards_conjugacy_loss_general}
\end{align}
where the expectation is taken jointly over draws of initial conditions and input sequences $(x_0, \bm{u}_{0:T-1})$.
\end{mydef}
Unlike the forward loss, the backward loss in \Cref{def:mstep_backwards_conjugacy_loss_general} does not suffer from the same
representation collapse issues; indeed, observe that the $t=0$ summand is equal to 
$\E_{x_0}\left[ \norm{ x_0 - D \circ E(x_0) }^2 \right]$,
which is a standard auto-encoder reconstruction loss~\cite[see e.g.,][]{hinton2006reducing}.

\begin{rmk}[Empirical losses]
\label{rmk:empirical_losses}
Both losses in \Cref{def:mstep_forward_conjugacy_loss_general}
and \Cref{def:mstep_backwards_conjugacy_loss_general} are described
in their population-level variants, with expectations taken over
jointly over $(x_0, \bm{u}_{0:T-1})$.
In practice though, we only observe a trajectory dataset
$\calD_N = \{ (\bm{x}_{0:T}^{(i)}, \bm{u}_{0:T-1}^{(i)}) \}_{i=1}^{N}$ (cf.\ the data model described in \Cref{assmp:data_model}).
From $\calD_N$ we form the empirical losses of 
\eqref{eq:mstep_forward_conjugacy_loss_general}
and 
\eqref{eq:mstep_backwards_conjugacy_loss_general}
as:
\begin{align}
    \widehat{J}_{\mathrm{fwd}}^{\;(T)}(E, f_z) &:= \frac{1}{NT} \sum_{i=1}^{N} \sum_{t=1}^{T} \norm{E(x_t^{(i)}) - \Phi^t_{f_z}(E(x_0^{(i)}),\bm{u}^{(i)}_{0:t-1})}^2, \\
    \widehat{J}_{\mathrm{bwd}}^{\;(T)}(E, D, f_z) &:= \frac{1}{N(T+1)} \sum_{i=1}^{N} \sum_{t=0}^{T} \norm{ x_t^{(i)} - D \circ \Phi^t_{f_z}(E(x_0^{(i)}),\bm{u}^{(i)}_{0:t-1})}^2.
\end{align}
Our analysis will proceed by first studying the properties of the
minimizers of the population-level losses, followed by concentration of measure
arguments to show that with high probability, 
$\widehat{J}_{\mathrm{fwd}}^{\;(T)}(E, f_z) \asymp J_{\mathrm{fwd}}^{(T)}(E, f_z)$
and $\widehat{J}_{\mathrm{bwd}}^{\;(T)}(E, D, f_z) \asymp J_{\mathrm{bwd}}^{(T)}(E, D, f_z)$ hold uniformly over our optimization variables (cf.~\Cref{sec:results:finite_trajectories}).
\end{rmk}

\subsection{Instantiations of Non-Linear Losses to Linear Systems}
\label{sec:problem_formulation:losses_LTI}

We now instantiate the
losses 
\eqref{eq:total_backwards_loss} and
\eqref{eq:total_forwards_loss} 
for our LTI system \eqref{eq:LDS}, with linear encoders, decoders, and latent dynamics.
Furthermore, we will utilize the distribution over initial conditions 
and control inputs $(x_0, \bm{u}_{0:T-1})$ specified in
the data model from \Cref{assmp:data_model}.

\begin{myfact}[Losses instantiated for linear encoders, decoders, and dynamics]
\label{prop:LTI_system_losses}
Represent the encoder as $E(x) = Ex$ for $E \in \R^{m \times n}$, 
decoder as $D(x) = Dx$ for $D \in \R^{n \times m}$, 
and latent dynamics as $f_z(z, u) = Lz + H u$ for $L \in \R^{m \times m}$ and $H \in \R^{m \times p}$.
Fix a rollout horizon $T \in \N_+$ and draw $(x_0, \bm{u}_{0:T-1})$
according to the distribution specified in \Cref{assmp:data_model} under the LTI dynamics \eqref{eq:LDS}.
We have:
\begin{align}
    J_{\mathrm{fwd}}^{(T)}(E, f_z) = J_{\mathrm{fwd}}^{(T)}(E, L, H) &:= \frac{\sigma_x^2}{T}\sum_{t=1}^{T} \norm{E A^t - L^t E}_F^2 + \frac{\sigma_u^2}{T} \sum_{k=0}^{T-1} \beta_{T}(k) \norm{ E A^k B - L^k H }_F^2, \label{eq:LTI_forward_loss} \\
    J_{\mathrm{bwd}}^{(T)}(E, D, f_z) = J_{\mathrm{bwd}}^{(T)}(E, D, L, H) &:= \frac{\sigma_x^2}{T+1}\sum_{t=0}^{T} \norm{A^t - DL^t E}_F^2 + \frac{\sigma_u^2}{T+1} \sum_{k=0}^{T-1} \beta_{T}(k) \norm{ A^k B - D L^k H }_F^2, \label{eq:LTI_backward_loss}
\end{align}
where $\beta_{T}(k) := \sum_{j=0}^{T-1-k} \alpha_T(j)$ for
$k \in \{0, \dots, T-1\}$.
\end{myfact}

The two concrete loss expressions
in \Cref{prop:LTI_system_losses}
directly expose the symmetry
and scaling properties
of the forward and backwards losses,
as discussed in the following remark.

\begin{rmk}[Loss function symmetries and scaling]
\label{remark:loss_function_symmetries}
The following symmetries for the forward and backward losses in \eqref{eq:LTI_forward_loss} and \eqref{eq:LTI_backward_loss} hold for any rollout horizon $T \in \N_+$:
\begin{align*}
    J^{(T)}_{\mathrm{fwd}}(E, L, H) &= J^{(T)}_{\mathrm{fwd}}(SE, SLS^{-1}, SH), \quad S \in O(m), \\
    J^{(T)}_{\mathrm{bwd}}(E, D, L, H) &= J^{(T)}_{\mathrm{bwd}}(SE, DS^{-1}, SLS^{-1}, SH), \quad S \in \mathrm{GL}(m).
\end{align*}
We also have the following scaling property
for the forward loss:
\begin{align*}
    J^{(T)}_{\mathrm{fwd}}(c E, L, c H) = c^2 J^{(T)}_{\mathrm{fwd}}(E, L, H), \quad c \in \R,
\end{align*}
\end{rmk}
\Cref{remark:loss_function_symmetries} highlights the need for regularization in both the forward and backward losses.
For the backward loss, although representation collapse is not a direct issue, regularization is still useful for fixing the scale of either the encoder or the decoder.
Without it, the encoder could become arbitrarily small while the decoder grows arbitrarily large to compensate, or vice versa.
While
for the backward loss we do have the design choice
to either regularize the encoder, the decoder, or both, to facilitate a more direct comparison with the forward loss we will only regularize the encoder.

\subsubsection{Losses over Finite Trajectories}
\label{sec:results:finite_trajectories}

We now consider both forward and backwards losses over finite
trajectories $\calD_{N} := \{(\bm{x}_{0:T}^{(i)}, \bm{u}_{0:{T-1}}^{(i)})\}_{i=1}^{N}$, where each trajectory is drawn i.i.d.\ from 
the trajectory distribution 
described in \Cref{assmp:data_model}.
For linear encoders, decoders, and dynamics,
the finite-sample losses in \Cref{rmk:empirical_losses} simplify to:
\begin{align}
    \wh{J}^{(T)}_{\mathrm{fwd}}(E, L, H) &:= \frac{1}{NT} \sum_{i=1}^{N} \sum_{t=1}^{T} \norm{ E x_t^{(i)} - \Phi^t_{f_z}( E x_0^{(i)}, \bm{u}_{0:t-1}^{(i)}) }^2, \quad f_z(z, u) = Lz + H u, \label{eq:fwd_loss_finite_sample} \\
    \wh{J}^{(T)}_{\mathrm{bwd}}(E, D, L, H) &:= \frac{1}{N(T+1)} \sum_{i=1}^{N} \sum_{t=0}^{T} \norm{ x_t^{(i)} - D\Phi_{f_z}^{t}( E x_0^{(i)}, \bm{u}_{0:t-1}^{(i)}) }^2. \label{eq:bwd_loss_finite_sample}
\end{align}
We now show that as long as the number of trajectories $N$ is large enough, 
we can interchange the finite-sample losses above with their population
counterparts.

\begin{restatable}[Approximate equivalence of finite-sample loss]{myprop}{populationandfinitesampleequivalence}
\label{prop:population_and_finite_sample_equivalence}
Fix a $\delta \in (0, 1)$ and suppose that $N \geq C \cdot \left[n + Tp + \log(2/\delta)\right]$, where $C > 0$ is a universal constant.
Then with probability at least $1-\delta$,
for every $(E, D, L, H)$, we have:
\begin{align*}
    \frac{1}{2} J_{\mathrm{fwd}}^{(T)}(E, L, H) &\leq \wh{J}_{\mathrm{fwd}}^{(T)}(E, L, H) \leq \frac{3}{2} J_{\mathrm{fwd}}^{(T)}(E, L, H), \\
        \frac{1}{2} J_{\mathrm{bwd}}^{(T)}(E, D, L, H) &\leq \wh{J}_{\mathrm{bwd}}^{(T)}(E, D, L, H) \leq \frac{3}{2} J_{\mathrm{bwd}}^{(T)}(E, D, L, H).
\end{align*}
Note the constants $1/2$ and $3/2$ can be made arbitrarily close to one at the expense of increasing $C$.
\end{restatable}

\Cref{prop:population_and_finite_sample_equivalence}
will allow us to transfer over
results proven for the population-loss over to the finite-sample
losses with high probability whenever $N \gtrsim n + Tp$,
which we will make use of shortly.

\subsection{Representation Collapse and Encoder-Only Regularization}
\label{sec:encoder_only_regularization}

The issue of representation collapse has a rich literature in self-supervised and unsupervised ML~\citep[see e.g.,][]{zbontar2021barlow,bardes2022vicreg,jing2022understanding}, and still remains an active area of research. 
In our context,
we will consider encoder-only regularizers $R(E)$ (as opposed to joint regularizers that also include the decoder) which ensure that
the pushforward $E \# \mu$ for some base measure $\mu$
is either (approximately) isotropic or has independent
coordinates.
Some (population-level) examples such  encoder-only regularizers are listed below:
\begin{align}
    R_{\mathrm{orth}}(E; \mu) &:= \E_{x \sim \mu} \norm{  I_m - J_E(x) J_E(x)^\T }_F^2, \quad J_E(x) := \frac{\rmd}{\rmd x} E(x), \label{eq:regularizer_orth} \\
    R_{\mathrm{sr}}(E; \mu) &:= \E_{a \sim \mathrm{Unif}(\mathbb{S}^{m-1})} \int \bigabs{ \E_{x \sim \mu}[ e^{j t \ip{a}{E(x)}} ] - e^{-t^2/2} }^2 e^{-t^2/2} \,\rmd t, \quad j := \sqrt{-1}, \label{eq:regularizer_SR} \\
    R_{\mathrm{vr}}(E; \mu, \beta) &:= \frac{1}{m}\sum_{i=1}^m
        (1 - (C_\mu(E))_{ii})^2
        +\frac{\beta}{m}\sum_{i\neq j}(C_\mu(E))_{ij}^2, \quad \beta > 0, \,\, C_\mu(E) := \Var_{x \sim \mu}( E(x) ). \label{eq:regularizer_VR}
\end{align}
The regularizer $R_\mathrm{orth}$ incentivizes orthogonality of the encoder via its Jacobian, $R_{\mathrm{sr}}$ is
a population-level variant of the SIGReg loss used in LeJEPA~\cite{balestriero2025lejepa},
and $R_{\mathrm{vr}}$ is a 
square-variance instantiation of
VICReg~\cite{bardes2022vicreg};
the base measure $\mu$ is a distribution on the state space
for all three regularizers. 
We introduce several regularizers here not to advocate for one particular approach,
but rather to illustrate that representation collapse in practice can be discouraged through
a variety of different mechanisms.

Towards setting the stage for our analysis for LTI systems, we next formally define the 
key mathematical properties we need
out of the encoder-only regularizers applied to linear encoders.
In our case, it suffices for $R(E)$ to encourage that the matrix $EE^\T$ is an approximate isometry.

\begin{mydef}[Regularizer admissibility]
\label{def:forward_admissible_regularizer}
An encoder-only regularizer $R : \R^{m \times n} \mapsto \R_{\geq 0}$ is \emph{admissible} 
if
(i) $R(E) = 0$ for $E^\T \in O(n, m)$, and
(ii)
there exists $\sigma_0 \in (0, 1]$,
$\sigma_1 \geq 1$,
and $\e_0 > 0$ such that
for all $E \in \R^{m \times n}$:
\begin{align}
    R(E) \leq \e_0 \implies \sigma_0 I_m \preccurlyeq EE^\T \preccurlyeq \sigma_1 I_m.
\end{align}
\end{mydef}

The admissibility properties described in
\Cref{def:forward_admissible_regularizer}
are broad enough to cover the regularizers
described previously in 
\Cref{eq:regularizer_orth,eq:regularizer_SR,eq:regularizer_VR}.

\begin{restatable}[Admissibility of $R_\mathrm{orth}$, $R_\mathrm{sr}$, $R_\mathrm{vr}$]{myprop}{regularizeradmissibility}
\label{prop:regularizer_admissibility}
For linear encoders $E(x)=Ex$ with $E\in\R^{m\times n}$,
the regularizers $R_\mathrm{orth}$, $R_{\mathrm{sr}}$, and $R_{\mathrm{vr}}$ (for every $\beta>0$)
with base measure $\mu_n=\sfN(0,I_n)$ are admissible encoder-only regularizers
(cf.~\Cref{def:forward_admissible_regularizer}).
\end{restatable}

The proof of \Cref{prop:regularizer_admissibility} is deferred to
\Cref{sec:appendix:proofs_for_regularizer_admissibility}.
With these losses and regularizer admissibility conditions in place, we are now ready to turn to our main results next.

\section{Main Results}
\label{sec:results}

With the loss notation from \Cref{sec:problem_formulation:losses} in place, we now state our main results
regarding the optimizers of the following latent dynamics modeling losses
that directly mirror those presented
\Cref{eq:loss_combo_bwd} and 
\Cref{eq:loss_combo_sigreg} from \Cref{sec:introduction:example}:
\begin{tcolorbox}[colback=blue!5, colframe=blue!55!black,
    boxrule=0.5pt, sharp corners, left=4pt, right=4pt, top=2pt, bottom=2pt]
\setlength{\abovedisplayskip}{0pt}
\setlength{\abovedisplayshortskip}{0pt}
\begin{align}
   \calL_{\mathrm{bwd}}^{(T)}(E, D, L, H) &:= J_{\mathrm{bwd}}^{(T)}(E, D, L, H) + \lambda_{\mathrm{fwd}} \cdot J_{\mathrm{fwd}}^{(T)}(E, L, H) + \lambda_{\mathrm{reg}} \cdot R(E),\label{eq:total_backwards_loss} \\
   \calL_{\mathrm{fwd}}^{(T)}(E, L, H) &:= J_{\mathrm{fwd}}^{(T)}(E, L, H) + \lambda_{\mathrm{reg}} \cdot R(E).\label{eq:total_forwards_loss}
\end{align}
\end{tcolorbox}

\subsection{Latent Stabilization via Multi-Step Backward Losses}
\label{sec:results:multi_step_backward}

We first show how minimizing 
the backward loss
$\calL_{\mathrm{bwd}}^{(T)}$ in 
\Cref{eq:total_backwards_loss}
yields a solution
to the latent stabilization problem (cf.~\Cref{prob:latent_stabilization}).
Specifically, we will learn a linear encoder
$E$, linear decoder $D$, and latent LTI dynamics $(L, H)$
by solving the optimization problem:
\begin{align}
    \left(\bar{E}_T, \bar{D}_T, \bar{L}_T, \bar{H}_T\right)\in\underset{\substack{E\in \R^{m\times d},\ D\in\mathbb{R}^{n\times m}, \\ L\in\mathbb{R}^{m\times m},\ H\in\mathbb{R}^{m\times p}}}{\argmin} \mathcal{L}^{(T)}_\mathrm{bwd}(E, D, L, H). \label{eq:bwd_solution}
\end{align}
While is it straightforward to check that
the optimizers of \eqref{eq:bwd_solution} are
attained, due to the loss function symmetries
noted in \Cref{remark:loss_function_symmetries},
the solution is \emph{not} unique.
Hence in the sequel, the pair
$(\bar{E}_T, \bar{D}_T, \bar{L}_T, \bar{H}_T)$
will refer to \emph{any} element in the solution set; our results will not depend on 
choosing a specific solution.

We begin by presenting the main set of assumptions 
used in our analysis.
The first assumption is straightforward, requiring the
non-stable subspace of $A$ to be both non-trivial (i.e., $\rho(A) \geq 1$)
and a strict subset of $\R^n$.
This assumption is used to focus on the interesting regime of the problem:
if $\rho(A) < 1$ then stabilization 
would trivially be solved by taking $K=0$;
on the other hand if all eigenvalues of $A$ were non-stable, then 
no lower-dimensional encoder (i.e., $m < n$)
can simultaneously satisfy exact invariance and detectability.
While our analysis could easily be modified to account for these two special cases, this would
come at the expense of readability, and hence we 
choose to omit them.

\begin{assmp}[Non-stable subspace dimension]
\label{assmp:non_trivial_unstable_subspace}
The non-stable subspace dimension
satisfies $1 \leq m_+ < n$.
\end{assmp}

Our next assumption concerns the choice of the latent dimension $m \in [n]$. It comes directly from our discussion in 
\Cref{sec:problem_formulation:invariance_subspaces}, specifically
\Cref{prop:compatibility_assumptions_equiv},
ensuring that there actually exists an encoder $E$ and
latent dynamics $(L, H)$ so that 
stabilizing $(L, H)$ is equivalent to stabilizing $(A, B)$.

\begin{assmp}[Latent dimension compatibility]
\label{assmp:latent_dimension_compatibility}
The latent dimension $m \in [n]$ satisfies $m \geq m_+$ and 
either (i) there exists a real-valued $\lambda \in \sigma_-(A)$ or (ii) the difference $m - m_+$ is even.
\end{assmp}

Our final assumption asks that the regularizer $R(E)$ is admissible.
\begin{assmp}[Regularizer admissibility]
\label{assmp:regularizer_assumption}
The regularizer $R(E)$
is admissible (cf.~\Cref{def:forward_admissible_regularizer})
with parameters $(\sigma_0, \sigma_1, \e_0)$.
\end{assmp}

To state our first result, we introduce notation for the following
\emph{one-sided subspace distance}:
given two subspaces $A, B \subseteq \R^n$ (with $\dim(A)$ not necessarily equal to $\dim(B)$), we define:
\begin{align}
    \Delta_{\mathrm{op}}(B ; A) := \opnorm{ P^\perp_A P_B }, \qquad \Delta_{F}(B ; A) := \norm{ P^\perp_A P_B }_F.
\end{align}
Recall that $\calG_+(A^\T)$ denotes the \emph{left} generalized eigenspaces of $A$ corresponding to the non-stable
eigenvalues $\sigma_+(A)$ (cf.~\Cref{eq:generalized_unstable_eigenspace}). Since $A$ is a real-valued matrix,
the subspace $\calG_+(A^\T)$ is conjugate invariant;
we denote the real part of this subspace as
$\calG^\R_+(A^\T) := \calG_+(A^\T) \cap \R^n$.
The following theorem is our first main result, which
provides quantitative error bounds
for any optimizer tuple of $\calL_{\mathrm{bwd}}^{(T)}$.
In the interest of readability, we present a version of the result where 
(a) the problem
dependent constants are not fully specified, and 
(b) we considered a fixed setting $\alpha_T(t) \equiv 1/T$
for the input variance (cf.~\Cref{assmp:data_model});
our more general theorem
(cf.~\Cref{thm:multistep_backwards_full})
provides the full description of the problem specific constants,
in addition to allowing for more general choices of $\alpha_T(t)$.

\begin{mythm}[Recovery error bounds---simplified version of \Cref{thm:multistep_backwards_full}]
\label{thm:informal_recovery}
Suppose that 
\Cref{assmp:stabilizability,assmp:non_trivial_unstable_subspace,assmp:latent_dimension_compatibility,assmp:regularizer_assumption} hold.
For the distribution over $(x_0, \bm{u}_{0:T-1})$ specified in
\Cref{assmp:data_model} with input variance $\alpha_T(t) \equiv 1/T$, 
there exists constants $C_i \equiv C_i(A, B, n, m, \sigma_x^2, \sigma_u^2)$, $i=1,\dots,4$, 
such that
any minimizer $(\bar{E}_T, \bar{D}_T, \bar{L}_T, \bar{H}_T)$ of \Cref{eq:bwd_solution} satisfies
the following error bounds:
\begin{enumerate}[label=(\roman*)]
    \item \underline{Approximate detectability of $(A, \bar{E}_T)$}:
    The row space of the learned encoder approximately contains the left non-stable subspaces of $A$, 
    in the sense that:
    \begin{gather}
        \Delta^2_F\left( \calG^\R_+(A^\T) ; \mathrm{row}(\bar{E}_T) \right) \leq \frac{C_1}{T}. \label{eq:informal_recovery:subspace_distance}
    \end{gather}

    \item \underline{Approximate memoryless reconstruction}: The latent model $(\bar{L}_T, \bar{H}_T)$ provides approximate memoryless reconstruction of the trajectories of $(A, B)$ encoded by $\bar{E}_T$:
    \begin{gather}
        \norm{\bar{E}_TA-\bar{L}_T\bar{E}_T}_F^2\leq \frac{C_2}{\lambda_\mathrm{fwd}}, \qquad \norm{ \bar{E}_T B - \bar{H}_T }_F^2 \leq \frac{C_3}{\lambda_{\mathrm{fwd}}}. \label{eq:informal_recovery:approximate_invariance}
    \end{gather}

    \item \underline{Approximate isometry of $\bar{E}_T\bar{E}_T^\T$}: The learned encoder is an approximate isometry:
    \begin{align}
        \lambda_{\mathrm{reg}} \geq \frac{C_4}{\e_0} \implies \sigma_0 I_m \preccurlyeq \bar{E}_T \bar{E}_T^\T \preccurlyeq \sigma_1 I_m. \label{eq:informal_recovery:approximate_isometry}
    \end{align}

\end{enumerate}

\end{mythm}

We emphasize that the problem dependent constants $C_i(A, B, n, m, \sigma_x^2, \sigma_u^2)$ appearing in \Cref{thm:informal_recovery} do \emph{not} depend on $T$, and hence the subspace distance \eqref{eq:informal_recovery:subspace_distance} and approximate forward conjugacy error terms \eqref{eq:informal_recovery:approximate_invariance}
above can be made arbitrarily small
by both increasing the multi-step horizon length $T$
and the regularization parameters $\lambda_{\mathrm{reg}}, \lambda_{\mathrm{fwd}}$.
The bounds (i)--(iii) in \Cref{thm:informal_recovery}
are approximate and quantitative variants
of the corresponding conditions (i)--(iii) described in \Cref{prop:latent_stability_iff_output_feedback_stability}.
To see this, first note that the subspace distance 
$\Delta_F\left( \calG^\R_+(A^\T) ; \mathrm{row}(\bar{E}_T) \right)$
from \eqref{eq:informal_recovery:subspace_distance}
provides a quantitative measure of the
deviation from the detectability of $(A, \bar{E}_T)$.
Indeed, if $\Delta_F\left( \calG^\R_+(A^\T) ; \mathrm{row}(\bar{E}_T) \right) = 0$, then
$\calG^\R_+(A^\T) \subseteq \mathrm{row}(\bar{E}_T)$;
taking orthogonal complements shows that 
\begin{align*}
    \ker(\bar{E}_T) \subseteq \left( \calG^\R_+(A^\T) \right)^\perp = \calG_-^\R(A) := \left( \bigoplus_{\lambda \in \sigma_-(A)} \calG_\lambda(A) \right) \cap \R^n,
\end{align*}
which implies that $(A,\bar{E}_T)$ is detectable.
Furthermore, the two residuals in
\eqref{eq:informal_recovery:approximate_invariance}
measure the degree to which the latent model
$f_z(z, u) = \bar{L}_T z + \bar{H}_T u$ 
provides a memoryless reconstruction (cf.~\Cref{prop:invariance_equiv_latent_dynamics}) 
of trajectories
of $(A, B)$ that are encoded by $\bar{E}_T$;
equivalently, \eqref{eq:informal_recovery:approximate_invariance} provides a quantitative measure of the $A$-invariance
of $\ker(\bar{E}_T)$.
Finally, the isometry lower bound $\bar{E}_T \bar{E}_T^\T \succcurlyeq \sigma_0 I_m$ in \eqref{eq:informal_recovery:approximate_isometry}
ensures that $\rank(\bar{E}_T) = m$.

The mechanism behind \Cref{thm:informal_recovery} is a separation of time scales in the backward loss. Intuitively,
suppose that some direction of the left non-stable subspace of $A$ is missing from $\mathrm{row}(\bar{E}_T)$. Then,
there exists initial conditions $x_0\in\mathrm{ker}(\bar{E}_T)$ that excite an non-stable mode of $A$,
while producing the latent state $z_0 =\bar{E}_T x_0 = 0$, and hence the decoded prediction $\bar{D}_T \Phi^t_{f_z}(\bar{E}_T x_0, \bm{u}_{0:t-1})$ 
carries no information about a component of $\Phi_f^t(x_0, \bm{u}_{0:t-1})$ that does not decay with $t$. 
The time-averaged sum in \eqref{eq:mstep_backwards_conjugacy_loss_general} therefore incurs an error proportional to the squared misalignment, 
with a constant that does not vanish as $T$ grows. Conversely, a latent model that captures the non-stable modes exactly, 
and discards only stable directions, incurs a reconstruction error that decays geometrically in $t$,
and is of order $1/T$ when averaged over time. 
Both the prediction horizon and backward loss is essential to this argument. At $T=0$, the backward loss reduces to the auto-encoder reconstruction loss, which under the isotropic initial condition distribution of \Cref{assmp:data_model} is minimized by the orthogonal projector onto \emph{any} $m$-dimensional subspace. Likewise, $J_\mathrm{fwd}^{(T)}$ and $R$ can vanish at encoders whose row space consists entirely of stable left eigenvectors of $A$; hence on their own, $J_\mathrm{fwd}^{(T)}$ and the regularizer $R$ are insufficient (cf.~\Cref{sec:results:insufficiency}),
and their role in $\calL_\mathrm{bwd}^{(T)}$ is instead to enforce forward conjugacy and encoder conditioning once the backward loss has selected the correct encoder subspace.

The conclusions of \Cref{thm:informal_recovery} state that the learned latent model \emph{approaches} the class of solutions characterized by \Cref{prop:latent_stability_iff_output_feedback_stability}. 
However for any finite horizon $T$, the residuals
\eqref{eq:informal_recovery:approximate_invariance}
and \eqref{eq:informal_recovery:approximate_invariance}
can be small but still non-zero. 
\Cref{prop:latent_stability_iff_output_feedback_stability} then no longer applies directly, and an arbitrary stabilizing gain $\bar{K}_T$ for $(\bar{L}_T, \bar{H}_T)$ need not stabilize the closed-loop
$A + B \bar{K}_T \bar{E}_T$. Our next two results shows that this gap is closed by choosing the stabilizing gain $\bar{K}_T$ to be a minimal
$H_\infty$-norm latent controller.
We do this in two parts:
(i) we first give a \emph{structural} condition
on an arbitrary tuple $(E_T, L_T, H_T)$
(not necessarily the minimizers of $\calL_{\mathrm{bwd}}^{(T)}$), providing
sufficient conditions on approximate 
detectability
and encoder invariance to ensure a minimum-norm $H_\infty$-design
stabilizes $(A, B)$, followed by
(ii) a setting of hyperparameters $(T, \lambda_{\mathrm{fwd}}, \lambda_{\mathrm{reg}})$ for the loss $\calL_{\mathrm{bwd}}^{(T)}$
that ensure the structural conditions in (i) hold.
We start with our structural result.
Again we present a simplified version where the 
problem specific constants are not explicitly specified; the more general
theorem (\Cref{thm:H_infty_closed_loop_stability})
is presented later.

\begin{mythm}[Minimum $H_\infty$-norm closed-loop stability---simplified version of \Cref{thm:H_infty_closed_loop_stability}]
\label{thm:resolvent_gain_bound}
Suppose \Cref{assmp:stabilizability,assmp:non_trivial_unstable_subspace} hold.
Given $\sigma_0 > 0$
and $M_0 \geq C_0(A, B, \sigma_0)$, 
for any parameters
$E_T \in \R^{m \times n}$,
$L_T \in \R^{m \times m}$ and $H_T \in \R^{m \times p}$ satisfying $\sigma_{\min}(E_T E_T^\T) \geq \sigma_0$ and:
\begin{align*}
    \max\left\{ \Delta_{\mathrm{op}}( \calG_+^\R(A^\T); \mathrm{row}(E_T) ), \opnorm{ E_T A - L_T E_T }, \opnorm{ E_T B - H_T } \right\} &\leq \e_{\mathrm{stab}}(A, B, \sigma_0, M_0, \opnorm{E_T} ),
\end{align*}
the following conditions hold:
\begin{enumerate}[label=(\roman*)]
    \item The pair $(L_T, H_T)$ is stabilizable.
    \item A minimum $H_\infty$-norm latent controller
    \begin{align}
        K_T \in \argmin_{K \in B_{p \times m}^{(\mathrm{op})}(M_0)}
        \sup_{z \in \mathbb{D}^c}
        \opnorm{(zI_m-L_T-H_TK)^{-1}}
    \end{align}
    exists,
    and every minimizer stabilizes $A+BK_TE_T$ with $H_\infty$-gain bound:
    \begin{align}
        \sup_{z \in \mathbb{D}^c}
        \opnorm{(zI_n-A-BK_TE_T)^{-1}} \leq C_1(A, B, \sigma_0, M_0, \opnorm{E_T}).
    \end{align}
\end{enumerate}
Here, $C_0, C_1, \e_{\mathrm{stab}}$ are 
constants depending only on their respective
inputs.
\end{mythm}

Next, combining the error bounds from
\Cref{thm:informal_recovery} 
with the closed-loop stability transfer
guarantees from \Cref{thm:resolvent_gain_bound} yields
the following end-to-end closed-loop guarantee
(see \Cref{thm:end_to_end_closed_loop} for the full result).

\begin{mycorollary}[End-to-end closed-loop guarantees---simplified version of \Cref{thm:end_to_end_closed_loop}]
\label{cor:latent_stability_transfer}
Suppose that
\Cref{assmp:stabilizability,assmp:non_trivial_unstable_subspace,assmp:latent_dimension_compatibility,assmp:regularizer_assumption} hold.
For $M_0 \geq C_0(A, B, \sigma_0)$ and the distribution over $(x_0, \bm{u}_{0:T-1})$ specified in
\Cref{assmp:data_model} with input variance $\alpha_T(t) \equiv 1/T$,  
suppose the loss
$\calL^{(T)}_{\mathrm{bwd}}$ hyperparameters 
satisfy $T \geq C_1$,
$\lambda_{\mathrm{fwd}} \geq C_2$, and
$\lambda_{\mathrm{reg}} \geq C_3$.
Then a minimum $H_\infty$-norm controller
\begin{align*}
    K_T \in \argmin_{K \in B_{p \times m}^{(\mathrm{op})}(M_0)}
    \sup_{z \in \mathbb{D}^c}
    \opnorm{(zI_m-\bar{L}_T-\bar{H}_T K)^{-1}}
\end{align*}
exists, and every minimizer
stabilizes $A + B \bar{K}_T \bar{E}_T$
with the following closed-loop $H_\infty$-norm bound:
\begin{align*}
    \sup_{z \in \mathbb{D}^c} \opnorm{ (zI_n - A - B \bar{K}_T \bar{E}_T)^{-1}} \leq C_4.
\end{align*}
Above, $C_i \equiv C_i(A, B, n, m, \sigma_x^2, \sigma_u^2, M_0, \sigma_0, \sigma_1)$ for $i=1,\dots,4$.
\end{mycorollary}

\Cref{cor:latent_stability_transfer} shows that for a suitable prediction horizon, the parameters $(\bar{E}_T, \bar{D}_T, \bar{L}_T, \bar{H}_T)$ learned from optimizing $\calL_{\mathrm{bwd}}^{(T)}$ yield a latent model
$(\bar{L}_T, \bar{H}_T)$ such that 
computing a minimum $H_\infty$-norm controller $\bar{K}_T$
induces a stabilizing closed-loop
$A+B\bar{K}_T\bar{E}_T$.
Up to the difference between the population and finite-trajectory loss (cf.~\Cref{rmk:empirical_losses}), this provides a solution to the latent stabilization problem (cf.~\Cref{prob:latent_stabilization}).

To bridge the gap between the population and finite-sample loss, 
recall from \Cref{sec:results:finite_trajectories}
that with high probability,
as long as the number of trajectory rollouts $N$ is sufficiently large,
the empirical loss function:
\begin{align*}
    \widehat{\calL}_{\mathrm{bwd}}^{(T)}(E, D, L, H) := \widehat{J}_{\mathrm{bwd}}^{(T)}(E, D, L, H) + \lambda_{\mathrm{fwd}} \cdot \widehat{J}_{\mathrm{fwd}}^{(T)}(E, L, H) + \lambda_{\mathrm{reg}} \cdot R(E)
\end{align*}
satisfies for all $(E, D, L, H)$ (cf.~\Cref{prop:population_and_finite_sample_equivalence}):
\begin{align*}
    \frac{1}{2} \calL_{\mathrm{bwd}}^{(T)}(E, D, L, H) \leq \widehat{\calL}_{\mathrm{bwd}}^{(T)}(E, D, L, H) \leq \frac{3}{2} \calL_{\mathrm{bwd}}^{(T)}(E, D, L, H).
\end{align*}
On this high probability event,
the optimizers $(\widehat{E}_T, \widehat{D}_T, \widehat{L}_T, \widehat{H}_T)$ 
of the empirical loss 
$\widehat{\calL}_{\mathrm{bwd}}^{(T)}$
can readily be shown to inherit the same
guarantees from \Cref{cor:latent_stability_transfer}
as the population minimizer
$(\bar{E}_T, \bar{D}_T, \bar{L}_T, \bar{H}_T)$,
with all constants $C_i$, $i =1,\dots,4$ adjusted by a universal constant.
This provides our final solution to \Cref{prob:latent_stabilization}.

\subsection{Insufficiency Results}
\label{sec:results:insufficiency}

In light of our results in the previous section,
a natural question is whether similar guarantees could be obtained from a decoder-free objective, such as $\calL_\mathrm{fwd}^{(T)}(E, L, H)$, or the single-step backward loss $\mathcal{L}_\mathrm{bwd}^{(1)}(E, D, L, H)$. 
We now present a pair of results that illustrate the insufficiency of both these objective functions for \Cref{prob:latent_stabilization}, 
showing that they admit minimizers which fail to capture the unstable subspace and cannot be used to generate stabilizing output-feedback controllers. 
We first look at the insufficiency of the forward loss.

\begin{restatable}[Insufficiency of multi-step forward loss]{mythm}{fwdcycleinsufficient}
\label{thm:fwd_cycle_insufficient}
Suppose $(A, B)$ is stabilizable (\Cref{assmp:stabilizability}),
the non-stable subspace of $A$ is non-trivial (\Cref{assmp:non_trivial_unstable_subspace}),
and
the regularizer $R$ is admissible (\Cref{assmp:regularizer_assumption}).
Let $m_{\mathrm{stab}} := \abs{\{ \lambda \in \eig(A) \mid \lambda \in \R \}}$
denote the {algebraic multiplicity} of the
real-valued eigenvalues of $A$.
Fix a latent dimension $m\in[n]$,
and suppose \emph{either} of the following conditions hold:
\begin{enumerate}[label=(\alph*)]
    \item There exists $\lambda \in \sigma_+(A)$ such that $\lambda \in \R$, $m \leq n-1$, and either $m$ is even or $m_{\mathrm{stab}} \geq 2$.
    \item There exists $\lambda \in \sigma_+(A)$ such that $\lambda \not\in \R$, $m \leq n-2$, and either $m$ is even or $m_{\mathrm{stab}} \geq 1$.
\end{enumerate}
Then there exists
$\bar{E} \in \R^{m \times n}$,
$\bar{L} \in \R^{m \times m}$, and
$\bar{H} \in \R^{m \times p}$ such that
for every $T \in \N_+$ and $\lambda_{\mathrm{reg}} \geq 0$:
\begin{enumerate}[label=(\roman*)]
    \item The tuple $(\bar{E}, \bar{L}, \bar{H})$ is a global minimizer of the forward loss $\calL^{(T)}_{\mathrm{fwd}}$ (cf.~\Cref{eq:total_forwards_loss}).
    \item The latent pair $(\bar{L}, \bar{H})$ is stabilizable.
    \item For every $K \in \R^{p \times m}$, the matrix $A + B K \bar{E}$ is \emph{not} Schur stable.
\end{enumerate}
\end{restatable}

The construction in the proof of \Cref{thm:fwd_cycle_insufficient} mirrors
the construction described in
\Cref{alg:exact_solution_given_A_B} from
\Cref{sec:problem_formulation:invariance_subspaces}, but makes sure to \emph{omit} at least one $\lambda \in \sigma_+(A)$ from the 
eigenvalues $\Lambda$ which are passed into
the real Schur decomposition.
The hypothesis enumerated in \Cref{thm:fwd_cycle_insufficient} ensure that it is possible to construct an admissible
multi-subset $\Lambda \subset \eig(A)$ such that
(i) $\abs{\Lambda} = m$, 
(ii) $\eig_+(A) \setminus \Lambda$ is non-empty, and 
(iii) $\Lambda$ is conjugate-invariant;
in particular, the latent space may be over-parameterized ($m \geq m_+$) and still miss some non-stable eigenvalues.
We emphasize that the failure is not due to representation collapse, as
the optimal encoder is perfectly conditioned, with $\bar{E} \bar{E}^\T = I_m$.

\Cref{thm:fwd_cycle_insufficient} shows that no sufficiently long horizon length can ensure that the decoder-free objective $\calL_{\mathrm{fwd}}^{(T)}$ always captures the non-stable subspace. 
The backward loss, by contrast, does capture the non-stable subspace, but only once the horizon $T$ is long enough (\Cref{thm:informal_recovery} and \Cref{cor:latent_stability_transfer}). 
We now show that the backward loss horizon cannot be reduced to a single step: at $T=1$, non-normality of the dynamics $A$ can drive the stabilizing gains of the learned latent model and those of the true system apart.
We study this through the following the two-dimensional state, scalar input family of systems (i.e., $n=2$, $p=1$):
\begin{align}
    A_\eta := \bmattwo{\lambda_u}{\eta}{0}{\lambda_s}, \qquad B := \cvectwo{b_1}{b_2}, \qquad \lambda_u > 1 > \abs{\lambda_s}, \quad b_1 \geq 0, b_2 > 0, \quad \eta > 0. \label{eq:nonnormal_example}
\end{align}
By the PBH test, the pair $(A_\eta, B)$ 
in \eqref{eq:nonnormal_example} is stabilizable for any $\eta > 0$.
We will utilize a latent dimension 
of $m = m_+ = 1$ for this system.
Since $m = p = 1$, the sets of gains that stabilize a latent model $(L, H)$ and the true system under output feedback through $E$ are the subsets of $\R$:
\begin{align}
    \calK_L(L, H) := \left\{ K \in \R \mid \abs{L + HK} < 1 \right\}, \qquad \calK_{A_\eta}(E) := \left\{ K \in \R \mid \rho(A_\eta + BKE) < 1 \right\}. \label{eq:gain_sets}
\end{align}

\begin{figure}[t]
    \centering
    \includegraphics[width=\linewidth]{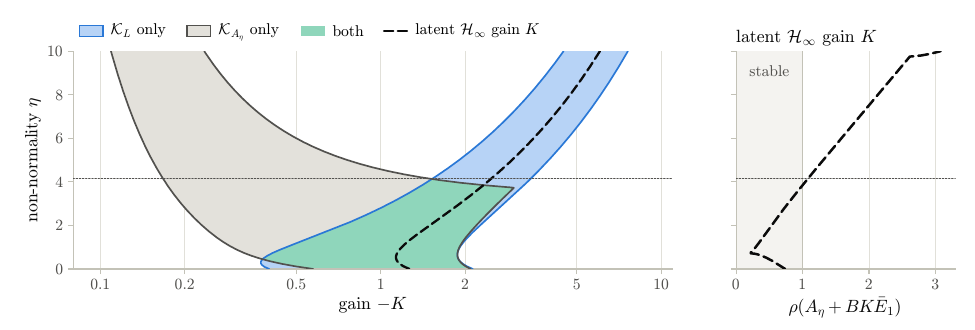}
    \caption{For the single-step backward loss $\calL_\mathrm{bwd}^{(1)}$,
    the set of stabilizing gains for the learned latent system separates from the set of output-feedback stabilizing gains for the LTI system \eqref{eq:nonnormal_example} as non-normality grows.
    In the right panel, the spectral radius $\rho(A_\eta+BK\bar{E}_1)$ of the resulting output-feedback closed-loop matrix is shown.
    }
    \label{fig:insufficiency_K_vs_b}
\end{figure}

\begin{restatable}[Insufficiency of the single-step backward loss for non-normal systems]{mythm}{bwdsinglestepnonnormal}
\label{thm:bwd_single_step_nonnormal}
Suppose the underlying LTI system $(A, B)$ is an instance of \Cref{eq:nonnormal_example},
and that the regularizer $R(E)$ is admissible (\Cref{def:forward_admissible_regularizer}).
There exists an $\eta_0 > 0$ such that for all $\eta \geq \eta_0$, any minimizer
$(\bar{E}_1, \bar{D}_1, \bar{L}_1, \bar{H}_1)$ for the $T=1$ backward loss $\calL_{\mathrm{bwd}}^{(1)}$ in \eqref{eq:total_backwards_loss} satisfies that
the latent and true stabilizing sets are \emph{disjoint}, i.e., $\calK_L(\bar{L}_1, \bar{H}_1) \cap \calK_{A_\eta}(\bar{E}_1) = \emptyset$. 
An explicit sufficient threshold $\eta_0$ is given by:
\begin{align}
    \eta_0 :=1+\max\left\{
        \frac{2(3+\abs{\lambda_u+\lambda_s})\norm{B}}{b_2},\frac{8C}{b_2^2}
        \left(\frac{2\norm{B}^2}{\sigma_x^2}
            +\frac{(2+\abs{\lambda_u+\lambda_s})^2}{\sigma_u^2\alpha_1(0)}\right)
    \right\}, \label{eq:eta_0}
\end{align}
where $C:=\frac{\sigma_x^2}{2}(3+\lambda_u^2+\lambda_s^2)
        +\sigma_u^2\alpha_1(0)b_2^2
        +\lambda_{\mathrm{fwd}}\sigma_x^2
            \max\{\lambda_u^2,1+\lambda_s^2\}$.
\end{restatable}

\Cref{thm:bwd_single_step_nonnormal}
states that for $\eta \geq \eta_0$ 
minimizing the $T=1$ backward loss 
is insufficient to ensure overlapping gain sets \eqref{eq:gain_sets}; this is illustrated numerically in \Cref{fig:insufficiency_K_vs_b}.
On the other hand,
\Cref{cor:latent_stability_transfer}
shows that minimizing loss
$\calL_{\mathrm{bwd}}^{(T)}$ with
$T$ sufficiently large suffices
for the two gain sets in \eqref{eq:gain_sets} to have a non-trivial intersection (which can be computed by minimizing the $H_\infty$-norm).
Specifically, we specialize the full theorem (\Cref{thm:end_to_end_closed_loop}) to
the system \eqref{eq:nonnormal_example}, with latent dimension $m=1$.
To obtain fully explicit bounds, take the input variance schedule
$\alpha_T(k) \equiv 1/T$,
let $(\sigma_0,\sigma_1,\e_0)$ be the admissibility parameters of $R$, and define:
\begin{align*}
    \Delta&:=\lambda_u-\lambda_s, \qquad
    \gamma:=\frac{1+\abs{\lambda_s}}{2}, \qquad
    \chi_\eta:=\sqrt{1+\frac{\eta^2}{\Delta^2}}, \\
    k_\eta&:=\frac{\lambda_u\sqrt{\Delta^2+\eta^2}}
        {b_1\Delta+b_2\eta}, \qquad
    g_\eta:=1+\frac{\lambda_u+\eta+\norm{B}k_\eta}{1-\abs{\lambda_s}}.
\end{align*}

\begin{figure}[t]
    \centering
    \includegraphics[width=\linewidth]{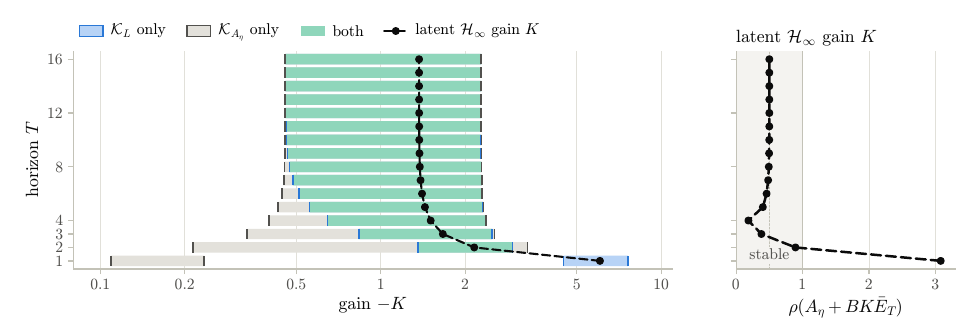}
    \caption{For the multi-step backward loss $\calL_\mathrm{bwd}^{(T)}$,
    the set of stabilizing gains for the learned latent system converges to the set of 
    output-feedback stabilizing gains for LTI system \eqref{eq:nonnormal_example} as $T$ grows (here the non-normality parameter $\eta$ is held fixed).
    In the right panel, the spectral radius $\rho(A_\eta+BK\bar{E}_T)$ of the resulting output-feedback closed-loop matrix is shown, indicating stable closed-loop behavior for $T \geq 2$.
    }
    \label{fig:insufficiency_K_vs_T}
\end{figure}

A tedious but elementary computation (which we omit)
involving the conditions listed in \Cref{thm:end_to_end_closed_loop}
shows that the following are sufficient choices of
hyperparameters for $\calL_{\mathrm{bwd}}^{(T)}$:
\begin{align*}
    T&\geq
        \frac{128\sigma_1g_\eta^2\chi_\eta^4
            (\sigma_x^2+\sigma_u^2\norm{B}^2)}
            {\sigma_0^2\sigma_x^2(1-\gamma^2)}
        (1+k_\eta+g_\eta\norm{B}k_\eta)^2, \\
    \lambda_{\mathrm{fwd}}
        &\geq\frac{64\sigma_1g_\eta^2\chi_\eta^4
            (\sigma_x^2+\sigma_u^2\norm{B}^2)}
            {\sigma_0^2(1-\gamma^2)}
        (1+k_\eta+g_\eta\norm{B}k_\eta)^2
        \max\left\{\frac{1}{\sigma_x^2},\frac{1}{\sigma_u^2}\right\},
        \\
    \lambda_{\mathrm{reg}}
        &\geq\frac{\chi_\eta^2(\sigma_x^2+\sigma_u^2\norm{B}^2)}
            {\e_0T(1-\gamma^2)}.
\end{align*}
For any resulting minimizer $(\bar{E}_T,\bar{D}_T,\bar{L}_T,\bar{H}_T)$,
a minimum $H_\infty$-norm latent controller over the gain ball
$\abs{K}\leq M_0=k_\eta/\sqrt{\sigma_0}$ exists, and every such controller
belongs to
$\calK_L(\bar{L}_T,\bar{H}_T)\cap\calK_{A_\eta}(\bar{E}_T)$.
We note that treating all parameters
other than $\eta$ as fixed,
the scaling of $T$, $\lambda_{\mathrm{fwd}}$, $\lambda_{\mathrm{reg}}$ are:
\begin{align*}
    T\geq \Omega(\eta^8), \qquad
    \lambda_{\mathrm{fwd}} \geq \Omega(\eta^8), \qquad
    \lambda_{\mathrm{reg}} \geq \Omega\left(\frac{\eta^2}{T}\right).
\end{align*}

\Cref{fig:insufficiency_K_vs_T} illustrates the two gain sets in \eqref{eq:gain_sets} converging as $T$ grows, for a fixed non-normality parameter $\eta$. Observe that 
the sufficient conditions computed above for 
the multi-step backward loss from
\Cref{thm:end_to_end_closed_loop} are 
conservative on this example, as closed-loop stability is achieved as early as $T = 2$;
we leave investigating a sharp characterization of the \emph{minimal} 
horizon $T$ necessary to future work.

\section{Experiments}
\label{sec:experiments}

We conduct a comprehensive experimental evaluation of optimizing both forward and backward losses for learning a latent dynamics model of the flexible arm described in \Cref{sec:introduction:example}.
We first describe the non-linear dynamics of the arm in \Cref{sec:experiments:arm}.
Then, in \Cref{sec:experiments:linear}
we consider learning linear latent models
for the linearized arm dynamics around its upright equilibrium, 
which matches the setting of our theory
in \Cref{sec:results}.
Next, in \Cref{sec:experiments:nonlinear}, we 
study the more practical setting where
we learn non-linear latent MLP models 
for the non-linear arm dynamics.
Finally, we conclude in \Cref{sec:experiments:multi} with an investigation into some of the optimization issues which arise
when training long-horizon multi-step losses.

\paragraph{Methodology.}
We collect training data by sampling initial conditions $x_0 \sim \sfN(0, I_n)$
and inputs $u_t \sim \sfN(0, \sigma_u^2/T)$, $\sigma_u = 0.5$
(cf.~\Cref{assmp:data_model}).
For each setup, we collect $5$ data sets, each with $N$ independent trajectories
of length $T$, and refer to each dataset as $\mathcal{D}_s$ for $s\in\{1,\dots,5\}$.
For experiments that vary $N$ we consider subsets of $\mathcal{D}_s$, whereas inputs are redrawn with variance $\sigma_u^2/T$ when varying $T$.
For the backward loss $\calL_{\mathrm{bwd}}^{(T)}$ (cf.~\eqref{eq:total_backwards_loss}),
we set $\lambda_{\mathrm{fwd}} = 1$,
and furthermore replace $J_{\mathrm{fwd}}^{(T)}$ by the one-step term $J_{\mathrm{fwd}}^{(1)}$ for linear models, and by the ``teacher-forcing'' loss
$\frac{1}{T}\sum_{t=0}^{T-1}\E\norm{E(x_{t+1}) - f_z(E(x_t), u_t)}^2$
for non-linear models. 
We find that these choices exhibit faster and more stable optimization compared with $J_{\mathrm{fwd}}^{(T)}$.
Linear models are trained with full-batch Adam followed by L-BFGS, 
and MLPs are trained via mini-batch Adam.
Since the optimization landscape changes with $T$, we solve the horizons in increasing order using the following protocol: at each $T$, 
we train from both random initializations and from warm-starting
with the solution at the previous horizon, and report, per data set, 
the model with the lowest objective (except for \Cref{sec:experiments:multi}, 
where we study the initialization scheme used and train models for each $T$ independently).
We say that a linear model $E$ \emph{captures} the non-stable modes of $A$ if $\Delta := \Delta_F(\calG^\R_+(A^\T); \mathrm{row}(E)) < 0.5$;
the specific threshold value is not critical as we find $\Delta$ to be bimodal in our runs.
\Cref{sec:appendix:experiments} contains a more detailed 
explanation of the experimental setup.

\subsection{Flexible Arm Model}
\label{sec:experiments:arm}

\begin{figure}[t]
    \centering
    \begin{tikzpicture}[
    x=1.45cm,y=1.45cm,font=\normalsize,
    >={Latex[length=2mm]},
    bar/.style={line width=3pt,draw=black!55,line cap=round},
    neighbor/.style={bar,draw=black!35,opacity=0.55},
    guide/.style={densely dashed,draw=black!40},
    force/.style={->,thick,draw=RoyalBlue!80!black},
    moment/.style={->,thick,draw=BrickRed!85!black},
    angle/.style={->,thin,draw=black!75}
]

\coordinate (J) at (0,0);
\coordinate (Jp) at (2.048,1.434);
\draw[neighbor] (J) -- ++(240:1.25);
\draw[bar] (J) -- (Jp);
\draw[neighbor] (Jp) -- ++(55:0.8);
\draw[guide] (J) -- ++(90:1.65);
\draw[guide] (J) -- ++(60:1.45);
\draw[angle] (90:1.12) arc[start angle=90,end angle=35,radius=1.12];
\node[fill=white,inner sep=1pt] at (72:1.32) {$\phi_i$};
\draw[angle] (60:0.68) arc[start angle=60,end angle=35,radius=0.68];
\node[fill=white,inner sep=1pt] at (49:0.88) {$q_i$};
\foreach \joint in {J,Jp}{
    \fill[black] (\joint) circle[radius=0.085];
}
\node[anchor=north east] at ($(J)+(-0.22,-0.15)$) {joint $i$};
\node[anchor=west] at ($(Jp)+(0.27,0)$) {joint $i+1$};

\coordinate (G) at ($(J)!0.5!(Jp)$);
\fill (G) circle[radius=0.045];
\draw[force] (G) -- ++(0,-1.5) node[below] {$\bar m g$};
\draw[thin,<->] ($(J)+(125:0.95)$) -- ($(Jp)+(125:0.95)$)
    node[midway,sloped,above] {$\ell$};
\draw[moment] (210:0.48) arc[start angle=210,end angle=-30,radius=0.48];
\node[anchor=east,text=BrickRed!85!black] at (-0.58,0.25) {$\tau_i$};
\node[anchor=west] at (3.8,1.05)
    {$\tau_i=\underbrace{u}_{\text{input}}-\underbrace{k_iq_i}_{\text{spring}}
      -\underbrace{(C_d\dot q)_i}_{\text{damping}}$};
\node[anchor=west] at (3.8,0.15) {$q_i=\phi_i-\phi_{i-1}$};
\end{tikzpicture}
    \caption{Local geometry of the flexible arm described in \eqref{eq:arm_dynamics}. Section $i$, of length $\ell$, connects joints $i$ and $i+1$, and the gravitation force $\bar m g$ acts at its center of mass. 
    The solid dots indicate torsional springs of stiffness $k_i:=\kappa\cdot (K_0)_{ii}$. 
    The torque $\tau_i$ combines the input, spring, and damping contributions at joint $i$, with $C_d:=c_0M+c_1K_{\mathrm{eff}}$. 
    }
    \label{fig:arm_free_body}
\end{figure}
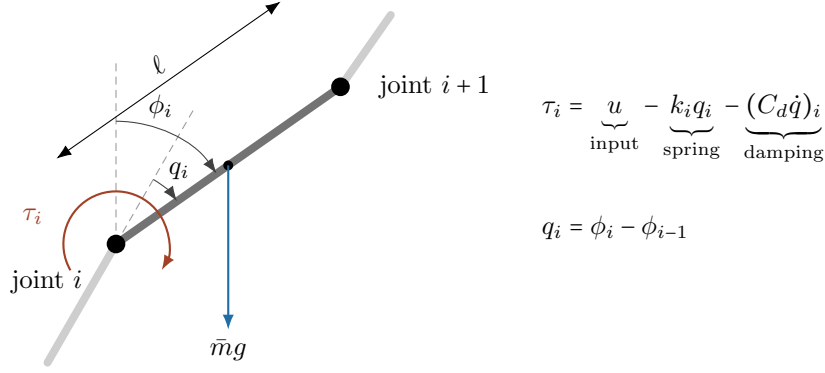

We consider a Hencky-type bar-chain model of a flexible arm~\cite{baroudi2019nonlinear,giorgio2018nonlinear}.
The arm is a planar chain of $n_s = 20$ rigid sections of length $\ell = 0.05$\,m and mass $\bar{m} = 0.04$\,kg, joined by torsional springs.
Let $q \in \R^{n_s}$ be the joint angles, where $q = 0$ corresponds to the upright configuration, so that the section angles are $\phi = Sq$ for the lower-triangular ones matrix $S$.
\Cref{fig:arm_free_body} illustrates the local geometry and the forces and torques on one section.
Denoting the configuration-dependent inertia as $M(q)\succ0$ and its centrifugal terms as $C(q, \dot q)\dot q$, the dynamics of the flexible arm model are:
\begin{align}
    M(q)\ddot{q} + C(q, \dot q)\dot q + (c_0M(q) + c_1K_{\mathrm{eff}})\dot{q} + \kappa K_0q - S^\T\left(w \odot \sin(Sq)\right) = \bm{1}u, \qquad q \in \R^{n_s}, \,\, u \in \R, \label{eq:arm_dynamics}
\end{align}
where $w_i := (n_s - i + \tfrac{1}{2})g\ell\bar{m}$ is the gravitational potential coefficient on joint $i$, 
$K_{\mathrm{eff}} := \kappa K_0 - S^\T\diag(w)S$
where $K_0$ is a diagonal matrix with entries
denoting the relative stiffness of the individual torsional springs, spaced uniformly in $[1, 1.5]$, 
$\kappa$ is a stiffness scale described shortly,
and $(c_0, c_1) := (0.4, 0.02)$ are the Rayleigh damping coefficients. 
The stiffness scale $\kappa$ lies midway between the two largest generalized eigenvalues of $(S^\T\diag(w)S, K_0)$ so that $K_{\mathrm{eff}}$ has 
exactly one negative eigenvalue, and hence gravity buckles (destabilizes) only the lowest bending mode.
The scalar input $u$ applies the same torque at every joint.
For $x = (q, \dot{q}) \in \R^{2n_s}$, the zero-order-hold linearization of \eqref{eq:arm_dynamics} at $q=0$ with $\Delta t = 0.05$\,s gives $(A, B)$ with $n = 40$, $p = 1$, 
and $\sigma_+(A) = \{1.129\}$, so $m_+ = 1$; 
the slowest stable eigenvalues are $0.874$ and $0.581 \pm 0.633j$, and all others have modulus less than or equal to $0.254$.
We measure the input in units of $2.14^{-1}$\,N$\cdot$m, which normalizes $\norm{B} = 1$. 
The input variance $\sigma_u^2$ and the LQR weight $R = I_p$ refer to these units. %
The non-linear arm is simulated with an exponential integrator with $50$ sub-steps per $\Delta t$.

\subsection{Learning Linear Latent Models for the Linearized Arm Dynamics}
\label{sec:experiments:linear}

\begin{figure}[ht!]
    \centering
    \includegraphics[width=\linewidth]{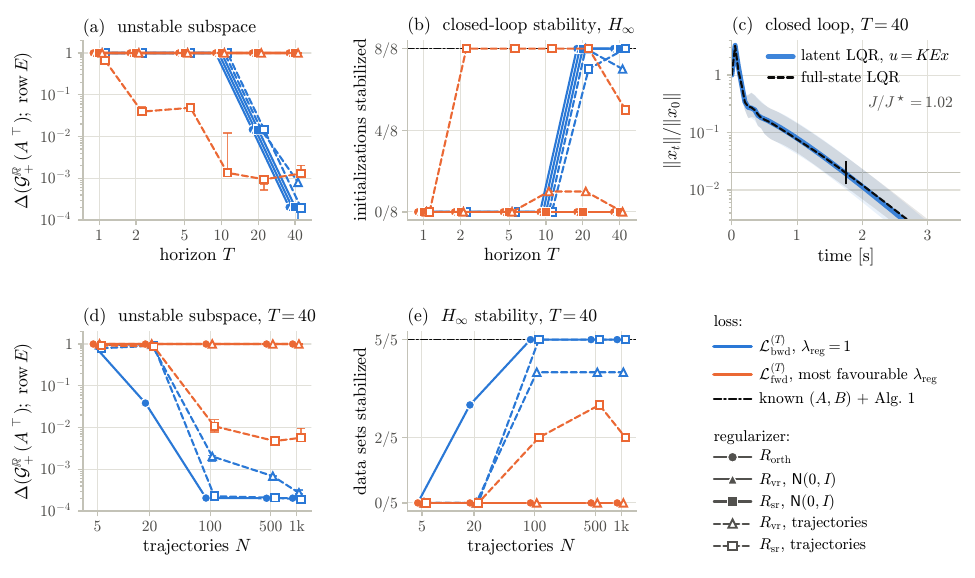}
    \caption{
    Linearized flexible arm controlled by linear latent models ($m = 4$). Panels (a, b) 
    show results from training on the population losses $\calL_{\mathrm{bwd}}^{(T)}$ 
    (blue, $\lambda_{\mathrm{reg}} = 1$) and $\calL_{\mathrm{fwd}}^{(T)}$ (orange, most favorable $\lambda_{\mathrm{reg}}$). 
    (a) Subspace distance $\Delta$ of the minimizer with the lowest objective among $8$ random initializations 
    (error bars span the $25$\%-$75$\% interquartile range over the $8$ initializations). 
    (b) Number of the $8$ initializations whose latent $H_\infty$ gain $K$ makes $A + BKE$ Schur stable.
    (c) 
    The ratio $\norm{x_t}/\norm{x_0}$
    for both the closed-loop of the latent LQR gain (cost $J$) learned from $\calL_{\mathrm{bwd}}^{(40)}$ with $R_{\mathrm{orth}}$, and of full state LQR (cost $J^\star$).
    The lines and shaded area illustrate the median and $25$\%-$75$\% interquartile range over $32$ initial states $x_0 \sim \sfN(0, I_n)$.
    The panels (d, e) are as (a, b), but for the empirical losses of $N$ trajectories at $T = 40$. 
    The black dash-dotted line in (b, e) signifies \Cref{alg:exact_solution_given_A_B} with known $(A, B)$, 
    which stabilizes at every $T$ and $N$. 
    Note that series are slightly offset horizontally for legibility.
    }
    \label{fig:exp_arm_recovery}
\end{figure}

We first collect training trajectories from
the linearization $(A, B)$ of the flexible arm \eqref{eq:arm_dynamics}
about its upright equilibrium, and learn linear latent dynamics and autoencoders, 
matching the setting of our 
theoretical results in \Cref{sec:results:multi_step_backward}.
We minimize $\calL_{\mathrm{bwd}}^{(T)}$ and $\calL_{\mathrm{fwd}}^{(T)}$ 
with $m = 4$ and $T \in \{1, 2, 5, 10, 20, 40\}$, 
both in its population form 
from $8$ random initializations in panels (a, b, c) of \Cref{fig:exp_arm_recovery},
or from up to $N=1000$ trajectories (cf.~\Cref{rmk:empirical_losses}) in panels (d, e) of \Cref{fig:exp_arm_recovery}, with two initializations per data set.
In \Cref{fig:exp_arm_recovery}, for the population case (panel b), we maintain $8$ streams of models (one for each initialization) across $T$; at each $T > 1$, 
we keep the lowest objective of two models (the reinitialization trained from scratch, or the warm-start on that stream's previous model).
We compare five regularizers $R$:
$R_{\mathrm{orth}}$, $R_{\mathrm{sr}}$, and $R_{\mathrm{vr}}$ ($\beta = 1$) with base measure $\sfN(0, I_n)$, 
which are admissible (cf.~\Cref{prop:regularizer_admissibility}), and 
$R_{\mathrm{sr}}$ and $R_{\mathrm{vr}}$ with base measure set to the empirical measure over all trajectory states $x_0, \dots, x_T$, 
as in~\cite{balestriero2025lejepa,maes2026leworldmodel}.
For $\calL_{\mathrm{bwd}}^{(T)}$ we fix $\lambda_{\mathrm{reg}} = 1$, 
and for $\calL_{\mathrm{fwd}}^{(T)}$ we use the $\lambda_{\mathrm{reg}} \in \{1, 10, 100, 1000\}$ that maximizes its stabilization rate (i.e., the value of $\lambda_\mathrm{reg}$ that produces the maximum number of stabilizing controllers, counted over all $T$).
From each learned latent dynamics $(L, H)$, we compute both a latent LQR gain ($Q = I_m$, $R = I_p$), and a latent $H_\infty$ gain which minimizes the $H_\infty$-norm from $w$ to $(z, u)$ for $z_{t+1} = (L + HK)z_t + w_t$. 
Here, the $z$-block is the objective of \Cref{thm:resolvent_gain_bound}, and the $u$-block 
minimization, which we apply with $u = KEx$,
replaces the 
constraint $K \in B^{\mathrm{(op)}}_{p \times m}(M_0)$.
We report $\Delta$ and whether the $H_\infty$ gain stabilizes $A + BKE$ in \Cref{fig:exp_arm_recovery}. 
The LQR-gain gives the same outcomes, see \Cref{tab:exp_gains}.
In \Cref{tab:exp_conjugacy}, we report the success rate 
on the non-linear arm from $8$ initial states (each given a random bend, no initial velocity) with the tip 
displaced by $2$ to $80$\,cm along the buckling mode. %
The reference is \Cref{alg:exact_solution_given_A_B} with known $(A, B)$, retaining the four slowest modes.

\begin{table}[ht]
\centering\footnotesize
\makeatletter\providecommand{\tabinput}[1]{\input{#1} }\makeatother
\begin{tabular}{lllcccccc}
\toprule
 & & & \multicolumn{6}{c}{population: \# stable LQR gains (of $8$), \# stable $H_\infty$ gains (of $8$), $\Delta$} \\
\cmidrule(l){4-9}
loss & regularizer & $\lambda_{\mathrm{reg}}$ & $T = 1$ & $2$ & $5$ & $10$ & $20$ & $40$ \\
\midrule
$\mathcal{L}_{\mathrm{bwd}}^{(T)}$ & $R_{\mathrm{orth}}$ & $1$ & $0, 0, 1$ & $0, 0, 1$ & $0, 0, 1$ & $0, 0, 1$ & $8, 8, 0.015$ & $8, 8, 2.0\cdot 10^{-4}$ \\
 & $R_{\mathrm{vr}}$, $\mathsf{N}(0,I)$ & $1$ & $0, 0, 1$ & $0, 0, 1$ & $0, 0, 1$ & $0, 0, 1$ & $8, 8, 0.015$ & $8, 8, 2.0\cdot 10^{-4}$ \\
 & $R_{\mathrm{sr}}$, $\mathsf{N}(0,I)$ & $1$ & $0, 0, 1$ & $0, 0, 1$ & $0, 0, 1$ & $0, 0, 1$ & $8, 8, 0.015$ & $8, 8, 2.1\cdot 10^{-4}$ \\
 & $R_{\mathrm{vr}}$, traj. & $1$ & $0, 0, 1$ & $0, 0, 1$ & $0, 0, 1$ & $0, 0, 1$ & $8, 8, 0.015$ & $8, 7, 8.0\cdot 10^{-4}$ \\
 & $R_{\mathrm{sr}}$, traj. & $1$ & $0, 0, 1$ & $0, 0, 1$ & $0, 0, 1$ & $0, 0, 1$ & $7, 7, 0.015$ & $8, 8, 2.0\cdot 10^{-4}$ \\
\midrule
$\mathcal{L}_{\mathrm{fwd}}^{(T)}$ & $R_{\mathrm{orth}}$ & $1$ & $0, 0, 1$ & $0, 0, 1$ & $0, 0, 1$ & $0, 0, 1$ & $0, 0, 1$ & $0, 0, 1$ \\
 & $R_{\mathrm{vr}}$, $\mathsf{N}(0,I)$ & $10$ & $0, 0, 1$ & $0, 0, 1$ & $0, 0, 1$ & $0, 0, 1$ & $0, 0, 1$ & $0, 0, 1$ \\
 & $R_{\mathrm{sr}}$, $\mathsf{N}(0,I)$ & $100$ & $0, 0, 1$ & $0, 0, 1$ & $0, 0, 1$ & $0, 0, 1$ & $0, 0, 1$ & $0, 0, 1$ \\
 & $R_{\mathrm{vr}}$, traj. & $10$ & $0, 0, 1$ & $0, 0, 1$ & $0, 0, 1$ & $1, 1, 1$ & $1, 1, 1$ & $0, 0, 1$ \\
 & $R_{\mathrm{sr}}$, traj. & $100$ & $0, 0, 0.66$ & $8, 8, 0.039$ & $8, 8, 0.049$ & $8, 8, 1.4\cdot 10^{-3}$ & $8, 8, 9.4\cdot 10^{-4}$ & $6, 5, 1.3\cdot 10^{-3}$ \\
\midrule[\heavyrulewidth]
 & & & \multicolumn{5}{c}{finite sample: \# stable LQR gains (of $5$), \# stable $H_\infty$ gains (of $5$), $\Delta$} & \\
\cmidrule(l){4-8}
 & & & $N = 5$ & $20$ & $100$ & $500$ & $1000$ &  \\
\midrule
$\widehat{\mathcal{L}}_{\mathrm{bwd}}^{(T)}$ & $R_{\mathrm{orth}}$ & $1$ & $0, 0, 0.96$ & $2, 3, 0.039$ & $5, 5, 2.1\cdot 10^{-4}$ & $5, 5, 2.0\cdot 10^{-4}$ & $5, 5, 2.0\cdot 10^{-4}$ &  \\
 & $R_{\mathrm{vr}}$, traj. & $1$ & $0, 0, 0.96$ & $0, 0, 0.97$ & $5, 4, 2.0\cdot 10^{-3}$ & $5, 4, 6.9\cdot 10^{-4}$ & $5, 4, 2.7\cdot 10^{-4}$ &  \\
 & $R_{\mathrm{sr}}$, traj. & $1$ & $0, 0, 0.79$ & $0, 0, 0.91$ & $5, 5, 2.3\cdot 10^{-4}$ & $5, 5, 2.1\cdot 10^{-4}$ & $5, 5, 1.9\cdot 10^{-4}$ &  \\
\midrule
$\widehat{\mathcal{L}}_{\mathrm{fwd}}^{(T)}$ & $R_{\mathrm{orth}}$ & $1$ & $0, 0, 1$ & $0, 0, 1$ & $0, 0, 1$ & $0, 0, 1$ & $0, 0, 1$ &  \\
 & $R_{\mathrm{vr}}$, traj. & $10$ & $0, 0, 1$ & $0, 0, 1$ & $0, 0, 1$ & $0, 0, 1$ & $0, 0, 1$ &  \\
 & $R_{\mathrm{sr}}$, traj. & $100$ & $0, 0, 0.92$ & $0, 0, 0.88$ & $2, 2, 0.011$ & $2, 3, 4.7\cdot 10^{-3}$ & $4, 2, 5.7\cdot 10^{-3}$ &  \\
\bottomrule
\end{tabular}
 
\caption{
    A comparison of stabilization performance for latent $H_\infty$ versus latent LQR control design,
    for both population losses (top, with $T$ varying)
    and finite-sample losses (bottom, with $T=40$ fixed, $N$ varying).
    The subspace distance $\Delta$ shown is for the 
    lowest-objective model over $8$ initialization (same as \Cref{fig:exp_arm_recovery}(a)) for the population loss, 
    and the median over $5$ datasets (each the better of $2$ initializations) for the finite-sample loss.
    Metrics for $\calL_{\mathrm{fwd}}^{(T)}$ and $\widehat{\calL}_{\mathrm{fwd}}^{(T)}$ 
    are shown for the value of $\lambda_\mathrm{reg}$ that produces the largest number of stabilizing models over all $T$.
}
\label{tab:exp_gains}%
\end{table}

\paragraph{Population minimizers of the backward loss.}
The results of \Cref{fig:exp_arm_recovery} empirically corroborate \Cref{thm:informal_recovery} and \Cref{thm:resolvent_gain_bound}. 
Indeed, we observe that
as we increase the reconstruction horizon $T$, 
the row space of the encoder trained using the population backward loss aligns with the unstable subspace for all choices of regularizer:
by $T=20$, the subspace distance is $\Delta \approx 0.015$, and by $T=40$ we have $\Delta \approx 2\cdot10^{-4}$ for all regularizers except for trajectory-VICReg, 
inducing closed-loop stability in 39 out of 40 
($5$ regularizers $\times$ $8$ initializations per regularizer) of the trained models (see \Cref{tab:exp_gains}). 
Of all the regularizers, trajectory-measure VICReg exhibits the largest number of stabilization failures at $T=40$ when regularizing $\mathcal{L}_\mathrm{bwd}^{(T)}$,
in both the population and finite-sample cases. We attribute this to a lack of admissibility (cf.~\Cref{def:forward_admissible_regularizer}), as the corresponding encoder's 
minimum squared singular value is significantly smaller than the other regularizers ($0.0015$ vs.\ $\ge0.05$ for admissible regularizers, cf.~\Cref{tab:exp_conjugacy} in the appendix).
Despite \Cref{tab:exp_conjugacy} indicating that the encoder and latent dynamics capture the unstable mode 
(subspace distance: $\Delta \approx 8\cdot10^{-4}$, unstable-mode conjugacy: $r_u \approx 3.9\cdot10^{-3}$, \Cref{tab:exp_gains,tab:exp_conjugacy} respectively), 
the forward conjugacy of the entire latent model is almost two orders of magnitude larger
than for the admissible regularizers ($\approx 2$ vs.\ $\approx 0.035$)  when computed in a scale-free manner
(cf.~\Cref{tab:exp_conjugacy}).
This suggests that trajectory-VICReg's stabilization failures are the result of the latent dynamics not
capturing invariance, 
likely due to the poor conditioning of the encoder impacting optimization, 
and emphasizing the need for \Cref{assmp:regularizer_assumption}. 

\paragraph{Population minimizers of the forward loss.}
\Cref{fig:exp_arm_recovery}'s results for the population forward loss $\mathcal{L}^{(T)}_{\mathrm{fwd}}$ support the conclusion of \Cref{thm:fwd_cycle_insufficient}:
we observe that despite
low forward conjugacy residual $r_A$ ($\le 6\cdot10^{-5}$ with admissible regularizers, see \Cref{tab:exp_conjugacy}), 
nearly all trained models fail to capture the unstable subspace for any value of $T$ and 
lead to unstable closed-loop behavior. The only exceptions are trajectory-measure SIGReg with $8/8$ for $2\le T\le 20$ and $5/8$ at $T=40$,
and trajectory-measure VICReg for only one initialization between $10\le T\le 20$. 
We emphasize that 
\Cref{thm:fwd_cycle_insufficient}
shows the \emph{existence} of a global minimizer that does not capture the unstable dynamics: this does not imply that every global minimizer 
fails to capture the unstable dynamics.
Indeed, it is in principle possible
to converge to stabilizing minimizers
due to the randomness of the
training process.
However, as seen in 
\Cref{fig:exp_arm_recovery}
and \Cref{tab:exp_gains},
the emergence of closed-loop stability
for the forward loss is not robust
to either $\lambda_{\mathrm{reg}}$
and regularizer choices.

\paragraph{Finite-sample agreement.}
The minimizers of the empirical losses $\widehat{\mathcal{L}}_\mathrm{bwd}^{(T)}$ and $\widehat{\mathcal{L}}_\mathrm{fwd}^{(T)}$ (see \Cref{rmk:empirical_losses}) 
exhibit the same trends as their population counterparts as the dataset size $N$ grows. 
Panel (d) of \Cref{fig:exp_arm_recovery} shows that the minimizing encoders of the empirical backward loss
exhibit similar alignment to the population-minimizing encoders for $N\ge 100$. 
For the empirical forward loss, trajectory-measure SIGReg is again the only regularizer to achieve closed-loop stability for some datasets, 
though it does not achieve comparable $H_\infty$ stabilization performance for any of the tested dataset sizes.

\paragraph{$H_\infty$ vs.\ LQR performance.}
We observe that latent LQR exhibits comparable stabilization performance to latent $H_\infty$ control. 
\Cref{tab:exp_gains} compares the number of stabilizing $H_\infty$ and LQR gains for both losses, and in the population and finite-sample cases.
Furthermore, the cost of latent LQR is $2\%$ above full-state LQR (\Cref{fig:exp_arm_recovery}c). 
This suggests that a similar analysis to \Cref{thm:resolvent_gain_bound} 
for a latent LQR control design procedure should be possible,
a task we leave for future work.

\begin{figure}[ht!]
    \centering
    \includegraphics[width=\linewidth]{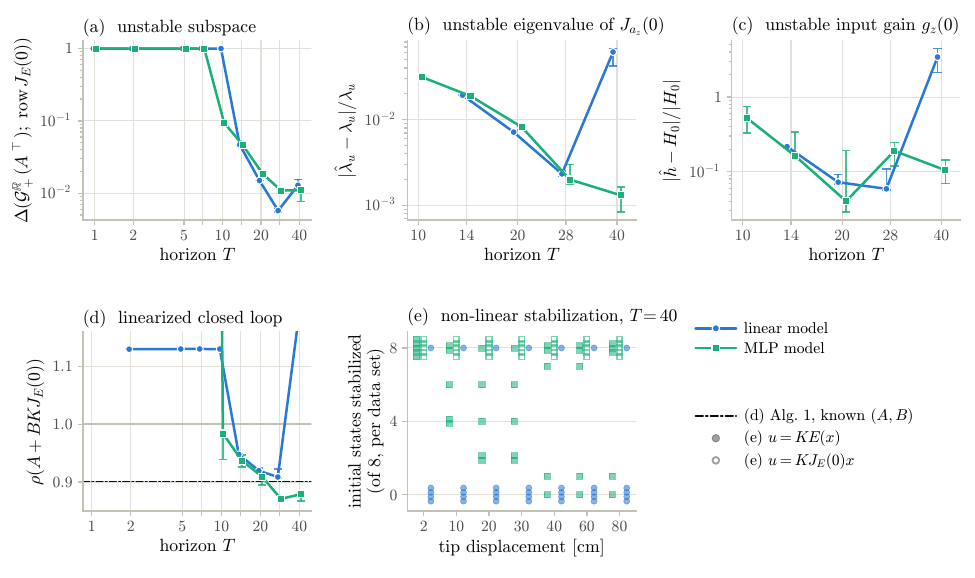}
    \caption{
    MLP (green) and linear (blue) latent models trained on data from the non-linear arm. 
    Solid lines are medians, and error bars are the $25\%$ to $75\%$ interquartile range.
    All models are trained using the empirical backwards objective $\widehat{\calL}_{\mathrm{bwd}}^{(T)}$ with latent dimension $m=4$ 
    (with the modifications to
    $\widehat{J}_{\mathrm{fwd}}^{(1)}$ discussed in the beginning of 
    \Cref{sec:experiments}).
    Panels (a)--(d) compare the linearization of the MLP model about the unstable equilibrium $q=0$ to the linear model.
    In (a), we compare the subspace distance of $\mathrm{row}(J_E(0))$,
    where $J_E(x) := \frac{\rmd}{\rmd x} E(x)$.
    Panels (b) and (c) compare relative errors of the unstable eigenvalue of $J_{a_z}(0)$ and induced input gain, 
    against $L_0$ and $H_0$ (shown only for models that capture the unstable mode, i.e. $\Delta < 0.5$).
    For panel (c), let $w$ be the left eigenvector of $J_{a_z}(0)$ for the unstable eigenvalue, scaled so that $w^\T J_E(0)D_0=1$.
    Thus, $\hat{h}:=w^\T g_z(0)$ is invariant to the learned latent representation, and equals $H_0$ when the unstable mode is recovered exactly (see \Cref{sec:appendix:unstable_mode}
    for more details).
    Panel (d) displays $\rho(A + BKJ_E(0))$ for the latent LQR gain $K$, and \Cref{alg:exact_solution_given_A_B} with known $(A, B)$ (dash-dotted).
    At $T=1$, the linear model is degenerate ($E\approx 0$) and has no stabilizing latent LQR gain. 
    MLP values for $T\le 7$ lie above the axis.
    For the models selected at $T = 40$, (e) displays for each data set (dot) the number of the $8$ initial states from which its controller stabilizes the non-linear arm, under $u = KE(x)$ (filled) and $u = KJ_E(0)x$ (hollow; the same law for the linear model). Data sets with the same count are stacked vertically.
    Each initial state releases the arm from rest with its tip displaced along the buckling mode (random sign), plus a random bend that leaves the tip in place. The same $8$ initial states are scaled to each tip displacement.
    }
    \label{fig:exp_arm_nonlinear}
\end{figure}

\subsection{Linear and Non-Linear Latent Models for the Non-Linear Arm}
\label{sec:experiments:nonlinear}

We now consider learning latent representations
from trajectories drawn from the non-linear arm
dynamics \eqref{eq:arm_dynamics}.
Here, we will consider learning both 
linear autoencoders and latent dynamics $(E, D, L, H)$, in addition to 
non-linear autoencoders $E(x), D(z)$ with 
control-affine latent dynamics $z_{t+1} = f_z(z_t, u_t) = a_z(z_t) + g_z(z_t)u_t$, 
all parameterized by multi-layer perceptron (MLP) networks with two hidden layers of width $64$, and $\tanh$ activations. 
The $E$, $D$, and $a_z$ networks are anchored so that $E(0) = 0$, $D(0) = 0$, and $a_z(0) = 0$, and $J_{a_z}(0)$ is initialized to a random orthogonal matrix. 
We let $m = 4$ and choose $R(E) = R_{\mathrm{orth}}(E; \hat{\mu})$ for the empirical distribution $\hat{\mu}$ of training states, 
along with $\lambda_{\mathrm{fwd}} = \lambda_{\mathrm{reg}} = 1$.
The dataset is $N = 2000$ trajectories of
the non-linear arm. %
Initial states and inputs 
are scaled by $\varepsilon = 0.01$, 
i.e., $(q_0, \dot q_0) \sim \sfN(0, \varepsilon^2 I_{40})$ (joint angles and rates of standard deviation $0.57^\circ$ and $0.57^\circ$/s) 
and $u_t \sim \sfN(0, \varepsilon^2\sigma_u^2/T)$.
For the rollouts, 
non-linearity is mild for $T \leq 20$, whereas by $T = 40$ the unstable mode saturates 
(qualitatively, $T=40$ corresponds to the arm falling over in about half of the trajectories).
At each horizon $T \in \{1, 2, 5, 7, 10, 14, 20, 28, 40\}$, every model is trained 
on the empirical backward loss
$\widehat{\calL}_{\mathrm{bwd}}^{(T)}$
from $8$ initializations for each of the $5$ data sets.
The MLP models are trained using mini-batches of $128$ trajectories. 
For linear control design on the non-linear latent dynamics, we compute a latent LQR gain $K$ from the linearization $(J_{a_z}(0), g_z(0))$ at the equilibrium. 
We leave an exploration of non-linear control design (e.g., iLQR+MPC)
on the latent model to future work.

\paragraph{Linear latent models of non-linear dynamics.} 
We observe that while fitting a linear autoencoder and latent dynamics to sufficiently non-linear 
data does not impact the encoder's ability to capture the unstable mode, 
it breaks forward conjugacy by forcing the linear latent dynamics to approximate non-linear dynamics, inducing latent 
controllers that do not 
stabilize the true system nor even its linearization.
In \Cref{fig:exp_arm_nonlinear}, panels (a, b, c), it can be seen that while the model initially captures the unstable subspace 
and its unstable dynamics, the accuracy of the latent dynamics model begins to decrease
between $T=28$ and $40$ (cf.\ panels b,c). 
At $T=40$, $E$ still captures the unstable subspace ($\Delta \approx 0.013$), but the latent dynamics are no longer forward conjugate.
\Cref{fig:exp_arm_nonlinear} shows that the relative error of the unstable eigenvalue
grows from $\approx 7.2\cdot 10^{-3}$ at
$T=20$ to $\approx 6.1\cdot10^{-2}$ at $T=40$. 
Concurrently, the linearized closed-loop behavior becomes unstable, and the induced 
controller fails to stabilize the non-linear system (\Cref{fig:exp_arm_nonlinear}, panels d, e, and \Cref{tab:exp_eps_tips}). 
The same pattern emerges at lower $T$ if we increase non-linearity by modulating the data scale $\varepsilon$.
As shown in \Cref{tab:exp_conjugacy_nonlinear}, while the linear model learned from 
non-linear data stabilizes the linearized system for $5/5$ data sets at $T=20$ and $\varepsilon=0.01$,
it only stabilizes for $1/5$ data sets  when $\varepsilon=0.1$ for the same $T$.
Thus, we observe that non-linearity in trajectory data, both due to horizon length and data scale,
degrades the stabilization performance of linear models by distorting the 
forward conjugacy property, despite the encoder capturing the unstable subspace.

\paragraph{Non-linear representations induce more accurate linearizations.}
Contrary to linear latent representations, the performance of the non-linear, 
control-affine latent models does not collapse for large $T$. In \Cref{fig:exp_arm_nonlinear}, panels (a, b, c),
the linearization of the encoder latent model maintains alignment with the unstable subspace, 
and the linearized latent dynamics capture the dynamics of the unstable mode, with performance that does not 
degrade as in the linear case. Furthermore, panel (d) shows that the closed-loop dynamics induced on the linearized arm are stable, so
the MLP encoder and dynamics produce an accurate linear model of the unstable dynamics about the equilibrium.
For highly non-linear training data ($T=40$), the linearized MLP model achieves significantly better performance at stabilizing the true non-linear system,
and is able to stabilize from all configurations tested, and for all data scales. 
Furthermore, its performance is
more robust to the data scale $\varepsilon$:
for $T\ge 20$, the linearization of the MLP model remains a stabilizing model for the linearized system for all $\varepsilon$ 
(see \Cref{tab:exp_conjugacy_nonlinear,tab:exp_eps_tips}).

\paragraph{Linear control design for non-linear systems.}
In panel (e) of \Cref{fig:exp_arm_nonlinear}, latent output-feedback control via the encoder and latent dynamics 
linearizations exhibits strong performance relative to the other methods tested. The largest tip displacements, 
for which latent LQR using the linearized encoder $J_E(0)$ stabilizes all initial conditions tested, correspond to an 
angle of almost $80^\circ$ from the upright position. Comparatively, latent LQR using the non-linear $E(x)$
performs worse when data scale is smaller (i.e., at $\varepsilon=0.01$, performance improves for data scales $\varepsilon=0.03,0.1$, cf.~\Cref{tab:exp_eps_tips}). 
These results indicate that when the collected dataset is too concentrated about the equilibrium,
a non-linear encoder may extrapolate poorly outside of the data region during non-linear output feedback $u_t=KE(x_t)$.
Performing output-feedback via the linearization $u_t=KJ_E(0)x_t$ effectively restricts $E$'s behavior
to what was learned in the data region.

\subsection{Learning Several Unstable Modes}
\label{sec:experiments:multi}

\begin{figure}[ht!]
    \centering
    \includegraphics[width=\linewidth]{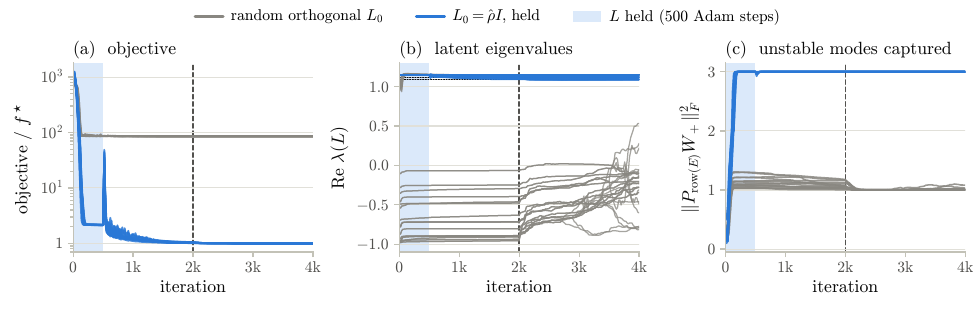}
    \caption{Training on empirical $\widehat{\calL}_{\mathrm{bwd}}^{(40)}$ 
    with three unstable modes, from the standard initialization (random orthogonal $L$, grey) and from $L = \hat{\rho}I_m$ held for the first $500$ Adam steps (blue), 
    for each of the 16 random initializations of one data set. The dashed line marks the switch from Adam to L-BFGS. 
    (a) Objective relative to the lowest one found, $f^\star$. (b) Real parts of the latent eigenvalues (dotted lines are $\sigma_+(A) = \{1.15, 1.12, 1.09\}$). 
    (c) $\norm{P_{\mathrm{row}(E)}W_+}_F^2$, where the columns of $W_+$ are the unit left eigenvectors of the unstable modes (equals $3$ when all three modes are captured).}
    \label{fig:exp_multi_unstable}
\end{figure}
In this experiment, 
we investigate whether gradient descent is able to find minimizers that capture all unstable modes of the system when there are multiple modes,
as the previous experiments,
\Cref{sec:experiments:linear}
and \Cref{sec:experiments:nonlinear},
only considered a single unstable buckling mode.
We consider the linearized arm model, and learn linear autoencoders and linear latent dynamics from trajectory data.
The optimization landscape for non-linear MLPs on non-linear data is significantly more complex, and we leave a more detailed
investigation to future work. 
We place $\kappa$ between the third and fourth largest generalized eigenvalues, so that three bending modes buckle, 
and set their modal stiffnesses so that $\sigma_+(A) = \{1.15, 1.12, 1.09\}$.
All other modes are unchanged.
We use mass-normalized coordinates $x = (M(0)^{1/2}q, M(0)^{1/2}\dot{q})$ ~\cite[Secs.~3.3,\ 3.6]{inman2006vibration},
in which the eigenvectors of different bending modes of \eqref{eq:arm_dynamics} are mutually-orthogonal,
and linear models with $m = m_+ = 3$;
in the original joint coordinates, the left eigenvectors of the two slowest unstable modes are approximately parallel, making optimization more challenging.
We choose $R_{\mathrm{orth}}$ for the regularizer, and sample $5$ data sets of $N = 2000$ trajectories,
with trajectory lengths $T$ in $\{10, 20, 30, 40\}$.
The initial conditions $x_0$ and i.i.d. control inputs $u_t$ are sampled
from $x_0\sim\mathsf{N}(0, I)$ and $u_t\sim\mathsf{N}(0, 0.5^2/T)$.
Since this experiment studies the initialization, each $T$ is trained separately, from $16$ initializations per data set, by $2000$ Adam steps followed by $2000$ L-BFGS steps.

\paragraph{Capturing unstable modes via growth-rate initialization.} 
We find that optimization of the backward loss
can converge to stationary points that do not capture the full unstable subspace, but where the gradient vanishes.
In our runs, a latent direction aligned with an unstable mode only when the corresponding eigenvalue of $L$ remained real and near or above that mode's eigenvalue early in training (\Cref{fig:exp_multi_unstable}b).
We therefore suggest an initialization, inspired by power iteration for identifying dominant eigenvalues,
that places the latent spectrum at the data's growth rate: keeping $(E, D, H)$ randomly initialized, we hold $L = \hat{\rho}I_m$ for the first $500$ Adam steps, 
where $\hat{\rho} := (\E\norm{x_T}^2/\E\norm{x_{T-1}}^2)^{1/2}$ is the empirical growth rate of the data.
This initialization captures all three modes %
in every run (cf. \Cref{fig:exp_multi_unstable}).

\begin{table}[ht!]
\centering
\small
\makeatletter\providecommand{\tabinput}[1]{\input{#1} }\makeatother
\begin{tabular}{lcccc}
\hline
initialization & $T = 10$ & $T = 20$ & $T = 30$ & $T = 40$ \\
\hline
random orthogonal $L$ & 0/80 & 0/80 & 0/80 & 0/80 \\
$L = 1.2\,I_m$ & 10/80 & 71/80 & 78/80 & 72/80 \\
$L = 1.2\,I_m$, held for $500$ steps & 80/80 & 80/80 & 80/80 & 80/80 \\
$L = \hat{\rho}\,I_m$, held for $500$ steps & 80/80 & 80/80 & 80/80 & 80/80 \\
 
\hline
\end{tabular}
\caption{Flexible arm with three unstable modes ($m = m_+ = 3$): We count runs (over $5$ data sets with $16$ random initializations each) that capture all three modes i.e., $\Delta(\calG^\R_+(A^\T); \mathrm{row}(E)) < 0.5$ (in every run, $\Delta \le 0.19$ or $\Delta \ge 0.82$). %
}
\label{tab:exp_multi_unstable}
\end{table}

\section{Conclusion and Future Work}
\label{sec:conclusion}

We studied when learned latent dynamics models enable feedback control design that stabilizes the original system.
For LTI systems, we identified encoder invariance and detectability as the structural conditions linking latent stabilization to stable zero dynamics.
We showed that minimizing a multi-step backward losses with sufficient long reconstruction horizon length approximately recovers this structure, and a minimum $H_\infty$-norm latent controller stabilizes the original system.
Our insufficiency results show that forward prediction losses and single-step backward reconstruction losses alone do not provide the same guarantee.
Experiments on a flexible-arm system illustrate these distinctions, and practical implications for non-linear models.

Our work suggests several exciting
directions for future work:
\begin{enumerate}[label=(\roman*)]
    \item \emph{Extensions of LTI analysis to partially observed settings:}
    When full-state measurements are unavailable, latent representations may need to incorporate histories of observations and inputs.
    A natural question is whether reconstruction losses involving these histories can recover the information needed for stabilization under suitable observability or detectability assumptions.
    
    \item \emph{An analysis for stabilization of non-linear systems:}
    Our experiments in \Cref{sec:experiments} suggest a natural starting point to study
    non-linear dynamics is to provide local asymptotic stability results near equilibrium points.
    Can we also go further and 
    transfer entire stability
    certificates from the latent dynamics back to the original system, going beyond convergence to the zero-manifold 
    discussed 
    in our prior work~\cite{lutkus2025latent}?
    An interesting technical question is identifying the correct non-linear extensions of encoder invariance and detectability that simultaneously
    arise from minimizers of non-linear $\calL_{\mathrm{bwd}}^{(T)}$, and also establish
    convergence to stable zero-dynamics.

    \item \emph{The optimization landscape of multi-step prediction-reconstruction losses:}
    Our theory characterizes properties of the global loss minimizers, while our experiments in \Cref{sec:experiments:multi} 
    illustrate the optimization challenges of gradient-based methods in finding global minimizers.
    Understanding how the reconstruction horizon and parameter initialization affect the optimization landscape of $\calL_{\mathrm{bwd}}^{(T)}$ could enable algorithm design for learning latent models with reliable convergence guarantees.
    
    \item \emph{Extensions beyond stabilization, including for safety barrier functions:}
    Latent representations for safety must preserve information about state constraints as well as dynamics.
    Safety detectability is shown to be an important property to ensure safety from latent design in the recent work of \citet{mootz2026latent}.
    How do we design prediction and reconstruction loss functions that provably ensure this property, and
    do multi-step losses also play a crucial role here?

\end{enumerate}

\section*{Acknowledgements}

We thank Leonardo Toso, James Anderson, and 
Florian D{\"{o}}rfler for their 
constructive and insightful feedback on our work.
This work was supported by the National Science Foundation through NSF CPS \#2434460 and an Okawa Foundation Research Grant.

\bibliographystyle{unsrtnat}
\bibliography{references}

\clearpage
\appendix
\begingroup
\renewcommand{\contentsname}{Contents}
\setcounter{tocdepth}{2}
\tableofcontents
\endgroup
\clearpage

\section{More Detailed Comparison to \citet{toso2026preserving}}
\label{sec:appendix:comparison_to_toso}

While the motivations of the concurrent work of \citet{toso2026preserving} are similar to ours, 
there are some key differences between our results.
Both our work and \cite{toso2026preserving} identify the insufficiencies of the forward prediction loss $\calL_{\mathrm{fwd}}^{(T)}$ for preserving unstable modes, 
in addition to identifying the role of detectability
of $(A, E)$ (cf.~\Cref{sec:problem_formulation:invariance_subspaces})
for latent stabilization.
However, the solution studied in \cite{toso2026preserving}
is quite different from ours; they propose using $H$-step action reconstruction (i.e., learning an inverse dynamics model)
in place of the SIGReg regularizer loss in $\calL_{\mathrm{fwd}}^{(T)}$, whereas our solution is based on
the multi-step backward conjugacy loss $\calL_{\mathrm{bwd}}^{(T)}$. 
This leads to
different analysis strategies; in particular our work studies 
optimizers
of losses that involve \emph{multi-step dynamics propagation},
which introduces a high degree of non-convexity.
Another key difference is in regards to the type of guarantees.
The main results of \cite{toso2026preserving}, in particular \cite[Thms.\ 1 and 2]{toso2026preserving}, show that
when the action-reconstruction loss is zero, the learned encoder must
satisfy the requisite detectability condition. 
However, this only shows the \emph{existence} of a latent controller $K$ that stabilizes $(A, B)$; it does not imply that such a controller can actually be recovered from the learned latent model
(cf.~\Cref{prop:latent_stability_iff_output_feedback_stability}, where the $A$-invariance of $\ker(E)$ is 
necessary to ensure latent stability transfer).
Furthermore, our results in \Cref{sec:results} (cf.~\Cref{cor:latent_stability_transfer})
prove that minimum $H_\infty$-norm optimization operating solely on the \emph{learned} latent model
$(\bar{L}_T, \bar{H}_T)$ produces a stabilizing controller for the original system, providing an end-to-end result for the 
latent stabilization problem (cf.~\Cref{prob:latent_stabilization}).

At a lower-level, there are several technical differences between our results. First, \cite[Thm.~1]{toso2026preserving} requires that 
the horizon $H$ is sufficient large, in addition to the latent state dimension $m$ satisfying $m \geq pH$ (using our notation) in order for the action-reconstruction loss to be identically zero, which is a required condition for their theorem to hold; it is noted (cf.~\cite[bottom of page 8]{toso2026preserving}) that in practice such action-reconstruction
losses are not actually enforced to be zero.
On the other hand, our main result \Cref{thm:informal_recovery}
applies when $m \geq m_+$, which is the minimum latent dimension size
needed to represent all the non-stable modes of $A$,
and does not require that $\calL_{\mathrm{bwd}}^{(T)}$ is zero at the 
optimal solution.
Second, \cite[Thm.~2]{toso2026preserving} requires that the matrix $A$ has no
marginally stable eigenvalues,
and also imposes extra coupling conditions between the unstable and stable modes of $A$ (cf.~\cite[Assumptions 1 and 2]{toso2026preserving}), which our theory does not require
due to our use of {invariant subspaces} 
and canonical variables (cf.~\Cref{sec:canonical_variables}).

With this in mind, we note that \cite{toso2026preserving} works with the more general partially observed LTI setting where $y_t = C x_t$ (instead of the full-state observation setting considered in our work), and
also conducts experiments for non-linear systems using visuomotor control; these are 
important future direction for our work, as discussed in \Cref{sec:conclusion}.
Furthermore, when the observation $y_t$ is high dimension (e.g., 
an RGB-D image of the scene), the action-reconstruction
regularization proposed in 
\cite{toso2026preserving}
is more computationally appealing
than doing pixel-level reconstruction, which is the natural
extension of $\calL_{\mathrm{bwd}}^{(T)}$. Hence, an interesting future direction is showing  end-to-end stabilization results for the practical action-reconstruction losses of \cite{toso2026preserving}, ideally
on both linear and non-linear systems.

\section{Analysis of Exact Invariance Constructions}
\label{sec:appendix:proofs_for_exact_invariance}

In this section we supply the missing proofs from \Cref{sec:problem_formulation:invariance_subspaces}.
We first restate and prove
\Cref{prop:invariance_equiv_latent_dynamics}.
\invarianceequivlatentdynamics*
\begin{proof}
$(\Rightarrow)$
Suppose such an $f_z$ exists. Considering one-step trajectories gives
$f_z(Ex,u)=E(Ax+Bu)$ for all $x,u$, so $f_z(0,0)=0$.
Hence for $v\in\ker(E)$, we have
$0 = f_z(Ev,0) = EAv$, implying $v \in \ker(EA)$ as well.
Thus $\ker(E)\subseteq\ker(EA)$, 
which means that $Av\in\ker(E)$ for every $v\in\ker(E)$; hence $\ker(E)$ is $A$-invariant.

$(\Leftarrow)$
$\ker(E)$ is $A$-invariant 
is equivalent to 
$A^\T$-invariance of $\mathrm{row}(E)$, which
implies $EA=LE$ for some $L\in\R^{m\times m}$.
Therefore $f_z(z,u)=Lz+EBu$ satisfies the required identity.
Any admissible $f_z$ must agree with this formula for every
$z=Ex\in\mathrm{col}(E)$ and $u\in\R^p$, proving linearity as well.
\end{proof}

Next we restate and prove
\Cref{prop:latent_stability_iff_output_feedback_stability}
\latentstabilityiffoutputfeedbackstability*
\begin{proof}
Let $T=[T_1\ T_2]$ be any invertible matrix whose columns in $T_1$ and $T_2$
form bases for $\mathrm{row}(E)$ and $\ker(E)$, respectively. Define:
\begin{align*}
    T^{-1}AT = \begin{bmatrix} A_{11} & 0 \\ A_{21} & A_{22} \end{bmatrix},
    \qquad T^{-1}B = \begin{bmatrix} B_1 \\ B_2 \end{bmatrix}.
\end{align*}
The zero upper-right block follows from $EA=LE$ and $ET_2=0$.
Since $E$ has full row rank, $R:=ET_1$ is invertible, and
$LR=RA_{11}$ and $H=RB_1$.
Consequently,
\begin{align*}
    T^{-1} (A+BKE) T
    = \begin{bmatrix}
        R^{-1}(L+HK)R & 0 \\
        A_{21} + B_2 K R & A_{22}
    \end{bmatrix}.
\end{align*}
The subspace $\ker(E)$ is $A$-invariant and unobservable, so detectability
implies that $A_{22}$ is Schur stable;
indeed if $A_{22} v = \lambda v$ for some
non-zero $v$ and $\abs{\lambda} \geq 1$, then
$A (T_2 v) = T_2 A_{22} v = \lambda (T_2 v)$ and
$E (T_2 v) = 0$; since $T_2 v \neq 0$, this would contradict detectability of $(A, E)$ by the PBH test.
Hence, $A+BKE$ is therefore Schur stable if and only if
$L+HK$ is Schur stable.

Hence, if $(L,H)$ is
stabilizable, we have just shown that $(A,B)$ is stabilizable. 
Conversely, suppose
$(A,B)$ is stabilizable, 
but suppose $(L, H)$ is not.
By PBH there exists non-zero $q$ and
$\abs{\lambda} \geq 1$ such that
$q^*L=\lambda q^*$ and $q^*H=0$. Since $E$ has full row rank,
$q^*E\neq 0$, while
\begin{align*}
    (q^*E)A=q^*LE=\lambda(q^*E), \qquad (q^*E)B=q^*H=0.
\end{align*}
This contradicts the PBH criterion for stabilizability of $(A,B)$,
completing the proof.
\end{proof}

Finally we restate and prove
\Cref{prop:compatibility_assumptions_equiv}
\compatibilityassumptionsequiv*
\begin{proof}
We first show that there exists an encoder $E \in \R^{m \times n}$ satisfying conditions (i)--(iii) of \Cref{prop:latent_stability_iff_output_feedback_stability} if and only if there exists an
$A$-invariant subspace $N$ of dimension $n-m$ such that the restriction $A|_N$ is Schur stable.
$(\Rightarrow)$:
Since $\ker(E)$ is $A$-invariant 
(condition (ii)), then
$\ker(E)$ is exactly the unobservable subspace
$U := \{ x \in \R^n \mid EA^t x = 0,\,\, \forall t \geq 0 \}$, since 
$Ex=0$ implies $EA^t x=0$ for all $t\geq0$
and conversely $x \in U$ implies $Ex = 0$.
Detectability (condition (i)) is therefore equivalent to Schur stability of $A|_{\ker(E)}$, and hence we can take
$N = \ker(E)$, which has
dimension $n - m$ by condition (iii).
$(\Leftarrow)$: Given such an $N$, choose the rows of $E$ to be a basis of $N^\perp$, which immediately
satisfies conditions (i)--(iii).

Hence, it suffices to show that
there exists an $A$-invariant subspace $N$ of dimension $n-m$ such that $A|_N$ is Schur stable
if and only if $m \geq m_+$ and either condition (i) or (ii) from the hypothesis holds.

$(\Rightarrow)$: The eigenvalues of $A|_N$ must be among the stable eigenvalues of $A$.
There are $n-m_+$ such eigenvalues (counted with algebraic
multiplicity), 
so $n-m\leq n-m_+$, or $m\geq m_+$.
If none of the stable eigenvalues are real, they occur in conjugate
pairs, both for $A$ and $A|_N$.
Thus $n-m_+$ and $n-m$ are both even, and hence $m-m_+$ must be even as well.

$(\Leftarrow)$: Suppose the stated conditions.
If $m=n$, take $N=\{0\}$.
Otherwise, we will utilize the following procedure
to construct a multi-set $\Lambda_r \subseteq \eig_-(A)$ of stable eigenvalues of size $n-m$. First, we initialize:
\begin{align*}
    \Lambda_r = \{\}, \qquad \Sigma_s = \mathrm{eig}_-(A).
\end{align*}
Then we follow the procedure:
\begin{enumerate}[label=(\arabic*)]
    \item If $n - m$ is odd, then move one real-valued $\lambda \in \Sigma_s$
    from $\Sigma_s$ into $\Lambda_r$.
    \item While $\abs{\Lambda_r} < n-m$:
    move a non-real conjugate pair from $\Sigma_s$ to $\Lambda_r$ if available, otherwise move any two
    real eigenvalues from $\Sigma_s$ to $\Lambda_r$. 
\end{enumerate}
Note that our conditions imply that if
there is no $\lambda \in \sigma_-(A)$ that is real-valued, then
we must have both $n-m_+$ and $m-m_+$ even, and hence $n-m$ must be even as well; this implies step (1) is always possible.
Furthermore, the loop invariant in step (2) is that
the difference $n-m-\abs{\Lambda_r}$ is always even
and both $\Lambda_r$ and $\Sigma_s$ are maintained to be conjugate invariant, and hence 
step (2) is always possible as $n-m \leq n-m_+$
(since $m \geq m_+$).
Now we apply the Schur decomposition in \eqref{eq:real_schur_decomposition}
to compute $(Q, S) = \mathrm{RealSchur}(A, \Lambda_r)$, obtaining:
\begin{align*}
    A=Q\bmattwo{S_{11}}{S_{12}}{0}{S_{22}}Q^\T,
    \qquad Q=\rvectwo{Q_1}{Q_2} \in O(n), \qquad
    Q_1\in\R^{n\times (n-m)}, \qquad \eig(S_{11})=\Lambda_r.
\end{align*}
Since $AQ_1=Q_1S_{11}$ and $S_{11}$ is Schur stable,
$N:=\mathrm{col}(Q_1)$ has the required properties.
\end{proof}

\section{Regularizer Admissibility for Commonly Used Regularizers}
\label{sec:appendix:proofs_for_regularizer_admissibility}

We restate and prove \Cref{prop:regularizer_admissibility}
regarding regularizer admissibility.
\regularizeradmissibility*
\begin{proof}
Fix $E\in\R^{m\times n}$.
We verify \Cref{def:forward_admissible_regularizer} for each regularizer,
using $\sigma_0=1/2$ and $\sigma_1=3/2$ throughout.

\paragraph{Orthogonality regularizer $R_{\mathrm{orth}}$.}
Since the encoder is linear, we have $J_E(x) = E$.
By \eqref{eq:regularizer_orth},
\begin{align*}
    R_{\mathrm{orth}}(E;\mu_n)=\norm{I_m-EE^\T}_F^2.
\end{align*}
This clearly vanishes when $E^\T\in O(n,m)$.
Taking $\e_0=1/4$ also gives:
\begin{align*}
    R_{\mathrm{orth}}(E;\mu_n)\leq\e_0
    \implies\opnorm{EE^\T-I_m}\leq\frac12
    \implies\frac12 I_m\preccurlyeq EE^\T \preccurlyeq\frac32 I_m,
\end{align*}
proving admissibility.

\paragraph{Variance--covariance regularizer $R_{\mathrm{vr}}$.}
Since the encoder is linear we have $Ex \sim \sfN(0, EE^\T)$ for $x \sim \mu_n$, and hence $C_{\mu_n}(E) = EE^\T$.
For every $\beta>0$, \eqref{eq:regularizer_VR} gives
\begin{align*}
    R_{\mathrm{vr}}(E;\mu_n,\beta)
    =\frac{1}{m}\sum_{i=1}^m(1-(EE^\T)_{ii})^2
        +\frac{\beta}{m}\sum_{i\neq j}(EE^\T)_{ij}^2 
    \geq\frac{\min\{1,\beta\}}{m}\norm{I_m - EE^\T}_F^2.
\end{align*}
Again, this regularizer vanishes when $E^\T \in O(n, m)$.
Taking $\e_0=\min\{1,\beta\}/(4m)>0$ ensures that:
\begin{align*}
    R_{\mathrm{vr}}(E;\mu_n,\beta)\leq\e_0 \implies \opnorm{EE^\T-I_m}\leq\frac12,
\end{align*}
proving admissibility.

\paragraph{SIGReg regularizer $R_{\mathrm{sr}}$.}
For each unit vector $a \in \bbS^{m-1}$, the Gaussian characteristic function gives
\begin{align*}
    \E_{x\sim\mu_n}[e^{jta^\T Ex}]=e^{-t^2 \norm{E^\T a}^2/2}.
\end{align*}
Evaluating the integral in \eqref{eq:regularizer_SR} therefore yields
\begin{align*}
    R_{\mathrm{sr}}(E;\mu_n) =
        \E_{a\sim\mathrm{Unif}(\mathbb S^{m-1})}[g(\norm{E^\T a}^2)], \qquad g(s) := \int_{-\infty}^{\infty} (e^{-st^2/2}-e^{-t^2/2})^2e^{-t^2/2}\,\rmd t ,\,\, s\geq0.
\end{align*}
Since $g(1)=0$, we have $g(\norm{E^\T a}^2) = g(1) = 0$ for all $a$ when
$E^\T \in O(n, m)$, which shows 
$R_{\mathrm{sr}}(E;\mu_n)$ vanishes for $E^\T \in O(n, m)$.

Next, we establish the bound:
\begin{align}
    g(s)\geq c\min\{(s-1)^2,1\},\qquad s\geq0,
    \qquad c:=\frac{3\sqrt{2\pi}}{4\cdot5^{5/2}}. \label{eq:regularizer_admissibility:gs_bound}
\end{align}
First, the mean value theorem 
on the function $s \mapsto e^{-st^2/2}$ yields:
\begin{align*}
    \abs{e^{-st^2/2}-e^{-t^2/2}}
        \geq\frac{t^2}{2}e^{-t^2}\abs{s-1}, \qquad s \in [0, 2],
\end{align*}
which implies that:
\begin{align*}
    g(s)\geq\frac{(s-1)^2}{4}
        \int_{-\infty}^{\infty} t^4e^{-5t^2/2}\,\rmd t
        =c(s-1)^2, \qquad s \in [0, 2].
\end{align*}
For $s\geq2$, the integrand defining $g(s)$ is non-decreasing in $s$,
so $g(s)\geq g(2)\geq c$, proving
\eqref{eq:regularizer_admissibility:gs_bound}.

Using \eqref{eq:regularizer_admissibility:gs_bound}, we have, putting $\Delta_E := EE^\T - I_m$ and supposing
$\Delta_E \neq 0$,
\begin{align*}
    R_{\mathrm{sr}}(E; \mu_n) &\geq c \E_a[ \min\{(\norm{E^\T a}^2 - 1)^2, 1\} ] = c \E_a[ \min\{ (a^\T \Delta_E a)^2, 1 \} ] \\
    &\geq c \E_a[ \min\{ (a^\T \Delta_E a)^2, (a^\T \Delta_E a)^2 / \opnorm{\Delta_E}^2 \} ] \\
    &= c \min\{1, 1/\opnorm{\Delta_E}^2 \} \E_a[ (a^\T \Delta_E a)^2 ] \\
    &\stackrel{(a)}{=} c \min\{1, 1/\opnorm{\Delta_E}^2 \} \frac{2\norm{\Delta_E}_F^2+(\Tr(\Delta_E))^2}{m(m+2)} \\
    &\geq \frac{2c}{m(m+2)} \min\{ \opnorm{\Delta_E}^2, 1 \}.
\end{align*}
Here the equality in (a) follows by diagonalizing $\Delta_E$ and using the spherical moments
$\E[a_i^4]=3/[m(m+2)]$ and
$\E[a_i^2a_j^2]=1/[m(m+2)]$ for $i\neq j$.
Now setting $\e_0=\frac{c}{2m(m+2)}$,
we have that:
\begin{align*}
    R_{\mathrm{sr}}(E;\mu_n)\leq\e_0
    \implies\opnorm{\Delta_E}\leq\frac12.
\end{align*}
Note that the excluded case $\Delta_E = 0$ also
satisfies this bound.
This proves admissibility of $R_{\mathrm{sr}}$.
\end{proof}

\section{Canonical Variables from Real Schur Decompositions}
\label{sec:canonical_variables}

In this section, we develop a set of oblique projectors
for $A$ based on the real Schur decomposition machinery
introduced in \Cref{sec:problem_formulation:invariance_subspaces}.
These projectors play a crucial role in the subsequent analysis, allowing us to capture both
the non-stable subspaces of $A$, in addition to 
parts of the stable subspaces, in a way that is compatible
with the conditions in \Cref{prop:latent_stability_iff_output_feedback_stability}.
To start, we will introduce the notion of a set of \emph{canonical variables} indexed by a multi-set of stable eigenvalues $\Lambda_s$.

\begin{mydef}[Canonical variables indexed by $\Lambda_s$]
\label{def:canonical_variables_v2}
Let $\Lambda_s \subseteq \eig_-(A)$
be a (possibly empty) multi-set of stable eigenvalues
satisfying conjugate invariance
(i.e., $\mathrm{mult}_{\Lambda_s}(\lambda) = \mathrm{mult}_{\Lambda_s}(\bar{\lambda})$ 
for all $\lambda \in \Lambda_s$).
Put $\Lambda_{r} := \eig(A) \setminus (\mathrm{eig}_+(A) \cup \Lambda_s)$.
Let $(Q^{\Lambda_s}, (R^{\Lambda_s})^\T) = \mathrm{RealSchur}(A^\T, \eig_+(A) \cup \Lambda_s)$, and write both 
$(Q^{\Lambda_s}, R^{\Lambda_s})$ in block matrix form satisfying
\begin{align}
    A = Q^{\Lambda_s} R^{\Lambda_s} (Q^{\Lambda_s})^\T = \begin{bmatrix} 
        Q^{\Lambda_s}_{u} & Q^{\Lambda_s}_{s} & Q^{\Lambda_s}_{r}
    \end{bmatrix} \begin{bmatrix}
        R^{\Lambda_s}_{uu} & 0 & 0 \\
        R^{\Lambda_s}_{su} & R^{\Lambda_s}_{ss} & 0 \\
        R^{\Lambda_s}_{ru} & R^{\Lambda_s}_{rs} & R^{\Lambda_s}_{rr}
    \end{bmatrix} 
    \begin{bmatrix}
        (Q^{\Lambda_s}_{u})^\T \\
        (Q^{\Lambda_s}_{s})^\T \\
        (Q^{\Lambda_s}_{r})^\T
    \end{bmatrix} , 
    \label{eq:schur_decomposition}
\end{align}
where the following properties hold:
\begin{enumerate}[label=(\alph*)]
    \item $Q^{\Lambda_s}_{u} \in O(n, m_+)$, $Q^{\Lambda_s}_{s} \in O(n, d_s)$,
    and $Q^{\Lambda_s}_{r} \in O(n, d_r)$,
    with 
    $d_s := \abs{\Lambda_s}$ and
    $d_r := n - m_+ - d_s$,
    \item $\eig(R^{\Lambda_s}_{uu}) = \eig_+(A)$,
    $\eig(R^{\Lambda_s}_{ss}) = \Lambda_s$, and
    $\eig(R^{\Lambda_s}_{rr}) = \Lambda_r$.
\end{enumerate}
Let $X^{\Lambda_s} \in \R^{(d_s + d_r) \times m_+}$ denote the unique solution to the following
Sylvester equation:
\begin{align}
    X^{\Lambda_s} R^{\Lambda_s}_{uu} - \bmattwo{R^{\Lambda_s}_{ss}}{0}{R^{\Lambda_s}_{rs}}{R^{\Lambda_s}_{rr}} X^{\Lambda_s} = \cvectwo{R^{\Lambda_s}_{su}}{R^{\Lambda_s}_{ru}}. \label{eq:sylvester_equation}
\end{align}
Such a solution exists and is unique since
$\sigma(R^{\Lambda_s}_{uu})$ and $ \sigma(R^{\Lambda_s}_{ss}) \cup \sigma(R^{\Lambda_s}_{rr})$ are disjoint by construction.
Partition $X^{\Lambda_s} = \cvectwo{ X^{\Lambda_s}_{s} }{ X^{\Lambda_s}_{r} }$,
where $X_s^{\Lambda_s} \in \R^{d_s \times m_+}$ and $X_r^{\Lambda_s} \in \R^{d_r \times m_+}$, and define the \emph{canonical variables}:
\begin{align}
    E^{\Lambda_s} &:=  
    \begin{bmatrix}
    (Q^{\Lambda_s}_{u})^\T \\
        (Q^{\Lambda_s}_{s})^\T 
    \end{bmatrix}, \quad D^{\Lambda_s} := \begin{bmatrix}
        Q_u^{\Lambda_s} + Q_r^{\Lambda_s} X_r^{\Lambda_s} & Q_s^{\Lambda_s}
    \end{bmatrix}, \quad L^{\Lambda_s} := \begin{bmatrix}
        R^{\Lambda_s}_{uu} & 0 \\
        R^{\Lambda_s}_{su} & R^{\Lambda_s}_{ss}
    \end{bmatrix}, \quad
    H^{\Lambda_s} := E^{\Lambda_s} B, \label{eq:canonical_variables} \\
        \Pi_+^{\Lambda_s} &:= D^{\Lambda_s} E^{\Lambda_s} , \qquad 
        \Pi_-^{\Lambda_s} := I_n - \Pi_+^{\Lambda_s} .
\end{align}
\end{mydef}

The following result lists the basic properties 
satisfied by the canonical variables
defined in \Cref{def:canonical_variables_v2}.
\begin{myprop}[Basic properties of canonical variables]
\label{prop:canonical_variables}
Suppose that \Cref{assmp:non_trivial_unstable_subspace} holds.
Let $\Lambda_s$ satisfy the conditions of
\Cref{def:canonical_variables_v2}, and
let $(E^{\Lambda_s}, D^{\Lambda_s}, L^{\Lambda_s}, H^{\Lambda_s}, \Pi_+^{\Lambda_s}, \Pi_-^{\Lambda_s})$ denote the corresponding canonical variables.
The following properties hold:
\begin{enumerate}[label=(\alph*)]
    \item (Right inverse): $E^{\Lambda_s} D^{\Lambda_s} = I_d$ for $d := m_+ + d_s$.
    \item (Oblique projectors): $\Pi_+^{\Lambda_s}, \Pi_-^{\Lambda_s}$ are oblique projectors, i.e., $(\Pi_{\pm}^{\Lambda_s})^2 = \Pi_{\pm}^{\Lambda_s}$.
    \item (Projector norm inequality):
    $\opnorm{ \Pi_-^{\Lambda_s} } \leq \opnorm{ \Pi_+^{\Lambda_s} }$.
    \item (Invariant subspace): 
    $E^{\Lambda_s} A = L^{\Lambda_s} E^{\Lambda_s}$, 
    and if $d_s = 0$ (i.e., $\Lambda_s = \{\}$) 
    then $A D^{\Lambda_s} = D^{\Lambda_s} L^{\Lambda_s}$.
    \item (Decay of stable decomposition): For every $k \geq 0$
    and $\gamma \in (\gamma_-, 1)$,
    \begin{align*}
        \Pi_+^{\Lambda_s} A^k &= D^{\Lambda_s} (L^{\Lambda_s})^k E^{\Lambda_s}, \qquad \opnorm{\Pi_-^{\Lambda_s} A^k}
        \leq C(A, \gamma) \gamma^k,
    \end{align*}
    where $C(A, \gamma)$ is a finite
    constant depending only on $A, \gamma$.
\end{enumerate}
\end{myprop}
\begin{proof}
We suppress the superscript $\Lambda_s$
on the canonical variables throughout the proof.

\paragraph{Part (a).}
Orthogonality of $Q=[Q_u\ Q_s\ Q_r] \in O(n)$ gives:
\begin{align*}
    ED
    = \cvectwo{Q_u^\T}{Q_s^\T}
      \begin{bmatrix}Q_u+Q_rX_r & Q_s\end{bmatrix}
    = I_d,
\end{align*}
which proves (a). 

\paragraph{Part (b).}
This is an immediate consequence of (a), since
$\Pi_+^2 = (DE)^2=D(ED)E=DE = \Pi_+$ and $\Pi_-^2 = (I_n-\Pi_+)^2=I_n-2\Pi_+ + \Pi_+^2 = I_n - \Pi_+ = \Pi_-$, proving (b).

\paragraph{Part (c).}
By \Cref{assmp:non_trivial_unstable_subspace}, $d=m_++d_s\geq1$.
If $d_r=0$, then $\Pi_+=I_n$ and $\Pi_-=0$, so the inequality is immediate.
Otherwise, let $F:=\rvectwo{X_r}{0}\in\R^{d_r\times d}$.
Grouping the $u$ and $s$ coordinates together, the definitions of $D,E$ give
\begin{align*}
    Q^\T\Pi_+Q=\bmattwo{I_d}{0}{F}{0},
    \qquad Q^\T\Pi_-Q=\bmattwo{0}{0}{-F}{I_{d_r}}.
\end{align*}
Since $Q$ is orthogonal, these representations imply
\begin{align*}
    \opnorm{\Pi_+}^2
    =\lambda_{\max}(I_d+F^\T F)=1+\opnorm{F}^2,\qquad
    \opnorm{\Pi_-}^2
    =\lambda_{\max}(I_{d_r}+FF^\T)=1+\opnorm{F}^2.
\end{align*}
Thus, when $d_r>0$, we have the stronger equality
$\opnorm{\Pi_-}=\opnorm{\Pi_+}=\sqrt{1+\opnorm{X_r}^2}$, proving (c).

\paragraph{Part (d).} 
The first two block rows of $Q^\T A$ give:
\begin{align*}
    EA
    = \begin{bmatrix} R_{uu} & 0 & 0 \\
        R_{su} & R_{ss} & 0 \end{bmatrix}
      \begin{bmatrix}
        Q_u^\T \\
        Q_s^\T \\
        Q_r^\T
      \end{bmatrix}
    = \bmattwo{R_{uu}}{0}{R_{su}}{R_{ss}}
      \cvectwo{Q_u^\T}{Q_s^\T}
    = LE,
\end{align*}
which proves the first identity.
If $d_s=0$, then $D=Q_u+Q_rX_r$ and $L=R_{uu}$, while the Sylvester equation \eqref{eq:sylvester_equation} reduces to $X_rR_{uu}-R_{rr}X_r=R_{ru}$.
Using $AQ_u=Q_uR_{uu}+Q_rR_{ru}$ and $AQ_r=Q_rR_{rr}$ from the Schur decomposition \eqref{eq:schur_decomposition}, we obtain
\begin{align*}
    AD
    &=A(Q_u+Q_rX_r)\\
    &=Q_uR_{uu}+Q_r(R_{ru}+R_{rr}X_r)\\
    &=Q_uR_{uu}+Q_rX_rR_{uu}
      =(Q_u+Q_rX_r)R_{uu}=DL.
\end{align*}
This proves the second identity and completes (d).

\paragraph{Part (e).}
Iterating $EA=LE$ yields $EA^k=L^kE$ for every $k\geq0$.
Thus
\begin{align*}
    \Pi_+A^k=DEA^k=DL^kE.
\end{align*}
To bound the decay rate, first expand $I_n-DE$ to write
\begin{align*}
    \Pi_-A^k=Q_r\bigl(Q_r^\T-X_rQ_u^\T\bigr)A^k.
\end{align*}
Now, define:
\newcommand{\stab}{\mathrm{stab}}
\begin{align*}
    Q_{\stab} := \rvectwo{Q_s}{Q_r}, \qquad R_{\stab} := \bmattwo{R_{ss}}{0}{R_{rs}}{R_{rr}}, \qquad R_{\stab,u} := \cvectwo{R_{su}}{R_{ru}},
\end{align*}
so that the Schur decomposition for $A$ and Sylvester equation for $X$ (cf.~\eqref{eq:sylvester_equation}) can be written as:
\begin{align*}
    A = \rvectwo{Q_u}{Q_{\stab}} \bmattwo{R_{uu}}{0}{R_{\stab,u}}{R_{\stab}} \cvectwo{Q_u^\T}{Q_{\stab}^\T}, \qquad  X R_{uu} - R_{\stab} X = R_{\stab,u}.
\end{align*}
Now defining $W := Q_{\stab}^\T - X Q_u^\T$, we show that $\mathrm{row}(W)$ is $A^\T$-invariant:
\begin{align*}
    WA &= (Q_{\stab}^\T - X Q_u^\T)\rvectwo{Q_u}{Q_{\stab}} \bmattwo{R_{uu}}{0}{R_{\stab,u}}{R_{\stab}} \cvectwo{Q_u^\T}{Q_{\stab}^\T} \\
    &= \rvectwo{-X}{I_{n-m_+}} \bmattwo{R_{uu}}{0}{R_{\stab,u}}{R_{\stab}} \cvectwo{Q_u^\T}{Q_{\stab}^\T} \\
    &= (-X R_{uu} + R_{\stab,u}) Q_u^\T + R_{\stab} Q_{\stab}^\T \\
    &= R_{\stab} (-X Q_u^\T + Q_{\stab}^\T ) && \text{[using Sylvester equation]}\\
    &= R_{\stab} W.
\end{align*}
Iterating this identity 
yields $W A^k = R_{\stab}^k W$ for all $k \geq 0$.
Observing that $\rvectwo{0}{I_{d_r}} W = Q_r^\T - X_r Q_u^\T$, we have:
\begin{align*}
    \Pi_- A^k = Q_r\bigl(Q_r^\T-X_rQ_u^\T\bigr)A^k = Q_r \rvectwo{0}{I_{d_r}} W A^k = Q_r \rvectwo{0}{I_{d_r}} R_{\stab}^k W.
\end{align*}
Since $\sigma(R_{\stab}) = \sigma_-(A)$, we have
$\rho(R_{\stab})=\gamma_-<1$.
Hence for any $\gamma\in(\gamma_-,1)$:
\begin{align*}
    M(R_{\stab},\gamma):=\sup_{k\geq0}\gamma^{-k}\opnorm{R_{\stab}^k}<\infty.
\end{align*}
Since $WW^\T=I_{n-m_+}+XX^\T$, we obtain
\begin{align*}
    \opnorm{\Pi_-A^k}
    \leq\opnorm{R_{\stab}^k}\opnorm{W}
    \leq M(R_{\stab},\gamma)\sqrt{1+\opnorm{X}^2}\,\gamma^k.
\end{align*}

To see that the prefactor
$M(R_{\stab},\gamma)\sqrt{1+\opnorm{X}^2}$
only depends on $A$ and $\gamma$,
let $\widetilde{\Lambda}_s$ denote any other admissible multi-set (not necessarily equal to $\Lambda_s$), and let
$A = \widetilde{Q} \widetilde{R} \widetilde{Q}^\T$
denote its corresponding Schur decomposition 
\eqref{eq:schur_decomposition}
and let $\widetilde{X}$ denote its corresponding
solution to the Sylvester equation \eqref{eq:sylvester_equation}.
First, we note that
$A \bar{Q}_{\stab} = \bar{Q}_{\stab} \bar{R}_{\stab}$ and
$\sigma(\bar{R}_{\stab}) = \sigma_-(A)$
for both $(\bar{Q}_{\stab}, \bar{R}_{\stab}) \in \{ (Q_{\stab}, R_{\stab}), (\widetilde{Q}_{\stab}, \widetilde{R}_{\stab}) \}$.
Thus $\mathrm{col}(\bar{Q}_{\stab})$ is an $A$-invariant subspace with only stable eigenvalues, whose dimension $n-m_+$ equals the total algebraic multiplicity of the stable eigenvalues of $A$.
Therefore:
\begin{align*}
    \mathrm{col}(Q_{\stab}) = \mathrm{col}(\widetilde{Q}_{\stab}) = \left( \bigoplus_{\lambda \in \sigma_-(A)} \calG_\lambda(A)  \right) \cap \R^n.
\end{align*}
Furthermore since both $Q, \widetilde{Q} \in O(n)$, $Q_u$ and $\widetilde{Q}_u$ span the orthogonal complement of the same subspace, and hence $\mathrm{col}(Q_u) = \mathrm{col}(\widetilde{Q}_u)$.
Consequently there exist $V \in O(n - m_+)$ and $U \in O(m_+)$ such that:
\begin{align*}
    \widetilde Q_{\stab}=Q_{\stab}V,\qquad
    \widetilde Q_u=Q_uU.
\end{align*}
Using these changes of basis, the blocks in the Sylvester equation \eqref{eq:sylvester_equation} satisfy:
\begin{align*}
    \widetilde R_{\stab}
      =V^\T R_{\stab}V, \qquad 
    \widetilde R_{uu}
      =U^\T R_{uu}U, \qquad
    \widetilde R_{\stab,u}
      =V^\T R_{\stab,u} U.
\end{align*}
It follows that $V^\T XU$ solves the Sylvester equation defining $\widetilde X$:
\begin{align*}
    (V^\T XU)\widetilde R_{uu}
      -\widetilde R_{\stab}(V^\T XU)
    =V^\T(XR_{uu}-R_{\stab}X)U
    =V^\T R_{\stab,u}U
      =\widetilde R_{\stab,u}.
\end{align*}
Since $\widetilde R_{uu}$ and $\widetilde R_{\stab}$ have disjoint spectra, this solution is unique, and hence $\widetilde X=V^\T XU$.
Orthogonal invariance of the operator norm therefore gives
\begin{align*}
    \opnorm{\widetilde X}=\opnorm{X},\qquad
    \opnorm{\widetilde R_{\stab}^{\,k}}
    =\opnorm{V^\T R_{\stab}^kV}
    =\opnorm{R_{\stab}^k}\quad\text{for every }k\geq0.
\end{align*}
Thus we see that
$M(\widetilde R_{\stab},\gamma)=M(R_{\stab},\gamma)$, and therefore
$M(R_{\stab},\gamma)\sqrt{1+\opnorm{X}^2}$
is unchanged by either the choice of admissible multi-set $\Lambda_s$ or the corresponding Schur bases, and depends only on $A$ and $\gamma$, finishing the proof of (e).
\end{proof}

We finish off this section with the following two instantiations of canonical variables (cf.~\Cref{def:canonical_variables_v2}) which we will repeatedly
use in the sequel.

\begin{mydef}[Canonical variables for non-stable modes]
\label{def:canonical_variables_unstable}
Let $E_0, D_0, L_0, H_0, \Pi_{\pm,0}$
denote canonical variables from \Cref{def:canonical_variables_v2}
for the choice of $\Lambda_s = \{ \}$.
Let $K_0 \in \R^{p \times m_+}$ be any
latent controller so that
$A + B K_0 E_0$ is Schur stable, and set
$\Gamma_0 := \sup_{z \in \mathbb{D}^c} \opnorm{ (zI_n - A - B K_0 E_0)^{-1}}$ (cf.~\Cref{sec:problem_formulation:invariance_subspaces}
for why such a $K_0$ exists).
\end{mydef}

\begin{mydef}[Canonical variables for latent modes]
\label{def:canonical_variables_latent}
Under \Cref{assmp:latent_dimension_compatibility}, 
let $E_m, D_m, L_m, H_m, \Pi_{\pm,m}$
denote canonical variables from \Cref{def:canonical_variables_v2}
for $\Lambda_s = \eig_-(A) \setminus \Lambda_r$,
where $\Lambda_r$ is computed in \Cref{alg:exact_solution_given_A_B}.
\end{mydef}

\section{Closed-Loop Stability Transfer Under Approximate Conditions}

In this section, we develop the perturbation analysis needed
to transfer closed-loop stability under only approximate invariance
and detectability.
Recall from \Cref{prop:latent_stability_iff_output_feedback_stability}
that if $E$ satisfies
(i) the pair $(A, E)$ is detectable,
(ii) $EA = LE$ for some $L \in \R^{m \times m}$, and
(iii) $\rank(E) = m$,  
then $K$ stabilizes $(L, H)$ (for $H = EB$) iff $KE$ stabilize $(A, B)$. 
As stated, this result requires conditions (i) and (ii)
to hold exactly, not allowing for e.g., approximate 
invariance $EA \approx LE$, or some suitable notion of almost detectable. The following result relaxes 
the requirements of \Cref{prop:latent_stability_iff_output_feedback_stability} and allows for
a set of approximate conditions to hold; the trade-off is that
it does not apply to every stabilizing controller for $(L, H)$, 
but instead requires the closed-loop $L+HK$ to have 
some bounded robustness margin, measured through its $H_\infty$-norm.

\begin{myprop}[Controller-dependent transfer of closed-loop stability]
\label{prop:controller_dependent_stability_transfer}
Suppose \Cref{assmp:stabilizability,assmp:non_trivial_unstable_subspace} hold,
and that
$E_T \in \R^{m \times n}$,
$L_T \in \R^{m \times m}$ and $H_T \in \R^{m \times p}$ satisfy
for some $\e\geq 0$ and $\sigma_0 > 0$,
\begin{align}
    \max\left\{
        \opnorm{P_{E_0^\T}P^\perp_{{E}_T^\T}},
        \opnorm{{E}_T A-{L}_T{E}_T},
        \opnorm{{E}_T B-{H}_T}
    \right\} \leq \e, \quad \sigma_{\min}({E}_T{E}_T^\T) \geq \sigma_0.
    \label{eq:controller_dependent_stability_transfer:requirement}
\end{align}
Let $K_T \in \R^{p\times m}$ be any controller such that
$L_T+H_TK_T$ is Schur stable, and define
\begin{align*}
    \Gamma_T := \sup_{z \in \mathbb{D}^c}
        \opnorm{\left(zI_m-L_T-H_TK_T\right)^{-1}}.
\end{align*}
If the following condition holds:
\begin{align}
    \e\left[\opnorm{\Pi_{+,0}}+\Gamma_T\left(\opnorm{K_T}
    +\opnorm{\Pi_{+,0}}\left(\sigma_0^{-1/2}
        +\Gamma_0\opnorm{B}\opnorm{K_T}\right)\right)\right]<1,
    \label{eq:controller_dependent_stability_transfer}
\end{align}
then $A+BK_TE_T$ is Schur stable, and
\begin{align}
    &\sup_{z \in \mathbb{D}^c}
        \opnorm{\left(zI_n-A-BK_TE_T\right)^{-1}} \nonumber \\
    &\qquad\qquad\leq
    \frac{\opnorm{\Pi_{+,0}}\left[\Gamma_0+
        \left(\sigma_0^{-1/2}+\Gamma_0\opnorm{B}\opnorm{K_T}\right)
        \Gamma_T\opnorm{E_T}\right]}
        {1-\e\left[\opnorm{\Pi_{+,0}}+\Gamma_T\left(\opnorm{K_T}
        +\opnorm{\Pi_{+,0}}\left(\sigma_0^{-1/2}
            +\Gamma_0\opnorm{B}\opnorm{K_T}\right)\right)\right]}.
    \label{eq:controller_dependent_H_infinity_bound}
\end{align}
\end{myprop}
\begin{proof}
From \Cref{prop:canonical_variables},
we know that
$E_0 A = L_0 E_0$ and $A D_0 = D_0 L_0$, which implies
that $A \Pi_{-,0} = \Pi_{-,0} A$.
Also, since $E_0 D_0 = I_{m_+}$, then we have
$E_0 \Pi_{-,0} = 0$. Hence, for any $z \in \C$,
\begin{align*}
    (zI_n-A-BK_0E_0)\Pi_{-,0}=\Pi_{-,0}(zI_n-A).
\end{align*}
From \Cref{def:canonical_variables_unstable},
the gain $K_0$ stabilizes $A+BK_0E_0$,
and hence for $z\in\mathbb{D}^c$ we can define:
\begin{align*}
    \mathcal{S}_0(z):=(zI_n-A-BK_0E_0)^{-1}\Pi_{-,0}, \qquad
    \opnorm{\mathcal{S}_0(z)}\leq\Gamma_0\opnorm{\Pi_{+,0}}.
\end{align*}
By the preceding identity, we have for all $z \in \mathbb{D}^c$,
\begin{align}
    \mathcal{S}_0(z)(zI_n-A)=\Pi_{-,0}. \label{eq:S0_projector_stable_identity}
\end{align}

Now we define:
\begin{align*}
    \Delta_A := E_T A - L_T E_T, \qquad \Delta_B := E_T B - H_T.
\end{align*}
For a fixed $z \in \mathbb{D}^c$ and unit-norm $v \in \C^n$, we observe:
\begin{align*}
    (z I_m - L_T - H_T K_T) E_T v = E_T (zI_n - A - BK_T E_T)v + \Delta_A v + \Delta_B K_T E_T v.
\end{align*}
Since $L_T + H_T K_T$ is assumed to be Schur stable,
$(z I_m - L_T - H_T K_T)$ is invertible, and hence:
\begin{align*}
    \norm{E_T v} &\leq \Gamma_T \left[ \opnorm{E_T} \norm{ (zI_n - A - BK_T E_T)v } + \opnorm{\Delta_A} + \opnorm{\Delta_B} \opnorm{K_T} \norm{E_Tv} \right] \\
    &\leq \Gamma_T \left[ \opnorm{E_T}\norm{ (zI_n - A - BK_T E_T)v }  + \e + \e \opnorm{K_T} \norm{E_Tv} \right] .
\end{align*}
Since \eqref{eq:controller_dependent_stability_transfer}
implies that $\e \Gamma_T \opnorm{K_T} < 1$, we can 
re-arrange terms to conclude:
\begin{align}
    \norm{E_T v} \leq \frac{\Gamma_T \left[\opnorm{E_T} \norm{ (zI_n - A - BK_T E_T)v } + \e \right]}{ 1 - \e \Gamma_T \opnorm{K_T}}. \label{eq:controller_dependent_stability_transfer:ET_v_bound}
\end{align}

The assumption $\sigma_{\min}(E_TE_T^\T) \geq \sigma_0$ implies that $E_T$ has full row rank and
$\opnorm{E_T^\dag}\leq\sigma_0^{-1/2}$;
here $(\cdot)^\dag$ denotes the Moore–Penrose pseudo-inverse of a matrix.
Also because $\Pi_{+,0}=D_0E_0$, we have
$\Pi_{+,0}P_{E_0^\T}=\Pi_{+,0}$.
Hence,
\begin{align*}
    \Pi_{+,0} = \Pi_{+,0} P_{E_T^\T} + \Pi_{+,0} P_{E_T^\T}^\perp = \Pi_{+,0} P_{E_T^\T} + \Pi_{+,0} P_{E_0^\T} P_{E_T^\T}^\perp = \Pi_{+,0} E_T^\dag E_T + \Pi_{+,0} P_{E_0^\T} P_{E_T^\T}^\perp ,
\end{align*}
where the last equality holds since $P_{M^\T} = M^\dag M$ for any matrix $M$.
Therefore,
\begin{align*}
    \norm{ \Pi_{+,0} v } &= \norm{ \Pi_{+,0} E_T^\dag E_Tv + \Pi_{+,0} P_{E_0^\T} P_{E_T^\T}^\perp v } \\
    &\leq \norm{ \Pi_{+,0} E_T^\dag E_Tv } + \norm{ \Pi_{+,0} P_{E_0^\T} P_{E_T^\T}^\perp v } \\
    &\leq \opnorm{ \Pi_{+,0} } \opnorm{ E_T^\dag } \norm{ E_T v } + \opnorm{ \Pi_{+,0} } \opnorm{ P_{E_0^\T} P_{E_T^\T}^\perp } \\
    &\leq \opnorm{ \Pi_{+,0} } \left( \sigma_0^{-1/2} \norm{ E_T v }  + \e \right).
\end{align*}
Now using \eqref{eq:S0_projector_stable_identity},
\begin{align*}
    \Pi_{-,0} v = \calS_0(z) (zI_n - A) v = \calS_0(z)\left[ (zI_n - A - B K_T E_T)v + B K_T E_T v \right],
\end{align*}
and therefore:
\begin{align*}
    \norm{ \Pi_{-,0} v } \leq \Gamma_0 \opnorm{\Pi_{+,0}} \left[ \norm{ (zI_n - A - BK_T E_T)v } + \opnorm{B} \opnorm{K_T} \norm{E_T v} \right] .
\end{align*}
Using $\Pi_{+,0}+\Pi_{-,0}=I_n$, we conclude that
\begin{align}
    1&\leq\norm{\Pi_{+,0}v}+\norm{\Pi_{-,0}v} \nonumber \\
    &\leq\opnorm{\Pi_{+,0}}\left[\e
        +\left(\sigma_0^{-1/2}+\Gamma_0\opnorm{B}\opnorm{K_T}\right)\norm{E_T v}
        +\Gamma_0\norm{ (zI_n - A - BK_T E_T)v }\right]. \label{eq:controller_dependent_stability_transfer:intermediate}
\end{align}

To finish the proof, substitute our bound on 
$\norm{ E_T v }$ from \eqref{eq:controller_dependent_stability_transfer:ET_v_bound}
into \eqref{eq:controller_dependent_stability_transfer:intermediate}, and then use the 
assumption \eqref{eq:controller_dependent_stability_transfer} to
re-arrange the inequality to lower bound 
$\norm{ (zI_n - A - BK_T E_T)v }$ as:
\begin{align*}
    \norm{ (zI_n - A - BK_T E_T)v } \geq
    \frac{1-\e\left[\opnorm{\Pi_{+,0}}+\Gamma_T\left(\opnorm{K_T}
        +\opnorm{\Pi_{+,0}}\left(\sigma_0^{-1/2}
            +\Gamma_0\opnorm{B}\opnorm{K_T}\right)\right)\right]}
        {\opnorm{\Pi_{+,0}}\left[\Gamma_0+
        \left(\sigma_0^{-1/2}+\Gamma_0\opnorm{B}\opnorm{K_T}\right)
        \Gamma_T\opnorm{E_T}\right]}.
\end{align*}
Here we discarded the nonnegative numerator term $\e^2\opnorm{\Pi_{+,0}}\Gamma_T\opnorm{K_T}$, and used $1-\e\Gamma_T\opnorm{K_T}\leq1$ to upper bound the denominator.
Since this bound holds for all $z\in\mathbb{D}^c$ and unit-norm $v\in\C^n$, $zI_n-A-BK_TE_T$ is invertible throughout $\mathbb{D}^c$.
Thus $A+BK_TE_T$ is Schur stable, and taking reciprocals gives \eqref{eq:controller_dependent_H_infinity_bound}.
\end{proof}

\Cref{prop:controller_dependent_stability_transfer} generalizes
the forward direction of \Cref{prop:latent_stability_iff_output_feedback_stability} in the following sense. Suppose that the approximation error
$\e$ in \eqref{eq:controller_dependent_stability_transfer:requirement}
can be taken as $\e=0$,
which together with the bound $\sigma_{\min}(E_TE_T^\T) \geq \sigma_0 > 0$, implies that conditions (i)--(iii) of \Cref{prop:latent_stability_iff_output_feedback_stability} hold.
Then the condition \eqref{eq:controller_dependent_stability_transfer}
immediately holds for \emph{any} stabilizing controller $K_T$ for $(L_T, H_T)$, and by the conclusion of \Cref{prop:controller_dependent_stability_transfer} we then have that
$A + B K_T E_T$ is Schur stable.

With \Cref{prop:controller_dependent_stability_transfer} in hand, the following result
(which is the full version of \Cref{thm:resolvent_gain_bound}) chooses one particular controller $K_T$---namely the minimum $H_\infty$-norm latent controller---and estimates the amount of approximation error $\e$ in \eqref{eq:controller_dependent_stability_transfer:requirement} which is tolerable for 
both $K_T$ to exist with finite $H_\infty$-norm,
and to also satisfy the conditions in 
\eqref{eq:controller_dependent_stability_transfer},
and hence guarantee that $A + B K_T E_T$ is Schur stable.

\begin{mythm}[Minimum $H_\infty$-norm closed-loop stability]
\label{thm:H_infty_closed_loop_stability}
Suppose \Cref{assmp:stabilizability,assmp:non_trivial_unstable_subspace} hold,
and that
$E_T \in \R^{m \times n}$,
$L_T \in \R^{m \times m}$ and $H_T \in \R^{m \times p}$ satisfy
for some $\e\geq 0$ and $\sigma_0 > 0$,
\begin{align*}
    \max\left\{
        \opnorm{P_{E_0^\T}P^\perp_{{E}_T^\T}},
        \opnorm{{E}_T A-{L}_T{E}_T},
        \opnorm{{E}_T B-{H}_T}
    \right\} \leq \e, \quad \sigma_{\min}({E}_T{E}_T^\T) \geq \sigma_0.
\end{align*}
Fix an $M_0\geq\opnorm{K_0E_0}/\sqrt{\sigma_0}$,
and suppose that:
\begin{align}
    \e \leq \frac12\min\Bigg\{
        &\frac{\sqrt{\sigma_0}}
            {\Gamma_0\left(\sqrt{\sigma_0}\opnorm{BK_0E_0}
                +1+\opnorm{K_0E_0}\right)}, \notag\\
        &\frac{\sqrt{\sigma_0}}
            {\sqrt{\sigma_0}\opnorm{\Pi_{+,0}}
                +2\Gamma_0\opnorm{E_T}
                \left[M_0+\opnorm{\Pi_{+,0}}\left(\sigma_0^{-1/2}
                    +\Gamma_0\opnorm{B}M_0\right)\right]}
        \Bigg\}. \label{eq:H_infty_closed_loop_stability:epsilon_condition}
\end{align}
Then the following conditions hold:
\begin{enumerate}[label=(\roman*)]
    \item The pair $(L_T, H_T)$ is stabilizable.
    \item The minimum $H_\infty$-norm latent controller
    \begin{align*}
        K_T \in \argmin_{K \in B_{p \times m}^{(\mathrm{op})}(M_0)}
        \sup_{z \in \mathbb{D}^c}
        \opnorm{(zI_m-L_T-H_TK)^{-1}}
    \end{align*}
    exists,
    and every minimizer makes both $L_T+H_TK_T$ and $A+BK_TE_T$ Schur stable
    and satisfies:
    \begin{align}
        \sup_{z \in \mathbb{D}^c}
        \opnorm{(zI_n-A-BK_TE_T)^{-1}} \leq 2\Gamma_0\opnorm{\Pi_{+,0}}
        \left[1+\frac{2\opnorm{E_T}^2}{\sqrt{\sigma_0}}
            \left(\sigma_0^{-1/2}+\Gamma_0\opnorm{B}M_0\right)\right]. \label{eq:H_infty_closed_loop_stability:gain}
    \end{align}
\end{enumerate}
\end{mythm}
\begin{proof}
Define the following matrices:
\begin{align*}
    \widehat K_T:=K_0E_0E_T^\dag, \qquad
    \widehat A_T:=A+B\widehat K_TE_T, \qquad
    \widehat L_T:=L_T+H_T\widehat K_T.
\end{align*}
The assumption
$\sigma_{\min}(E_TE_T^\T) \geq \sigma_0$ yields
$\opnorm{E_T^\dag}\leq \sigma_0^{-1/2}$. Consequently, using our assumption on $M_0$,
\begin{align*}
    \opnorm{\widehat K_T}\leq\opnorm{K_0E_0}/\sqrt{\sigma_0}\leq M_0, \qquad
    \opnorm{\widehat K_TE_T}\leq\opnorm{K_0E_0}.
\end{align*}
Let $A_0:=A+BK_0E_0$ and $\Delta:=\widehat A_T-A_0$.
Since
$P_{E_T^\T} = E_T^\dag E_T$ and 
$E_0P_{E_0^\T}=E_0$,
\begin{align*}
    \Delta
    =-BK_0E_0 P^\perp_{E_T^\T}
    =-BK_0E_0P_{E_0^\T}P_{E_T^\T}^{\perp} \implies \opnorm{\Delta}
    \leq\e\opnorm{BK_0E_0},
\end{align*}
where the implication uses
our assumption $\opnorm{ P_{E_0^\T}P_{E_T^\T}^{\perp} } \leq \e$.

Set $\psi:=1-\e\Gamma_0\opnorm{BK_0E_0}$; by our assumption 
\eqref{eq:H_infty_closed_loop_stability:epsilon_condition}
we have $\psi\geq1/2$.
Since $A_0$ is Schur stable, for every $z\in\mathbb{D}^c$ the resolvent
$R_0(z):=(zI_n-A_0)^{-1}$ exists and satisfies $\opnorm{R_0(z)}\leq\Gamma_0$.
We factorize:
\begin{align*}
    zI_n-\widehat A_T
    =(zI_n-A_0)\left(I_n-R_0(z)\Delta\right).
\end{align*}
The perturbation bound on $\opnorm{\Delta}$ yields:
\begin{align*}
    \opnorm{R_0(z)\Delta}
    \leq\e\Gamma_0\opnorm{BK_0E_0}
    =1-\psi\leq\frac12<1.
\end{align*}
Therefore, the Neumann series
$\left(I_n-R_0(z)\Delta\right)^{-1}
    =\sum_{j=0}^{\infty}\left(R_0(z)\Delta\right)^j$
converges
with the norm bound:
\begin{align*}
    \opnorm{(I_n-R_0(z)\Delta)^{-1}} \leq \frac{1}{1-\opnorm{R_0(z)\Delta}} \leq \frac{1}{\psi}.
\end{align*}
It follows that $zI_n-\widehat A_T$ is invertible, with
\begin{align*}
    \widehat R_T(z):=(zI_n-\widehat A_T)^{-1}
    =\left(I_n-R_0(z)\Delta\right)^{-1}R_0(z),
    \qquad \opnorm{\widehat R_T(z)}\leq\frac{\Gamma_0}{\psi}.
\end{align*}
Since this holds for every $z\in\mathbb{D}^c$, $\widehat A_T$ is Schur stable.

Let $\Delta_A:=E_TA-L_TE_T$ and $\Delta_B:=E_TB-H_T$.
With these variables we can write:
\begin{align*}
    (zI_m-\widehat L_T)E_T
    =E_T(zI_n-\widehat A_T)+\Delta_A+\Delta_B\widehat K_TE_T.
\end{align*}
Since $E_T$ has full row rank, we have
$E_T E_T^\dag = I_m$. Hence,
multiplying on the right by $\widehat R_T(z)E_T^\dag$ yields
\begin{align*}
    (zI_m-\widehat L_T)E_T\widehat R_T(z)E_T^\dag
    =I_m+ \underbrace{(\Delta_A+\Delta_B\widehat K_TE_T)  \widehat R_T(z)E_T^\dag}_{=: J_T(z)}, \qquad 
    \opnorm{J_T(z)} \leq
        \frac{\e\Gamma_0(1+\opnorm{K_0E_0})}{\sqrt{\sigma_0}\,\psi}.
\end{align*}
The assumption \eqref{eq:H_infty_closed_loop_stability:epsilon_condition} implies:
\begin{align*}
    \psi-\frac{\e\Gamma_0(1+\opnorm{K_0E_0})}{\sqrt{\sigma_0}}
    =1-\e\Gamma_0\left(\opnorm{BK_0E_0}
        +\frac{1+\opnorm{K_0E_0}}{\sqrt{\sigma_0}}\right)\geq\frac12 \implies \opnorm{J_T(z)}<1.
\end{align*}
Thus 
$I_m + J_T(z)$ is invertible
with norm bound $\opnorm{(I_m + J_T(z))^{-1}} \leq 1/(1 - \opnorm{J_T(z)})$.
This implies that
$zI_m-\widehat L_T$ is invertible
with inverse $E_T\widehat R_T(z)E_T^\dag(I_m+J_T(z))^{-1}$.
Consequently,
\begin{align*}
    \sup_{z \in \mathbb{D}^c}\opnorm{(zI_m-\widehat L_T)^{-1}}
    \leq\frac{\Gamma_0\opnorm{E_T}}
        {\sqrt{\sigma_0}\left(\psi-\e\Gamma_0(1+\opnorm{K_0E_0})/\sqrt{\sigma_0}\right)}
    \leq\frac{2\Gamma_0\opnorm{E_T}}{\sqrt{\sigma_0}}.
\end{align*}
This proves (i) and exhibits a feasible controller with finite objective in (ii).

For each fixed $z \in \mathbb{D}^c$, the extended resolvent norm is
$1/\sigma_{\min}(zI_m-L_T-H_TK)$, with $1/0:=+\infty$,
and is lower semicontinuous in $K$.
Its supremum over $z \in \mathbb{D}^c$ is therefore lower semicontinuous as well.
Compactness of 
$B_{p \times m}^{(\mathrm{op})}(M_0)$ implies that a minimizer exists, and the construction of $\widehat{K}_T \in B_{p \times m}^{(\mathrm{op})}(M_0)$
shows that all minimizers $K_T$ attain a finite gain, and are thus Schur stabilizing
with gain bounded by:
\begin{align*}
    \Gamma_T:=\sup_{z \in \mathbb{D}^c}
        \opnorm{(zI_m-L_T-H_TK_T)^{-1}}
        \leq \sup_{z \in \mathbb{D}^c} \opnorm{(zI_m - \widehat{L}_T)^{-1}} 
        \leq\frac{2\Gamma_0\opnorm{E_T}}{\sqrt{\sigma_0}}.
\end{align*}
To conclude the proof, we apply
\Cref{prop:controller_dependent_stability_transfer}
to the tuple $(E_T, L_T, H_T, K_T)$.
The assumption \eqref{eq:H_infty_closed_loop_stability:epsilon_condition} ensures that:
\begin{align*}
    &\e\left[\opnorm{\Pi_{+,0}}+\Gamma_T\left(
        \opnorm{K_T}+\opnorm{\Pi_{+,0}}\left(\sigma_0^{-1/2}
            +\Gamma_0\opnorm{B}\opnorm{K_T}\right)\right)\right] \\
    &\quad\leq\e\left[\opnorm{\Pi_{+,0}}
        +\frac{2\Gamma_0\opnorm{E_T}}{\sqrt{\sigma_0}}\left(
        M_0+\opnorm{\Pi_{+,0}}\left(\sigma_0^{-1/2}
            +\Gamma_0\opnorm{B}M_0\right)\right)\right]
        \leq\frac12<1,
\end{align*}
which satisfies the requirement 
\eqref{eq:controller_dependent_stability_transfer} in \Cref{prop:controller_dependent_stability_transfer}.
We then conclude that
$A + B K_T E_T$ is Schur stable
and the gain bound 
\eqref{eq:H_infty_closed_loop_stability:gain} holds, proving (ii).
\end{proof}

\section{Analysis of Multi-Step Backwards Loss Optimization}

We start by introducing an idealized multi-step backward loss 
for analysis purposes.

\begin{mydef}[Unconstrained multi-step backward loss]\label{df:uncon_loss}
    We define the unconstrained backward loss to be the loss in \Cref{def:mstep_backwards_conjugacy_loss_general} but with \emph{unconstrained prediction functions} $z_t : \R^m \times \R^{pt} \mapsto \R^m$ replacing the state transition function $\Phi_{f_z}^t$ defined in \Cref{eq:mstep_forward_conjugacy_loss_general}:
    \begin{gather}
        J_\mathrm{uncon}^{(T)}(E, D) := \min_{z_0, \dots, z_{T}}\frac{1}{T+1}\sum_{t=0}^{T}\mathbb{E}_{x_0,\mathbf{u}_{0:t-1}}\left[ \norm{\Phi_f(x_0,\mathbf{u}_{0:t-1})-D \circ z_t(E (x_0), \mathbf{u}_{0:t-1})}^2 \right].
    \end{gather}
\end{mydef}

\begin{myfact}[Unconstrained loss instantiated for linear encoders, decoders, and dynamics]
\label{prop:unconstrained_loss}
Let the encoder $E(x) = Ex$ for $E \in \R^{m \times n}$, 
decoder $D(x) = Dx$ for $D \in \R^{n \times m}$.
Under the LTI dynamics \eqref{eq:LDS}
and distribution for $(x_0, \bm{u}_{0:T-1})$ 
specified in \Cref{assmp:data_model},
    \begin{align*}
        J_{\mathrm{uncon}}^{(T)}(E, D) &= \frac{\sigma_x^2}{T+1}\sum_{t=0}^{T} \norm{ P_D A^t P^\perp_{E^\T} }_F^2 + \frac{\sigma_x^2}{T+1} \sum_{t=0}^{T} \norm{P^\perp_D A^t}_F^2 + \frac{\sigma_u^2}{T+1} \sum_{k=0}^{T-1} \beta_{T,k} \norm{P^\perp_D A^k B}_F^2 \\
        &= \frac{\sigma_x^2}{T+1}\sum_{t=0}^{T} \norm{ A^t P^\perp_{E^\T} }_F^2 + \frac{\sigma_x^2}{T+1} \sum_{t=0}^{T} \norm{ P^\perp_{D} A^t P_{E^\T} }_F^2 + \frac{\sigma_u^2}{T+1} \sum_{k=0}^{T-1} \beta_{T,k} \norm{P^\perp_D A^k B}_F^2 ,
    \end{align*}
where $\beta_T(k)$ is as defined in \Cref{prop:LTI_system_losses}.
\end{myfact}
The unconstrained loss plays a crucial part of our analysis,
as it removes the dependence on the variables $(L_T, H_T)$
in $J_{\mathrm{bwd}}^{(T)}$, where
the former which enters through repeated matrix powers $L_T^k$.
At the same time, the unconstrained loss
relates to $J_{\mathrm{bwd}}^{(T)}$ via the
inequality
$J_{\mathrm{uncon}}^{(T)}(E, D) \leq J_{\mathrm{bwd}}^{(T)}(E, D, L, H)$ for any $(E, D, L, H)$.

Our final assumption concerns the variance of the inputs $\bm{u}_{0:T-1}$ entering the system.
As our analysis relies on inputs and initial condition propagation
over time, we need to ensure that the signal of the later inputs
does not drown out the beginning inputs. 
This is formalized through the following weighted sum.

\begin{mydef}[Definition of weighted geometric sum $S_T$]
\label{def:S_T_decay}
For $\gamma \in (0, 1)$ and $T \in \N_+$, we define
the weighted geometric sum
$S_T(\gamma) := \sum_{k=0}^{T-1} \beta_T(k) \gamma^{2k}$,
where $\beta_T(k)$ is as defined in \Cref{prop:LTI_system_losses}.
\end{mydef}

In a moment, we will impose decay conditions on $S_T$.
First though, we utilize this definition to evaluate the loss of the canonical variables $(E_m, D_m, L_m, H_m)$
from \Cref{def:canonical_variables_latent} on
the backward loss $\calL_{\mathrm{bwd}}^{(T)}$.

\begin{myprop}[Loss of canonical variables]
\label{prop:canonical_variable_loss}
Suppose that
\Cref{assmp:non_trivial_unstable_subspace,assmp:latent_dimension_compatibility,assmp:regularizer_assumption} hold.
Let the tuple $(E_m, D_m, L_m, H_m)$ denote the canonical variables 
from \Cref{def:canonical_variables_latent}.
For every $T\in\N_+$ and $\gamma\in(\gamma_-,1)$:
\begin{align}
    J_{\mathrm{bwd}}^{(T)}(E_m, D_m, L_m, H_m)
    &\leq \frac{C^2(A,\gamma)}{T+1}
        \left[\frac{\sigma_x^2(n-m)}{1-\gamma^2}
            +\sigma_u^2\norm{B}_F^2 S_T(\gamma)\right], \\
    J_{\mathrm{fwd}}^{(T)}(E_m, L_m, H_m) &= 0, \qquad R(E_m) = 0.
\end{align}
\end{myprop}
\begin{proof}
By \Cref{prop:canonical_variables}, for every $k\geq0$,
\begin{align*}
    E_mA^k=L_m^kE_m, \qquad H_m=E_mB, \qquad
    D_mL_m^kE_m=\Pi_{+,m}A^k.
\end{align*}
Consequently, the reconstruction residuals are
\begin{align*}
    A^k-D_mL_m^kE_m=\Pi_{-,m}A^k, \qquad 
    A^kB-D_mL_m^kH_m=\Pi_{-,m}A^kB.
\end{align*}
The identity $E_mD_m=I_m$ implies that $D_mE_m$ is a rank-$m$ projector,
so $\Pi_{-,m}=I_n-D_mE_m$ has rank $n-m$.
Using the decay estimate from \Cref{prop:canonical_variables}, we obtain
\begin{align*}
    \norm{\Pi_{-,m}A^k}_F^2
    &\leq(n-m)\opnorm{\Pi_{-,m}A^k}^2
    \leq(n-m)C^2(A,\gamma)\gamma^{2k}, \\
    \norm{\Pi_{-,m}A^kB}_F^2
    &\leq\opnorm{\Pi_{-,m}A^k}^2\norm{B}_F^2
    \leq C^2(A,\gamma)\gamma^{2k}\norm{B}_F^2.
\end{align*}
Substituting into the backwards loss formula in \Cref{prop:LTI_system_losses} gives
\begin{align*}
    J_{\mathrm{bwd}}^{(T)}(E_m,D_m,L_m,H_m)
    &=\frac{\sigma_x^2}{T+1}\sum_{k=0}^T\norm{\Pi_{-,m}A^k}_F^2
        +\frac{\sigma_u^2}{T+1}\sum_{k=0}^{T-1}
            \beta_T(k)\norm{\Pi_{-,m}A^kB}_F^2 \\
    &\leq\frac{C^2(A,\gamma)}{T+1}
        \left[\sigma_x^2(n-m)\sum_{k=0}^T\gamma^{2k}
            +\sigma_u^2\norm{B}_F^2
                \sum_{k=0}^{T-1}\beta_T(k)\gamma^{2k}\right] \\
    &\leq\frac{C^2(A,\gamma)}{T+1}
        \left[\frac{\sigma_x^2(n-m)}{1-\gamma^2}
            +\sigma_u^2\norm{B}_F^2 S_T(\gamma)\right].
\end{align*}
The last step uses the 
definition of $S_T(\gamma)$ from \Cref{def:S_T_decay}.
Next, $E_mA=L_mE_m$ and $E_mB=H_m$ implies via \Cref{prop:LTI_system_losses} that
$J_{\mathrm{fwd}}^{(T)}(E_m, L_m, H_m) = 0$.
Finally,
$E_m$ consists of rows of an orthogonal matrix (cf.~\eqref{eq:canonical_variables} from \Cref{def:canonical_variables_v2}),
and hence regularizer admissibility (cf.~\Cref{def:forward_admissible_regularizer}) yields $R(E_m)=0$.
\end{proof}

With \Cref{prop:canonical_variable_loss} in place, we are now ready to state and prove our main recovery bounds for minimizers of the multi-step backward loss.
The following result is the 
detailed version of \Cref{thm:informal_recovery}.

\begin{mythm}[Recovery bounds for multi-step backward loss minimization]
\label{thm:multistep_backwards_full}
Suppose that
\Cref{assmp:stabilizability,assmp:non_trivial_unstable_subspace,assmp:latent_dimension_compatibility,assmp:regularizer_assumption} hold.
Any minimizer $(\bar{E}_T, \bar{D}_T, \bar{L}_T, \bar{H}_T)$ of
$\calL_{\mathrm{bwd}}^{(T)}$ (cf.~\Cref{eq:bwd_solution})
satisfies, for any $\gamma \in (\gamma_-, 1)$:
\begin{align*}
    \norm{ P_{E_0^\T} P^\perp_{\bar{E}_T^\T } }_F^2
        &\leq \frac{C_T}{\sigma_x^2\eta_{T+1}}, \quad C_T:=\frac{C^2(A,\gamma)}{T+1}
        \left[\frac{\sigma_x^2(n-m)}{1-\gamma^2}
             +\sigma_u^2\norm{B}_F^2 S_T(\gamma)\right], \\
    \norm{\bar{E}_TA-\bar{L}_T\bar{E}_T}_F^2
        &\leq \frac{TC_T}{\lambda_{\mathrm{fwd}}\sigma_x^2}, \qquad
    \norm{\bar{E}_TB-\bar{H}_T}_F^2
        \leq \frac{TC_T}{\lambda_{\mathrm{fwd}}\sigma_u^2\beta_T(0)},
\end{align*}
where $\eta_{T+1} := \lambda_{\min}\left( \frac{1}{T+1} \sum_{t=0}^{T} (L_0^\T)^t L_0^t \right) > 0$.
Furthermore, if $ \lambda_{\mathrm{reg}}\geq C_T/\e_0$ then we also have:
\begin{align*}
    \sigma_0 I_m
        &\preccurlyeq \bar{E}_T\bar{E}_T^\T \preccurlyeq \sigma_1 I_m.
\end{align*}
\end{mythm}
\begin{proof}
The canonical variables $(E_m,D_m,L_m,H_m)$ covering all latent modes from \Cref{def:canonical_variables_latent} are feasible for the optimization in
\Cref{eq:bwd_solution}. 
Hence by optimality 
of $(\bar{E}_T, \bar{D}_T, \bar{L}_T, \bar{H}_T)$ and \Cref{prop:canonical_variable_loss},
\begin{align*}
    \calL_{\mathrm{bwd}}^{(T)}(\bar{E}_T,\bar{D}_T,\bar{L}_T,\bar{H}_T)
    \leq\calL_{\mathrm{bwd}}^{(T)}(E_m,D_m,L_m,H_m) 
    =J_{\mathrm{bwd}}^{(T)}(E_m,D_m,L_m,H_m)\leq C_T.
\end{align*}
Since all three terms in
$\calL_{\mathrm{bwd}}^{(T)}(\bar{E}_T,\bar{D}_T,\bar{L}_T,\bar{H}_T)$  are nonnegative,
$C_T$ therefore bounds each term.

To prove the subspace estimate, we  utilize the unconstrained loss
(cf.~\Cref{df:uncon_loss}) as a point of comparison.
Specifically, using $J_{\mathrm{uncon}}^{(T)}(E, D) \leq J_{\mathrm{bwd}}^{(T)}(E, D, L, H)$ for any $(E, D, L, H)$ combined with \Cref{prop:unconstrained_loss}:
\begin{align*}
    C_T
    &\geq J_{\mathrm{bwd}}^{(T)}(\bar{E}_T,\bar{D}_T,\bar{L}_T,\bar{H}_T)
    \geq J_{\mathrm{uncon}}^{(T)}(\bar{E}_T,\bar{D}_T)
    \geq\frac{\sigma_x^2}{T+1}\sum_{t=0}^T\norm{A^tP^\perp_{\bar{E}_T^\T}}_F^2.
\end{align*}
Now the subset $(E_0, L_0, H_0)$ of canonical variables for the non-stable modes (cf.~\Cref{def:canonical_variables_unstable})
satisfies
$P_{E_0^\T} = E_0^\T E_0$,
$E_0E_0^\T = I_{m_+}$, 
and $E_0A^k=L_0^kE_0$ for all $k \geq 0$.
Therefore:
\begin{align*}
    \frac{1}{T+1}\sum_{t=0}^T\norm{A^tP^\perp_{\bar{E}_T^\T}}_F^2
    &\geq\frac{1}{T+1}\sum_{t=0}^T\norm{P_{E_0^\T} A^tP^\perp_{\bar{E}_T^\T}}_F^2 
    =\frac{1}{T+1}\sum_{t=0}^T\norm{E_0A^tP^\perp_{\bar{E}_T^\T}}_F^2 \\
    &=\frac{1}{T+1}\sum_{t=0}^T\norm{L_0^tE_0P^\perp_{\bar{E}_T^\T}}_F^2
    \geq\eta_{T+1}\norm{E_0P^\perp_{\bar{E}_T^\T}}_F^2 = \eta_{T+1}\norm{ P_{E_0^\T} P^\perp_{\bar{E}_T^\T}}_F^2.
\end{align*}
Combining these inequalities
shows that 
$\norm{ P_{E_0^\T} P^\perp_{\bar{E}_T^\T}}_F^2 \leq \frac{C_T}{\sigma_x^2 \eta_{T+1}}$.

Next, evaluating the $T$-step forward loss gives
(cf.~\Cref{prop:LTI_system_losses}):
\begin{align*}
   \lambda_{\mathrm{fwd}}
        J_{\mathrm{fwd}}^{(T)}(\bar{E}_T,\bar{L}_T,\bar{H}_T)
        =
\lambda_{\mathrm{fwd}}\frac{\sigma_x^2}{T}\sum_{t=1}^{T} \norm{\bar{E}_T A^t - \bar{L}_T^t \bar{E}_T}_F^2 + \lambda_{\mathrm{fwd}}\frac{\sigma_u^2}{T} \sum_{k=0}^{T-1} \beta_{T}(k) \norm{ \bar{E}_T A^k B - \bar{L}_T^k \bar{H}_T }_F^2
        \leq C_T,
\end{align*}
from which we conclude:
\begin{align*}
    \norm{\bar{E}_TA-\bar{L}_T\bar{E}_T}_F^2 \leq \frac{TC_T}{\lambda_{\mathrm{fwd}} \sigma_x^2}, \qquad \norm{\bar{E}_TB-\bar{H}_T}_F^2 \leq \frac{TC_T}{\lambda_{\mathrm{fwd}} \sigma_u^2 \beta_T(0)}.
\end{align*}

Finally, since $R(\bar{E}_T)\leq\frac{C_T}{\lambda_{\mathrm{reg}}}$
if $\lambda_{\mathrm{reg}}\geq C_T/\e_0$, then $R(\bar{E}_T)\leq\e_0$,
and condition (ii) of \Cref{def:forward_admissible_regularizer} yields
$\sigma_0 I_m\preccurlyeq\bar{E}_T\bar{E}_T^\T\preccurlyeq\sigma_1 I_m$.
\end{proof}

\begin{rmk}[Uniform lower bound of $\eta_T$]
\label{remark:uniform_lower_bound_eta_T}
We remark that the quantity
$\eta_{T+1} := \lambda_{\min}\left( \frac{1}{T+1} \sum_{t=0}^{T} (L_0^\T)^t L_0^t \right)$ which appears in 
\Cref{thm:multistep_backwards_full}
actually satisfies the following
time-independent \emph{uniform lower bound}:
\begin{align}
    \inf_{T \in \N_+} \eta_{T+1} := \eta_0 > 0.
\end{align}
This is a consequence of the canonical
variable construction \Cref{def:canonical_variables_unstable}, since the eigenvalues of $L_0$ are all at least
marginally stable; see \Cref{lemma:nonstable_gramian}.
Furthermore, whenever all the eigenvalues of 
$\sigma_+(A)$ are strictly unstable,
then actually $\eta_{T}$ actually grows
\emph{exponentially} in $T$, 
which further improves the error estimates
supplied in \Cref{thm:multistep_backwards_full}.
\end{rmk}

Our next objective is to 
combine the error estimates from
\Cref{thm:multistep_backwards_full}
with the stability transfer guarantees
from \Cref{thm:H_infty_closed_loop_stability},
providing an end-to-end bound 
ensuring closed-loop stability of the control design.
To do this, we need to impose some regularity
structure on the input variances from 
\Cref{assmp:data_model}. We do this by imposing
a regularity condition on the weighted geometric sum
$S_T$ defined in \Cref{def:S_T_decay}.

\begin{assmp}[Sufficient decay of $S_T$]
\label{assmp:S_T_decay}
We assume that for any $\gamma \in (0, 1)$, 
$S_T$ from \Cref{def:S_T_decay} satisfies
$\lim_{T \to \infty} S_T(\gamma)/T = 0$.
From this, we define $\Phi(\e, \gamma) := \min\{ T_0 \in \N_+ \mid S_T(\gamma)/T \leq \e, \,\, \forall T \geq T_0 \}$.
\end{assmp}

\begin{rmk}[Examples of admissible $\alpha_{T}(k)$ sequences]
\label{remark:admissible_sequences}
Some examples of valid $\alpha_{T}(k)$ sequences for the input variance of $u_k$
which satisfy \Cref{assmp:S_T_decay}
include the following:
(a) \emph{Constant variance}: We can set $\alpha_{T}(k) = 1/T$,
in which case 
$\beta_{T}(k) = (T-k)/T$,
$S_T(\gamma)/T = T^{-1} \sum_{k=0}^{T-1} (1-k/T) \gamma^{2k} \leq 1/(T (1-\gamma^2))$, and
$\Phi(\e, \gamma) \leq \ceil{1/(\e(1-\gamma^2))}$.
(b) \emph{Decaying variance}: We can set $\alpha_{T}(k) = 1/(k+1)^\eta$ for $\eta > 0$
in which case $\beta_{T}(k) = \sum_{j=1}^{T-k} j^{-\eta} := H^{(\eta)}_{T-k}$ and $S_T(\gamma)/T \leq H^{(\eta)}_{T}/(T(1-\gamma^2))$.
By the standard asymptotic scaling of generalized harmonic numbers:
\begin{align*}
    S_T(\gamma)/T \lesssim \begin{cases}
        \frac{1}{(1-\eta)(1-\gamma^2) T^\eta} & \text{if } \eta \in (0, 1) \\
        \frac{\log{T}}{T(1-\gamma^2)} &\text{if } \eta = 1 \\
        \frac{\zeta(\eta)}{T(1-\gamma^2)} &\text{if } \eta > 1,
    \end{cases} \implies \Phi(\e, \gamma) \lesssim \begin{cases}
        [ \e(1-\eta)(1-\gamma^2) ]^{-1/\eta} &\text{if } \eta \in (0, 1) \\
        [ \e(1-\gamma^2) ]^{-1} \log(1/[ \e (1-\gamma^2) ]) &\text{if } \eta = 1 \\
        \zeta(\eta)/[\e (1-\gamma^2)] &\text{if } \eta > 1,
    \end{cases}
\end{align*}
where $\zeta(s) := \sum_{n=1}^{\infty} n^{-s}$ is the Riemann zeta function.
\end{rmk}

With the decay assumption in place, the following result
is our main end-to-end guarantee, and provides the
detailed result for \Cref{cor:latent_stability_transfer}.

\begin{mycorollary}[End-to-end closed-loop guarantees]
\label{thm:end_to_end_closed_loop}
Suppose that
\Cref{assmp:stabilizability,assmp:non_trivial_unstable_subspace,assmp:latent_dimension_compatibility,assmp:regularizer_assumption,assmp:S_T_decay} hold.
Fix any $\gamma \in (\gamma_-, 1)$ and
$M_0\geq\opnorm{K_0E_0}/\sqrt{\sigma_0}$, and define:
\begin{align}
    \e_{\mathrm{stab}}:=\frac12\min\Bigg\{
        &\frac{\sqrt{\sigma_0}}
            {\Gamma_0\left(\sqrt{\sigma_0}\opnorm{BK_0E_0}
                +1+\opnorm{K_0E_0}\right)}, \notag\\
        &\frac{\sqrt{\sigma_0}}
            {\sqrt{\sigma_0}\opnorm{\Pi_{+,0}}
                +2\Gamma_0\sqrt{\sigma_1}
                \left[M_0+\opnorm{\Pi_{+,0}}\left(\sigma_0^{-1/2}
                    +\Gamma_0\opnorm{B}M_0\right)\right]}
        \Bigg\}.
\end{align}
Let $\eta_0>0$ be as in \Cref{remark:uniform_lower_bound_eta_T},
and let $C_T$ be as in \Cref{thm:multistep_backwards_full}.
If the following conditions on the
loss $\calL_{\mathrm{bwd}}^{(T)}$ hyperparameters hold:
\begin{align}
    T &\geq \max\left\{
        \frac{2C^2(A,\gamma)(n-m)}{(1-\gamma^2)\eta_0\e_{\mathrm{stab}}^2},\,
        \Phi\left(\frac{\sigma_x^2\eta_0\e_{\mathrm{stab}}^2}
            {2C^2(A,\gamma)\sigma_u^2\norm{B}_F^2},\gamma\right)
        \right\}, \label{eq:end_to_end_horizon}\\
    \lambda_{\mathrm{fwd}} &\geq
        \frac{TC_T}{\e_{\mathrm{stab}}^2\min\{\sigma_x^2,\sigma_u^2\beta_T(0)\}},
        \qquad
    \lambda_{\mathrm{reg}} \geq \frac{C_T}{\e_0},
        \label{eq:end_to_end_weights}
\end{align}
then any minimizer $(\bar{E}_T, \bar{D}_T, \bar{L}_T, \bar{H}_T)$ of
$\calL_{\mathrm{bwd}}^{(T)}$ (cf.~\Cref{eq:bwd_solution})
satisfies the following conditions:
\begin{enumerate}[label=(\roman*)]
    \item The pair $(\bar{L}_T, \bar{H}_T)$ is stabilizable.
    \item A minimum $H_\infty$-norm controller
    \begin{align*}
        \bar{K}_T \in \argmin_{K \in B_{p \times m}^{(\mathrm{op})}(M_0)}
        \sup_{z \in \mathbb{D}^c}
        \opnorm{(zI_m-\bar{L}_T-\bar{H}_TK)^{-1}}
    \end{align*}
    exists, and every minimizer makes 
    $A+B\bar{K}_T\bar{E}_T$ Schur stable
    with the closed-loop $H_\infty$-norm bound:
    \begin{align*}
        \sup_{z \in \mathbb{D}^c}
        \opnorm{(zI_n-A-B\bar{K}_T\bar{E}_T)^{-1}}
        \leq 2\Gamma_0\opnorm{\Pi_{+,0}}
        \left[1+\frac{2\sigma_1}{\sqrt{\sigma_0}}
            \left(\sigma_0^{-1/2}+\Gamma_0\opnorm{B}M_0\right)\right].
    \end{align*}
\end{enumerate}
\end{mycorollary}
\begin{proof}
First, the assumption \eqref{eq:end_to_end_horizon} 
and definition of $\Phi$ (cf.~\Cref{assmp:S_T_decay})
implies:
\begin{align*}
    \frac{C^2(A,\gamma)\sigma_x^2(n-m)}{(T+1)(1-\gamma^2)}
        \leq\frac{\sigma_x^2\eta_0\e_{\mathrm{stab}}^2}{2}, 
        \qquad \frac{C^2(A,\gamma)\sigma_u^2\norm{B}_F^2S_T(\gamma)}{T+1}
    \leq C^2(A,\gamma)\sigma_u^2\norm{B}_F^2
        \frac{S_T(\gamma)}{T}
    \leq\frac{\sigma_x^2\eta_0\e_{\mathrm{stab}}^2}{2}.
\end{align*}
Adding these two inequalities together yields $C_T\leq\sigma_x^2\eta_0\e_{\mathrm{stab}}^2$.
Now we apply \Cref{thm:multistep_backwards_full}
using $\eta_{T+1}\geq\eta_0$ and the assumption \eqref{eq:end_to_end_weights} to conclude
(since the operator norm is bounded by the Frobenius norm):
\begin{align*}
    \max\left\{
        \opnorm{P_{E_0^\T}P^\perp_{\bar{E}_T^\T}},
        \opnorm{\bar{E}_TA-\bar{L}_T\bar{E}_T},
        \opnorm{\bar{E}_TB-\bar{H}_T}
    \right\} \leq\e_{\mathrm{stab}}, \qquad
    \sigma_0 I_m\preccurlyeq
        \bar{E}_T\bar{E}_T^\T&\preccurlyeq\sigma_1 I_m.
\end{align*}
We now invoke \Cref{thm:H_infty_closed_loop_stability}
for the pair $(\bar{E}_T, \bar{L}_T, \bar{H}_T)$.
The bounds above show that the requisite
\eqref{eq:H_infty_closed_loop_stability:epsilon_condition} holds, 
and therefore the result follows from the conclusions of \Cref{thm:H_infty_closed_loop_stability}, noting that 
$\opnorm{\bar{E}_T} \leq \sqrt{\sigma_1}$.
\end{proof}

\subsection{Analysis of Finite-Sample Losses}
\label{sec:proofs:finite_sample}

We first describe a simple Gaussian approximate isometry result
which our analysis crucially builds on.
\begin{myprop}[Gaussian approximate isometry]
\label{prop:gaussian_isometry}
Let $w_1, \dots, w_q \sim \sfN(0, I_d)$ be iid draws.
Fix $\delta \in (0, 1)$.
As long as $q \gtrsim d + \log(2/\delta)$, we have that with probability at least $1-\delta$,
for every $d_1 \in \N_+$ and $M \in \R^{d_1 \times d}$, 
\begin{align}
    \frac{1}{2} \norm{M}_F^2 \leq \frac{1}{q} \sum_{i=1}^{q} \norm{ M w_i }^2 \leq \frac{3}{2} \norm{M}_F^2. \label{eq:gaussian_isometry}
\end{align}
\end{myprop}
\begin{proof}
From standard results on the concentration of singular 
values for random Gaussian matrices~\cite{vershynin2018high}, we have that if $q \gtrsim d + \log(2/\delta)$, 
then with probability at least $1-\delta$,
\begin{align*}
    \frac{1}{2} I_d \preccurlyeq \frac{1}{q} \sum_{i=1}^{q} w_iw_i^\T \preccurlyeq \frac{3}{2} I_d.
\end{align*}
The result then follows on this event since for any $M \in \R^{d_1 \times d}$,
\begin{align*}
    \frac{1}{q} \sum_{i=1}^{q} \norm{ M w_i }^2 = \frac{1}{q} \sum_{i=1}^{q} w_i^\T M^\T M w_i = \Tr\left( M \cdot \frac{1}{q} \sum_{i=1}^{q} w_iw_i^\T \cdot M^\T \right).
\end{align*}
\end{proof}

We now restate and prove \Cref{prop:population_and_finite_sample_equivalence}.
\populationandfinitesampleequivalence*
\begin{proof}
We prove the result for the forward loss in \eqref{eq:fwd_loss_finite_sample};
the backward loss in \eqref{eq:bwd_loss_finite_sample} follows from a nearly identical argument.

We observe that for any $t \in [T]$, the function
$(x_0, \bm{u}_{0:t-1}) \mapsto E \Phi_{f}^t( x_0, \bm{u}_{0:t-1}) - \Phi_{f_z}^t( E x_0, \bm{u}_{0:t-1} )$
is a \emph{linear function} of the input, and
consequently of the \emph{whitened} vector $\bm{w}_{0:T} := (\bar{x}_0, \bm{\bar{u}}_{0:T-1})$,
where
$\bar{x}_0 := \frac{1}{\sigma_x} x_0$ and
$\bar{u}_t := \frac{1}{\sigma_u \sqrt{\alpha_{T}(t)}} u_t$.
Let $M_t(E, L, H) \in \R^{m \times (n + Tp)}$ denote the matrix such that:
\begin{align*}
    M_t(E, L, H) \bm{w}_{0:T} = E \Phi_{f}^t( x_0, \bm{u}_{0:t-1}) - \Phi_{f_z}^t( E x_0, \bm{u}_{0:t-1} ).
\end{align*}
Noting that $\mathrm{Law}( \bm{w}_{0:T} ) = \sfN(0, I_{n + Tp})$
and supposing that $N \gtrsim n + Tp + \log(2/\delta)$,
from \Cref{prop:gaussian_isometry} there is an event $\calE$ with probability at least $1-\delta$
such that on $\calE$,
\begin{align*}
    \wh{J}_{\mathrm{fwd}}^{(T)}(E, L, H) = \frac{1}{T}\sum_{t=1}^{T} \frac{1}{N}\sum_{i=1}^{N} \norm{ M_t(E, L, H) \bm{w}_{0:T}^{(i)} }^2 \leq \frac{3}{2T} \sum_{t=1}^{T} \norm{ M_t(E, L, H) }_F^2 = \frac{3}{2} J_{\mathrm{fwd}}^{(T)}(E, L, H),
\end{align*}
and similarly on $\calE$ we have $\wh{J}_{\mathrm{fwd}}^{(T)}(E, L, H) \geq \frac{1}{2} J_{\mathrm{fwd}}^{(T)}(E, L, H)$.
\end{proof}

\section{Insufficiency Results for Non Multi-Step Backward Losses}
\label{sec:proofs:insufficiency_results}

We first state and reprove
\Cref{thm:fwd_cycle_insufficient}.
\fwdcycleinsufficient*
\begin{proof}
We construct a multi-set $\Lambda \subset \eig(A)$ by the following procedure:

\begin{enumerate}[label=(\arabic*)]
    \item Initialize $\Lambda \gets \{ \}$ and $\Sigma \gets \eig(A)$.
    
    \item In case (a), choose a real $\lambda_\star\in\sigma_+(A)$ and remove
    one occurrence of it from $\Sigma$; in case (b), choose a nonreal
    $\lambda_\star\in\sigma_+(A)$ and remove 
    one occurrence of $\lambda_\star$ and $\bar{\lambda}_\star$ each from $\Sigma$.

    \item If $m$ is odd, move one real eigenvalue from $\Sigma$ to $\Lambda$.
    
    \item While $\abs{\Lambda}<m$: move a nonreal conjugate pair from $\Sigma$
    to $\Lambda$ if one is available, otherwise, move any two real
    eigenvalues from $\Sigma$ to $\Lambda$.
\end{enumerate}
For this procedure, we note that
step (2) is valid by the assumption that
either case (a) or case (b) holds.
Step (3) is valid since 
in case (a) we assume if $m$ is odd
then $m_{\mathrm{stab}} \geq 2$, 
and in case (b) we assume if $m$ is odd then
$m_{\mathrm{stab}} \geq 1$.
Step (4) is valid since 
at the beginning of the step (i.e., before the first iteration), we start with the three loop
invariants
(i) $m - \abs{\Lambda}$ is an even number, 
(ii) $\abs{\Sigma} \geq m - \abs{\Lambda}$, and
(iii) both $\Lambda, \Sigma$ are conjugate invariant,
and we maintain these loop invariants throughout as we always move two eigenvalues from $\Sigma$ to $\Lambda$ (either two real or a complex conjugate pair).
The procedure terminates with a conjugate-invariant multi-set $\Lambda \subset \eig(A)$
satisfying
$\abs{\Lambda}=m$ and
$\mathrm{mult}_{\Lambda}(\lambda_\star)<\mathrm{alg}_A(\lambda_\star)$.

With $\Lambda$ in place, we 
apply the real Schur decomposition
$(Q, S) = \mathrm{RealSchur}(A^\T, \Lambda)$
from \eqref{eq:real_schur_decomposition}
to obtain:
\begin{align*}
    A^\T=Q\bmattwo{S_{11}}{S_{12}}{0}{S_{22}}Q^\T,
    \qquad Q=\rvectwo{Q_1}{Q_2}\in O(n),
    \qquad \eig(S_{11})=\Lambda, \qquad \eig(S_{22}) = \eig(A) \setminus \Lambda,
\end{align*}
where $Q_1\in\R^{n\times m}$.
Define
\begin{align*}
    \bar{E}:=Q_1^\T, \qquad
    \bar{L}:=S_{11}^\T, \qquad
    \bar{H}:=\bar{E}B.
\end{align*}
The first block column of the Schur decomposition gives
$A^\T Q_1=Q_1S_{11}$, so
$\bar{E}A=\bar{L}\bar{E}$.

For (i), iterating these identities gives
$\bar{E}A^k=\bar{L}^k\bar{E}$ and
$\bar{E}A^kB=\bar{L}^k\bar{H}$ for every $k\geq0$.
Thus \Cref{prop:LTI_system_losses} yields
$J_{\mathrm{fwd}}^{(T)}(\bar{E},\bar{L},\bar{H})=0$ for every $T$.
Furthermore, since $\bar{E}\bar{E}^\T=I_m$,
regularizer admissibility (cf.~\Cref{def:forward_admissible_regularizer})
yields $R(\bar{E})=0$.
Hence the tuple $(\bar{E}, \bar{L}, \bar{H})$ 
satisfies $\calL^{(T)}_{\mathrm{fwd}}(\bar{E}, \bar{L}, \bar{H}) = 0$ for every
$T\in\N_+$ and $\lambda_{\mathrm{reg}}\geq0$, and is a global minimizer
because the loss is non-negative.

For (ii), suppose $(\bar{L},\bar{H})$ is not stabilizable.
The PBH test implies there exists $w\in\C^m\setminus\{0\}$ and $\mu\in\C$
with $\abs{\mu}\geq1$ such that
$w^*\bar{L}=\mu w^*$ and $w^*\bar{H}=0$.
The row vector $v^*:=w^*\bar{E}$ is nonzero, since
$\bar{E}\bar{E}^\T=I_m$, and satisfies
\begin{align*}
    v^*A=w^*\bar{L}\bar{E}=\mu v^*, \qquad
    v^*B=w^*\bar{H}=0.
\end{align*}
This contradicts stabilizability of $(A,B)$ by the PBH test.

For (iii), for every $K\in\R^{p\times m}$,
\begin{align*}
    Q^\T(A+BK\bar{E})Q
        =\bmattwo{\bar{L}+\bar{H}K}{0}
            {S_{12}^\T+Q_2^\T BK}{S_{22}^\T}.
\end{align*}
This matrix is block lower triangular, so its spectrum contains
$\eig(S_{22})=\eig(A)\setminus\Lambda$.
In particular, it contains $\lambda_\star$ (the excluded non-stable eigenvalue 
in the construction of $\Lambda$)
independently of $K$.
Since $\abs{\lambda_\star}\geq1$, 
$A + BK\bar{E}$ is not Schur stable.
\end{proof}

We next state and reprove
\Cref{thm:bwd_single_step_nonnormal}.
\bwdsinglestepnonnormal*
\begin{proof}
The proof strategy works as follows. We
first construct a candidate pair $(E', D', L', H')$
such that its one-step backward loss
is $O(\eta)$, and hence every minimizer
must also have loss $O(\eta)$.
Then, we show that any tuple
$(E, L, H)$ which minimizes the one-step backward
loss while simultaneously 
ensuring that
$\calK_L(L, H) \cap \calK_{A_\eta}(E)$ is non-empty
must then satisfy an inequality of the form
$\Omega(\eta^2) \leq O(\eta)$.
Thus, for all $\eta$ exceeding a critical cross-over point $\eta_0$,  
no minimizer of the one-step backward loss can
have a non-empty gain intersection, proving the result.

Write $\alpha:=\alpha_1(0)$ and $\tau:=2+\abs{\lambda_u+\lambda_s}$, and assume $\eta \geq 1$.
We define our candidate pair
$(E', D', L', H')$ as:
\begin{align*}
    E'&:=\frac{[1,\sqrt\eta]}{\sqrt{1+\eta}}, \qquad
    D':=\sqrt{1+\eta}\cvectwo{1}{0}, \qquad
    L':=\sqrt\eta, \qquad H':=E' B.
\end{align*}
Since $E'(E')^\T=1$, admissibility gives $R(E')=0$.
Direct calculation yields
\begin{align*}
    \norm{I_2-D' E'}_F^2&=1+\eta, \\
    \norm{A_\eta-L' D' E'}_F^2
        &=(\lambda_u-\sqrt\eta)^2+\lambda_s^2, \\
    \norm{B-D' H'}^2&=b_2^2(1+\eta), \\
    \norm{E' A_\eta-L' E'}^2
        &=\frac{(\lambda_u-\sqrt\eta)^2+\lambda_s^2\eta}{1+\eta},
        \qquad E' B-H'=0.
\end{align*}
Applying the identities from \Cref{prop:LTI_system_losses} with
$\beta_1(0)=\alpha_1(0)=\alpha$ gives:
\begin{align*}
    \calL_{\mathrm{bwd}}^{(1)}(E',D',L',H')
    &=\frac{\sigma_x^2}{2}\left(\norm{I_2-D'E'}_F^2
        +\norm{A_\eta-D'L'E'}_F^2\right)
        +\frac{\sigma_u^2\alpha}{2}\norm{B-D'H'}^2 \\
    &\qquad+\lambda_{\mathrm{fwd}}\left(
        \sigma_x^2\norm{E'A_\eta-L'E'}^2
        +\sigma_u^2\alpha\norm{E'B-H'}^2\right)
        +\lambda_{\mathrm{reg}}R(E') \\
    &=\frac{\sigma_x^2}{2}\left(1+\eta+(\lambda_u-\sqrt\eta)^2+\lambda_s^2\right)
        +\frac{\sigma_u^2\alpha}{2}b_2^2(1+\eta) + \lambda_{\mathrm{fwd}}\sigma_x^2
        \frac{(\lambda_u-\sqrt\eta)^2+\lambda_s^2\eta}{1+\eta} \\
    &\leq\frac{\sigma_x^2}{2}(1+2\eta+\lambda_u^2+\lambda_s^2)
        +\frac{\sigma_u^2\alpha}{2}b_2^2(1+\eta) + \lambda_{\mathrm{fwd}}\sigma_x^2
        \max\{\lambda_u^2,1+\lambda_s^2\}
    \leq C\eta.
\end{align*}
Here we used $(\lambda_u-\sqrt\eta)^2\leq\lambda_u^2+\eta$ and $\eta\geq1$.

Now let $(E,D,L,H)$ be any minimizer of $\calL^{(1)}_{\mathrm{bwd}}$ and put
\begin{align*}
    r_0:=\norm{I_2-DE}_F, \qquad
    r_A:=\norm{A_\eta-LDE}_F, \qquad
    r_B:=\norm{B-DH}.
\end{align*}
Optimality and nonnegativity of the loss terms imply
\begin{align}
    \frac{\sigma_x^2}{2}(r_0^2+r_A^2)
        +\frac{\sigma_u^2\alpha}{2}r_B^2
        \leq \calL_{\mathrm{bwd}}^{(1)}(E, D, L, H) \leq C\eta. \label{eq:nonnormal_residual_budget}
\end{align}
Now suppose there exists a $K \in \calK_L(L, H) \cap \calK_{A_\eta}(E)$.
Then $\abs{L+HK}<1$ and
\begin{align*}
    \abs{\Tr(A_\eta+BKE)} = \abs{ \lambda_1(A_\eta + BKE) + \lambda_2(A_\eta + BKE) } \leq \abs{\lambda_1(A_\eta+BKE)} + \abs{\lambda_2(A_\eta+BKE)} < 2.
\end{align*}
Therefore $\abs{KEB}<2+\abs{\lambda_u+\lambda_s}=\tau$.
Because $L,H,K$ and $EB$ are scalars, we have the exact identity
\begin{align*}
    A_\eta B=(A_\eta-LDE)B+(L+HK)DEB-(KEB)DH.
\end{align*}
Using $\norm{DEB}\leq(1+r_0)\norm{B}$ and
$\norm{DH}\leq\norm{B}+r_B$, Cauchy--Schwarz, and \eqref{eq:nonnormal_residual_budget} yields:
\begin{align*}
    \norm{ A_\eta B } &\leq \opnorm{ A_\eta - L D E } \norm{B} + \abs{L + HK} \norm{ D E B } + \abs{KEB} \norm{DH} \\
    &\leq r_A \norm{B} + (1+r_0)\norm{B} + \tau (\norm{B} + r_B) \\
    &= (1+\tau)\norm{B}+\norm{B}(r_0+r_A)+\tau r_B \\
    &\leq (1+\tau)\norm{B}
        +\left(\frac{\sigma_x^2}{2}(r_0^2+r_A^2)
        +\frac{\sigma_u^2\alpha}{2}r_B^2 \right)^{1/2}
            \left(\frac{4\norm{B}^2}{\sigma_x^2}
                +\frac{2\tau^2}{\sigma_u^2\alpha}\right)^{1/2} \\
    &\leq (1+\tau)\norm{B}
        +\sqrt{2C\eta}
            \left(\frac{2\norm{B}^2}{\sigma_x^2}
                +\frac{\tau^2}{\sigma_u^2\alpha}\right)^{1/2}.
\end{align*}
On the other hand, $\lambda_u > 1$, $\eta \geq 1$, $b_1 \geq 0$, and $b_2>0$ imply:
\begin{align*}
    \norm{A_\eta B}^2
        =\bignorm{\cvectwo{\lambda_u b_1+\eta b_2}{\lambda_s b_2}}^2 = (\lambda_u b_1 + \eta b_2)^2 + (\lambda_s b_2)^2 
        \geq (\eta b_2)^2 \implies \norm{A_\eta B} \geq \eta b_2.
\end{align*}
For $\eta\geq\eta_0$ with $\eta_0$ from \eqref{eq:eta_0},
each of the two terms in the preceding upper bound is strictly less
than $\eta b_2/2$, a contradiction.
\end{proof}

\subsection{Insufficiency of Reduced-Order Balanced Truncations}
\label{sec:appendix:balanced_truncation}

Suppose we have a controllable system $(A, B)$ of the form:
\begin{align}
    A = \bmattwo{A_u}{0}{0}{A_s}, \qquad B = \cvectwo{B_u}{B_s}, \quad B_u B_s^\T = 0, \label{eq:balanced_truncation_system}
\end{align}
where $A_u\in\R^{m_+\times m_+}$ contains all non-stable modes and $A_s$ is Schur stable. We will consider a balanced truncation where we explicitly keep all non-stable modes, 
and apply a balanced truncation to the stable sub-system $(A_s, B_s)$.
The additional condition $B_uB_s^\T=0$ ensures that the finite-horizon controllability Gramian is block-diagonal,
The finite-horizon observability Gramian is also block-diagonal because $A$ is block-diagonal and the measurement map is the identity.
Thus balanced truncation can be performed separately on the two blocks;
as the horizon tends to infinity, 
finite-horizon balancing of the Gramians for the entire unstable system as studied in \citet{flinois2015projection} will tend to this
separate balancing procedure.

Define the corresponding controllability and
observability Gramians (recall that the measurement map is $y_t = x_t$ in this paper) for $(A_s, B_s)$:
\begin{align*}
    W_c := \sum_{k=0}^{\infty} A_s^k B_sB_s^\T (A_s^k)^\T, \qquad W_o := \sum_{k=0}^{\infty} (A_s^k)^\T A_s^k.
\end{align*}
Let $n_s$ denote the dimension of $A_s$.
Since $(A,B)$ is controllable and $A$ is block-diagonal, the pair $(A_s,B_s)$ is also controllable, so that $W_c$ is positive definite; $W_o$ is positive definite because the stable state is fully observed.
To construct balanced coordinates, take the eigendecomposition
\begin{align*}
    W_c^{1/2} W_o W_c^{1/2} = U\Sigma^2 U^\T,
    \qquad \Sigma = \operatorname{diag}(\sigma_1,\ldots,\sigma_{n_s}),
    \qquad \sigma_1 \geq \cdots \geq \sigma_{n_s} > 0,
\end{align*}
where $U \in O(n)$ is orthogonal and the $\sigma_i$ are the Hankel singular values of $(A_s,B_s,I_{n_s})$.
Define the invertible coordinate transformation
\begin{align*}
    T_s := \Sigma^{1/2}U^\T W_c^{-1/2},
    \qquad T_s^{-1} = W_c^{1/2}U\Sigma^{-1/2}.
\end{align*}
Then the controllability and observability Gramians in the coordinates $\xi_s=T_sx_s$ coincide:
\begin{align*}
    T_sW_cT_s^\T = T_s^{-\T}W_oT_s^{-1} = \Sigma.
\end{align*}
To retain $r$ stable coordinates, partition the balanced realization as
\begin{align*}
    T_s A_s T_s^{-1} = \bmattwo{A_{11}}{A_{12}}{A_{21}}{A_{22}},
    \qquad T_sB_s = \cvectwo{B_1}{B_2},
    \qquad T_s^{-1} = \begin{bmatrix} C_1 & C_2 \end{bmatrix},
\end{align*}
where $A_{11}\in\R^{r\times r}$. Balanced truncation discards the last $n_s-r$ coordinates, yielding
\begin{align*}
    z_{s,t+1} = A_{11}z_{s,t}+B_1u_t,
    \qquad \hat{x}_{s,t}=C_1z_{s,t}.
\end{align*}
For a total latent dimension $m\in\{m_+,\ldots,n\}$, we take $r=m-m_+$, define $P_r:=\begin{bmatrix}I_r&0\end{bmatrix}\in\R^{r\times n_s}$, and append the unchanged non-stable subsystem:
\begin{align*}
    E=\bmattwo{I_{m_+}}{0}{0}{P_r T_s},
    \qquad L=\bmattwo{A_u}{0}{0}{A_{11}},
    \qquad H=\cvectwo{B_u}{B_1}.
\end{align*}

\paragraph{A destabilizing closed-loop example.}
We apply this balanced truncation
procedure to the system:
\begin{align*}
    A = \begin{bmatrix} 6/5 & 0 & 0 \\ 
        0 & 0 & 1/2 \\
        0 & 1/2 & 0 \end{bmatrix}, \qquad
    B = \begin{bmatrix} 1 & 0 \\ 0 & 1 \\ 0 & 0 \end{bmatrix}.
\end{align*}
This pair is controllable and satisfies $B_uB_s^\T=0$.
Its eigenvalues are $6/5,1/2,-1/2$, so $m_+=1$.
We choose $m=2$, and hence retain $r=1$ stable coordinate.
The stable subsystem and its Gramians are
\begin{align*}
    A_s=\bmattwo{0}{1/2}{1/2}{0},
    \qquad B_s=\bmattwo{0}{1}{0}{0},
    \qquad W_c=\bmattwo{16/15}{0}{0}{4/15},
    \qquad W_o=\frac43 I_2.
\end{align*}
Writing $\varphi:=(5/4)^{1/4}$, the balancing equations give
\begin{align*}
    U=I_2,\qquad
    \Sigma=\bmattwo{8/(3\sqrt5)}{0}{0}{4/(3\sqrt5)},
    \qquad T_s=\bmattwo{\varphi}{0}{0}{\sqrt2\,\varphi}.
\end{align*}
Retaining the first balanced stable coordinate therefore yields
\begin{align*}
    E=\begin{bmatrix}1&0&0\\0&\varphi&0\end{bmatrix},
    \qquad L=\bmattwo{6/5}{0}{0}{0},
    \qquad H=\bmattwo{1}{0}{0}{\varphi}.
\end{align*}
The gain
\begin{align*}
    K=\bmattwo{-6/5}{0}{0}{9/(10\varphi)}
    \qquad\implies\qquad
    L+HK=\bmattwo{0}{0}{0}{9/10},
\end{align*}
so it stabilizes the reduced model. However,
\begin{align*}
    A+BKE=\begin{bmatrix}0&0&0\\0&9/10&1/2\\0&1/2&0\end{bmatrix},
    \qquad
    \sigma(A+BKE)=\left\{0,\frac{9+\sqrt{181}}{20},\frac{9-\sqrt{181}}{20}\right\}.
\end{align*}
Since $(9+\sqrt{181})/20 \approx 1.1227 >1$, the full closed-loop is unstable.
This feedback stabilizes the originally unstable coordinate but destabilizes the originally stable two-dimensional sub-system.
Observe that the encoder has full row rank and $(A,E)$ is detectable, but the invariance condition fails:
\begin{align*}
    EA-LE=\begin{bmatrix}0&0&0\\0&0&\varphi/2\end{bmatrix}\ne0.
\end{align*}
The example shows that a stabilizing reduced-model gain need not stabilize the original system, even when all non-stable modes are retained exactly,
illustrating the importance of 
the approach taken in
\Cref{sec:problem_formulation:invariance_subspaces}.

\section{Minimum Eigenvalue of Non-Stable Gramians}

In this section we provide a self-contained proof
of the following result, which establishes a uniform
lower bound over time of the Gramian matrix for a
non-stable system.

\begin{mylemma}
\label{lemma:nonstable_gramian}
Let $M \in \R^{d \times d}$ be any matrix such that 
$\abs{\lambda} \geq 1$ for every $\lambda \in \sigma(M)$. 
Then:
\begin{align*}
    \inf_{T \in \N_+} \lambda_{\min}\left( \frac{1}{T} \sum_{t=0}^{T-1} (M^\T)^t M^t \right) > 0.
\end{align*}
\end{mylemma}
\begin{proof}
Let $M = S J S^{-1}$ denote the Jordan decomposition of $M$.
We have $M^k = S J^k S^{-1}$, and for a vector $v$,
\begin{align*}
    \norm{ M^k v }^2 = v^* S^{-*} (J^k)^* S^* S J^k S^{-1} v \geq \sigma_{\min}^2(S) \norm{ J^k S^{-1} v }^2,
\end{align*}
and hence 
\begin{align*}
    \lambda_{\min}\left( \frac{1}{T} \sum_{k=0}^{T-1} (M^\T)^k M^k \right) &= \lambda_{\min}\left( \frac{1}{T} \sum_{k=0}^{T-1} S^{-*} (J^k)^* S^* S J^k S^{-1} \right) \\
    &\geq \sigma_{\min}^2(S) \lambda_{\min}\left( S^{-*} \left( \frac{1}{T} \sum_{k=0}^{T-1} (J^k)^* J^k \right) S^{-1} \right) \\
    &\geq \frac{1}{\mathrm{cond}^2(S)} \lambda_{\min}\left( \frac{1}{T} \sum_{k=0}^{T-1} (J^k)^* J^k \right) = \frac{1}{\mathrm{cond}^2(S)} \min_{i =1, \dots, b} \lambda_{\min}\left( \frac{1}{T}\sum_{k=0}^{T-1} (J_i^k)^* J_i^k \right),
\end{align*}
where $J = \mathrm{blkdiag}(J_1, \dots, J_b)$.
Hence, our problem reduces to studying Jordan blocks $J_d(\lambda)$
associated with eigenvalues $\lambda$ satisfying $\abs{\lambda} \geq 1$.
Let us first assume that $\abs{\lambda} > 1$.
Since $\rho(J_d(\lambda)^{-1}) < 1$, 
there exists a $t_\lambda$ such that for all $t \geq t_\lambda$,
$\opnorm{ J_d(\lambda)^{-t} } \leq 1$.
Hence for any unit $v \in \C^d$ and $t \geq t_\lambda \geq 1$,
\begin{align*}
    1 = \norm{v} = \norm{ J_d(\lambda)^{-t} J_d(\lambda)^t v} \leq \norm{ J_d(\lambda)^t v }.
\end{align*}
Thus for $T \geq 2 t_\lambda$,
\begin{align*}
    \frac{1}{T} \sum_{t=0}^{T-1} \norm{ J_d(\lambda)^t v }^2 \geq \frac{1}{T} \sum_{t = t_\lambda}^{T-1} \norm{ J_d(\lambda)^t v }^2 \geq \frac{T-t_\lambda}{T} \geq \frac{1}{2}.
\end{align*}
On the other hand, when $T < 2 t_\lambda$,
\begin{align*}
    \frac{1}{T} \sum_{t=0}^{T-1} \norm{ J_d(\lambda)^t v }^2 \geq \frac{1}{T} > \frac{1}{2t_\lambda},
\end{align*}
and hence $\lambda_{\min}\left( \frac{1}{T}\sum_{k=0}^{T-1} (J_d(\lambda)^k)^* J_d(\lambda)^k \right) \geq \min\{1/2, 1/(2t_\lambda)\}$.
Thus the challenge is when $\abs{\lambda}=1$.
Suppose this is the case, and let us write
$J_d(\lambda) = \lambda I_d + N$ where $N_{i,i+1}=1$ is the nilpotent
matrix satisfying $N^d = 0$. Using the Binomial theorem we have:
\begin{align*}
    \lambda^{-t} J_d(\lambda)^t = \sum_{k=0}^{d-1} {t \choose k} \lambda^{-k} N^k. 
\end{align*}
We now use \Cref{prop:polynomial_lower_bound} to lower bound:
\begin{align*}
    \frac{1}{T} \sum_{t=0}^{T-1}\norm{ J_d(\lambda)^t v }^2 &= \frac{1}{T} \sum_{t=0}^{T-1}\bignorm{ \sum_{k=0}^{d-1} {t \choose k} \lambda^{-k} N^k v }^2 = \sum_{i=1}^{d} \left(\frac{1}{T} \sum_{t=0}^{T-1}\bigabs{ \sum_{k=0}^{d-1} {t \choose k} \lambda^{-k} (N^k v)_i  }^2\right) \\
    &\geq c_d \sum_{i=1}^{d} \abs{v_i}^2 = c_d ,
\end{align*}
from which we conclude 
$\lambda_{\min}\left( \frac{1}{T}\sum_{k=0}^{T-1} (J_d(\lambda)^k)^* J_d(\lambda)^k \right) \geq c_d$.
\end{proof}

The following auxiliary proposition was
utilized in the proof of \Cref{lemma:nonstable_gramian}.
\begin{myprop}
\label{prop:polynomial_lower_bound}
Fix a degree $d \in \N$. There exists a constant $c_d > 0$,
only depending on $d$, such that for all $T \in \N_+$ and
degree at most $d$ polynomials $p$,
\begin{align*}
    \frac{1}{T}\sum_{t=0}^{T-1} \abs{p(t)}^2 \geq c_d \abs{p(0)}^2.
\end{align*}
\end{myprop}
\begin{proof}
Let $\phi(t) := (1, t, t^2, \dots, t^d) \in \C^{d+1}$
and construct $\alpha \in \C^{d+1}$ such 
that $q(s) := p(T s) = \alpha^* \phi(s)$.
Therefore,
\begin{align*}
    \frac{1}{T}\sum_{t=0}^{T-1} \abs{p(t)}^2 = \frac{1}{T}\sum_{t=0}^{T-1} \abs{q(t/T)}^2 = \alpha^* \left[ \frac{1}{T}\sum_{t=0}^{T-1} \phi(t/T) \phi(t/T)^* \right] \alpha =: \alpha^* H_T \alpha.
\end{align*}
We have that $H_T$ converges to the limit of the Riemann sum:
\begin{align*}
    H_T \to H_\infty := \int_0^1 \phi(s) \phi(s)^* \rmd s.
\end{align*}
Now for any $\alpha \neq 0$ we have 
\begin{align*}
    \alpha^* H_\infty \alpha = \int_0^1 \bigabs{\sum_{k=0}^{d} \alpha_k s^k}^2 \rmd s > 0,
\end{align*}
since any non-zero polynomial cannot be equal to zero almost everywhere on $[0, 1]$. This means that $\lambda_{\min}(H_\infty) > 0$, and hence 
there is a $T_d \geq 1$ such that for $T \geq T_d$, we have
$\lambda_{\min}(H_T) \geq \lambda_{\min}(H_\infty)/2$.
Hence we can now proceed as follows.
For $T < T_d$,
\begin{align*}
    \frac{1}{T} \sum_{t=0}^{T-1} \abs{p(t)}^2 \geq \frac{1}{T} \abs{p(0)}^2 \geq \frac{1}{T_d} \abs{ p(0) }^2.
\end{align*}
Now when $T \geq T_d$,
\begin{align*}
    \frac{1}{T} \sum_{t=0}^{T-1} \abs{p(t)}^2 = \alpha^* H_T \alpha \geq \lambda_{\min}(H_\infty)/2 \cdot \norm{\alpha}^2 \geq \lambda_{\min}(H_\infty)/2 \cdot \abs{\alpha_0}^2 = \lambda_{\min}(H_\infty)/2 \cdot \abs{p(0)}^2.
\end{align*}
Hence we can take $c_d = \min\{1/T_d, \lambda_{\min}(H_\infty)/2\}$.
\end{proof}

\section{Experimental Details}
\label{sec:appendix:experiments}

\makeatletter\providecommand{\tabinput}[1]{\input{#1} }\makeatother
\newcommand{\exptab}[5]{\IfFileExists{#1}{\begin{table}[ht]\centering\footnotesize\begin{tabular}{#2}\toprule #3\\\midrule\input{#1} \bottomrule\end{tabular}\caption{#4}\label{#5}\end{table}}{}}

We provide tables that include further experimental details and results.
\Cref{tab:exp_protocols} summarizes the protocols of \Cref{sec:experiments}.
For the linear models of \Cref{sec:experiments:linear}, \Cref{tab:exp_gains} in the main text reports the number of initializations (population losses) and data sets (finite-sample losses)
whose LQR and $H_\infty$ gains stabilize the linearized arm, along with the subspace distance $\Delta$.
At $T = 40$, \Cref{tab:exp_conjugacy} complements it for both losses, in the population case and for $N \in \{100, 500, 1000\}$ trajectories,
by reporting the smallest eigenvalue $\lambda_{\min}(EE^\T)$, the forward conjugacy $r_A$ of the latent model in an orthonormal encoder basis, the conjugacy $r_u$ and $r_h$ of the unstable mode,
and the fraction of initial states from which the $H_\infty$ gain stabilizes the non-linear arm at each tip displacement.
For the models of \Cref{sec:experiments:nonlinear}, \Cref{tab:exp_conjugacy_nonlinear} reports, for data scales $\varepsilon \in \{0.01, 0.03, 0.1\}$ and horizons $T \in \{10, 20, 40\}$,
the deviation of the training data from linearized rollouts and the fraction of trajectories in which the arm falls,
along with the same quantities for the linear models and the MLP linearizations at the equilibrium, and the number of data sets whose latent LQR gain stabilizes the linearized arm.
Finally, \Cref{tab:exp_eps_tips} reports, for the same data scales and horizons, the fraction of initial states from which the linear and MLP controllers stabilize the non-linear arm at each tip displacement.

\paragraph{Unstable eigenvalue and input gain of a learned model.}
\label{sec:appendix:unstable_mode}
Panels (b) and (c) of \Cref{fig:exp_arm_nonlinear} compare the linearization $(E, L, H) = (J_E(0), J_{a_z}(0), g_z(0))$ of a learned model at the equilibrium with the canonical variables of \Cref{def:canonical_variables_unstable}.
Since $m_+ = 1$, these are $L_0 = \lambda_u := 1.129$, a left and a right eigenvector $E_0$ and $D_0$ of $A$ for $\lambda_u$ with $E_0 D_0 = 1$, and $H_0 = E_0 B$.
In order to compare the unstable components of $(E, L, H)$ to $(E_0, L_0, H_0)$ in a way that
both (a) handles the dimension mismatch and (b) is invariant to 
latent coordinate transformations, we consider the following procedure.
We let $\hat{\lambda}_u$ denote the eigenvalue of $L$ nearest to $\lambda_u$, and
measure the relative error $\abs{ \hat{\lambda}_u - \lambda_u } / \lambda_u$ as a way of
assessing the closeness of the unstable component of $L$. 
Now for $H$, we let $w$ denote
the left eigenvector for $\hat{\lambda}_u$, i.e., 
$w^\T L = \hat{\lambda}_u w^\T$, and we assume that it is possible to scale $w$ so that $w^\T E D_0 = 1$.
We observe that if 
$EA=LE$ held exactly, then we would have:
\begin{align*}
    w^\T EA = w^\T L E = \hat{\lambda}_u w^\T E \implies \lambda_u = \lambda_u w^\T E D_0 = w^\T E A D_0 = \hat{\lambda}_u w^\T E D_0 = \hat{\lambda}_u.
\end{align*} 
That is, we would have $w^\T E$ is a left eigenvector of $A$ with eigenvalue $\lambda_u$.
Now, since the row-space of $E_0$ spans the left unstable eigenspace of $A$, then we must have $w^\T E = \alpha E_0$ for some $\alpha$.
Multiplying both sides both $D_0$ yields
$1 = w^\T E D_0 = \alpha E_0 D_0 = \alpha$, and hence $w^\T E = E_0$.
Therefore, if we also have $EB = H$:
\begin{align*}
    w^\T H = w^\T E B =  E_0 B = H_0.
\end{align*}
Hence, we use the relative error $\abs{ w^\T H - H_0 } / \abs{H_0}$ as a measure for closeness of the unstable component of $H$.

\begin{table}[ht]
\centering\footnotesize
\setlength{\tabcolsep}{5pt}
\setlength{\aboverulesep}{0pt}\setlength{\belowrulesep}{0pt}%
\setlength{\extrarowheight}{3pt}%
\newcommand{\pb}{\rule[-4pt]{0pt}{0pt}}%
\newcolumntype{I}{!{\color{black!30}\vrule width 0.4pt}}%
\newcommand{\hl}{\noalign{{\color{black!30}\hrule height 0.4pt}}}%
\begin{tabularx}{\linewidth}{@{}>{\raggedright}p{1.9cm}<{\pb}
    I>{\raggedright\hsize=1.1\hsize}X<{\pb}
    I>{\raggedright\hsize=0.85\hsize}X<{\pb}
    I>{\raggedright\hsize=1.1\hsize}X<{\pb}
    I>{\raggedright\arraybackslash\hsize=0.95\hsize}X<{\pb}@{}}
\toprule
 & \multicolumn{2}{cI}{\Cref{fig:exp_arm_recovery}} & \multicolumn{1}{cI}{\multirow{2}{*}{\Cref{fig:exp_arm_nonlinear}}} & \multicolumn{1}{c@{}}{\multirow{2}{*}{\Cref{fig:exp_multi_unstable}, \Cref{tab:exp_multi_unstable}}} \\
\cmidrule(lr){2-3}
 & \multicolumn{1}{cI}{population} & \multicolumn{1}{cI}{finite sample} & \multicolumn{1}{cI}{} & \multicolumn{1}{c@{}}{} \\
\midrule
plant & linearized arm & linearized arm & non-linear arm, coordinates $x/\varepsilon$ & three-mode linearized arm, mass coordinates \\ \hl
data & none (exact expectations) & $5$ data sets, $N \in \{5, 20, 100, 500, 1000\}$, nested & $5$ data sets, $N = 2000$, $\varepsilon = 0.01$ ($0.03$ and $0.1$ in \Cref{tab:exp_conjugacy_nonlinear}) & $5$ data sets, $N = 2000$ \\ \hl
initializations & $8$ & $2$ per data set and $N$ & MLP: $8$ per data set, one continued; linear: $8$ per data set & $16$ per data set and strategy, all counted \\ \hl
horizons & in increasing order, fresh and continued initializations & as population & as population & each $T$ separately \\ \hl
optimizer & full-batch Adam ($3000$ steps) + L-BFGS ($\le 2000$), double precision & as population & MLP: minibatch Adam ($30{,}000$ steps); linear: as population & full-batch Adam ($2000$) + L-BFGS ($2000$) \\ \hl
forward term in $\calL_{\mathrm{bwd}}^{(T)}$ & $J_{\mathrm{fwd}}^{(1)}$ & $J_{\mathrm{fwd}}^{(1)}$ & MLP: teacher forcing; linear: $J_{\mathrm{fwd}}^{(1)}$ & $J_{\mathrm{fwd}}^{(1)}$ \\ \hl
$\lambda_{\mathrm{reg}}$ & $\calL_{\mathrm{bwd}}^{(T)}$: $1$; $\calL_{\mathrm{fwd}}^{(T)}$: most favorable & as population & $1$ & $1$ \\ \hl
reported & lowest objective (line), every initialization (dots) & lowest objective per data set & lowest objective per data set & every run \\ \hl
stability test & $\rho(A + BKE) < 1$; non-linear arm: the same $8$ initial states scaled to tip displacements from $2$ to $80$\,cm & as population & $\rho(A + BKJ_E(0)) < 1$, $K$ from $(J_{a_z}(0), g_z(0))$; non-linear arm: the same $8$ initial states scaled to tip displacements from $2$ to $80$\,cm & -- \\
\bottomrule
\end{tabularx}
\caption{Protocols of the experiments of \Cref{sec:experiments}. The data are regenerated at every horizon with the same initial states.}
\label{tab:exp_protocols}
\end{table}

{\setlength{\tabcolsep}{5pt}%
\exptab{tables/linear_conjugacy.tex}{lllccccccccccc}
{loss & regularizer & $N$ & $\lambda_{\min}$ & $r_A$ & $r_u$ & $r_h$ & \multicolumn{7}{c}{$H_\infty$, non-linear arm, tip [cm]} \\ \cmidrule(l){8-14} & & & & & & & $2$ & $10$ & $20$ & $30$ & $40$ & $60$ & $80$}
{Linear models at $T = 40$ ($\lambda_{\mathrm{reg}}$ as in \Cref{fig:exp_arm_recovery}): for the population losses (pop.), the lowest-objective initialization;
for the empirical losses, the reported models of the $5$ data sets with $N$ trajectories ($\Delta$ and stabilization of $A + BKE$: \Cref{tab:exp_gains}).
$\lambda_{\min}$ is the smallest eigenvalue of $EE^\T$, and $r_A := \norm{\hat E A - \hat L \hat E}_F$, where $(\hat E, \hat L, \hat H) := (S^{-1}E, S^{-1}LS, S^{-1}H)$ with $S := (EE^\T)^{1/2}$ is the same model with orthonormal encoder rows.
$r_u := \norm{e^\T A - \hat{\lambda}_u e^\T}$ and $r_h := \abs{e^\T B - w^\T H}$, where $\hat{\lambda}_u$ is the eigenvalue of $L$ nearest $\lambda_u = 1.129$, $w$ is its left eigenvector, and $e := E^\T w$, scaled so that $\norm{e} = 1$;
both residuals are invariant under any invertible change of latent coordinates.
We report $r_u$ and $r_h$ for models with $\Delta < 0.5$.
For the empirical losses, $\lambda_{\min}$, $r_A$, $r_u$ and $r_h$ are medians over the data sets.
Last columns: fraction of the $8$ initial states (averaged over the data sets) from which the $H_\infty$ gain stabilizes the non-linear arm, per tip displacement.%
\label{tab:exp_gauge}\label{tab:exp_finite}}%
{tab:exp_conjugacy}}

\exptab{tables/nonlinear_conjugacy.tex}{llrrlccccccc}
{$\varepsilon$ & $T$ & dev.\ [\%] & fallen [\%] & model & captured & $\Delta$ & $\lambda_{\min}$ & $r_A$ & $r_u$ & $r_h$ & stable}
{Non-linear arm data at scale $\varepsilon$, and the linearization $(E, L, H) := (J_E(0), J_{a_z}(0), g_z(0))$ of the models selected at horizon $T$
(for the linear model, its parameters). Data: deviation from linearized rollouts, and fraction of trajectories in which a section passes $90^\circ$.
Models: data sets (of $5$) that capture the unstable mode, medians over the data sets of $\Delta$, $\lambda_{\min}$ and $r_A$,
medians over the capturing data sets of $r_u$ and $r_h$ (all as in \Cref{tab:exp_conjugacy}), and data sets whose latent LQR gain stabilizes $A + BKE$
(stabilization of the non-linear arm: \Cref{tab:exp_eps_tips}).%
\label{tab:exp_eps}}%
{tab:exp_conjugacy_nonlinear}

\exptab{tables/nonlinear_eps_tips.tex}{lllccccccc}
{$\varepsilon$ & $T$ & controller & \multicolumn{7}{c}{tip displacement [cm]} \\ \cmidrule(l){4-10} & & & $2$ & $10$ & $20$ & $30$ & $40$ & $60$ & $80$}
{Non-linear arm data at scale $\varepsilon$ (cf.~\Cref{tab:exp_conjugacy_nonlinear}): fraction of the $8$ initial states, averaged over the $5$ data sets,
from which the controllers of the models selected at horizon $T$ stabilize the non-linear arm, per tip displacement
(panel (e) of \Cref{fig:exp_arm_nonlinear} shows $\varepsilon = 0.01$, $T = 40$ per data set).}
{tab:exp_eps_tips}

\end{document}